\documentclass[12pt]{article}

\usepackage[utf8]{inputenc}
\usepackage[T2A]{fontenc}

\usepackage[english]{babel}
\usepackage{amsfonts,amsmath, amssymb}
\usepackage{xcolor}

\usepackage{enumerate}
\usepackage{epsf,epsfig,amsfonts,a4wide}
\usepackage{amsfonts}
\usepackage{amsmath,amssymb,amsthm}
\usepackage{url}

\usepackage{amssymb}
\usepackage{amsthm}
\usepackage{epsf,epsfig,amsfonts,graphicx}
\usepackage{a4wide}
\newtheorem{Theorem}{Theorem}[section]
\newtheorem{Lemma}[Theorem]{Lemma}
\newtheorem{Proposition}[Theorem]{Proposition}
\theoremstyle{definition}
\newtheorem{Ex}[Theorem]{Example}
\newtheorem{Definition}[Theorem]{Definition}
\newtheorem{RemarkDefinition}[Theorem]{Remark and Definition}
\newtheorem{Comment}[Theorem]{Comment}
\newtheorem{Remark}[Theorem]{Remark}

\newcommand{\R}{\mathbb{R}}
\newcommand{\Id}{\textrm{\rm Id}}
\newcommand{\re}[0]{\mathrm{Re \,}}
\newcommand{\im}[0]{\mathrm{Im \,}}
\newcommand{\ddd}{\mathrm{d}}
\newcommand{\const}{\mathrm{const}}

\newcommand{\rank}{\operatorname{rank}}
\newcommand{\Symp}{\operatorname{Symp}}
\newcommand{\spec}{\operatorname{Spec}}
\newcommand{\Sing}{\operatorname{Sing}}
\newcommand{\Span}{\operatorname{Span}}

\newcommand{\Lie}[0]{\mathcal {L}}
\def\I{\mathfrak{I}}
\def\H{\mathfrak{H}}
\def\M{\mathfrak{M}}

\title{Hidden toric symmetry and structural stability \\ of singularities in integrable systems}
\author{Elena A. Kudryavtseva
\thanks{Moscow State University, Moscow, Russia;
Moscow Center for Fundamental and Applied Mathematics, Moscow, Russia.
E-mail address: eakudr@mech.math.msu.su. ORCID: 0000-0002-3745-2663}
}
\date{}

\begin{document}
\maketitle

\begin{abstract}
The goal of the paper is to develop a systematic approach to the study of (perhaps degenerate) singularities of integrable systems and their structural stability. 
Following N.T.~Zung, as the main tool, we use  ``hidden''  system-preserving torus actions near singular orbits.
We give sufficient conditions for the existence of such actions and show that they are persistent under integrable perturbations. We find toric symmetries for several infinite series of singularities and prove, as an application, structural stability of Kalashnikov's parabolic orbits with resonances in the real-analytic case. We also classify all Hamiltonian $k$-torus actions near a singular orbit on a symplectic manifold $M^{2n}$ (or on its complexification) and prove that the normal forms of these actions are persistent under small perturbations.
As a by-product, we prove an equivariant version of the Vey theorem (1978) about local symplectic normal form of nondegenerate singularities.


\end{abstract}

\section{Introduction} \label {sec:1}

Let $(M^{2n},\Omega,f_1,\dots,f_n)$ be an integrable Hamiltonian system with $n$ degrees of freedom, where the momentum map $F=(f_1,\dots,f_n)$ is proper.
Consider the Hamiltonian $\R^n$-action on $M$ generated by the momentum map $F$.
We will call orbits of this action simply {\em orbits}.
Consider the singular fibration (called the {\em Liouville fibration}), whose fibers are connected components of the level sets $F^{-1}(a)$, $a\in\R^n$.

By a {\em local} (respectively {\em semilocal}) {\em singularity} of such a singular fibration we will mean the fibration germ at a singular orbit (respectively fiber).
We recall that a point $m_0\in M$ is called a {\em singular} (or {\em critical}) point of this fibration if $\rank \ddd F(m_0)<n$. An orbit (or a fiber) is called {\em singular} if it contains a singular point.
The minimal rank of singular points belonging to a fiber is called {\em rank} of the fiber.

Topology and geometry of singular Lagrangian fibrations have been studied from local \cite{vey, var:giv82, dui, ler87, eli, ito91, zung96, zung96a, bau:zung97, zung00, zung03, zung:mir04, zung04, bgk}, semilocal \cite {fom, ler:uma2, zung96, zung96a, BF, ler00, kud:lep11, kud12} and global viewpoints \cite{fom, BF, bol:fom:ric, zung04}. 
By using Hamiltonian torus actions (and more general group actions) near singularities, N.T.~Zung developed a systematic approach for study and classification of singularities of singular Lagrangian fibrations associated with integrable systems \cite{zung96, zung96a, bau:zung97, zung00, zung:habil, zung03, zung:mir04, zung04}.
In particular, Zung applied such actions for a classification of local and semilocal nondegenerate singularities in the works \cite{zung96,zung96a} where the notions of action-angle coordinates and a system-preserving Hamiltonian torus action (introduced in the Liouville-Arnold-Mineur theorem near regular fibres) were extended to the case when singularities are present (see also \cite[Proposition 4]{fom}, \cite[Theorem 3.2]{BF}).

Concerning local singularities, it is known that a compact rank-$r$ orbit of an $n$-degrees of freedom integrable system always admits a locally-free $F$-preserving Hamiltonian $(S^1)^r$-action on some neighbourhood of this orbit, provided that the orbit is either regular (Liouville theorem), or singular and has one of the following types:
\begin{itemize}
\item[--] nondegenerate (Ito \cite{ito91} for the real-analytic case, Zung \cite {zung96, zung96a} for the general case, Fomenko \cite[Proposition 4]{fom}, \cite[Theorem 3.2]{BF} for the twisted hyperbolic case with $n=2$ and $r=1$),
\item[--] not too degenerate (Bau and Zung \cite[Theorem 2.1]{bau:zung97}),
\item[--] having finite type (Zung \cite {zung03}, \cite [Theorem 3.7]{zung04}, cf.\ \S\ref{subsec:adj}) in real-analytic case,
\item[--] $r=n-1$ and $\dim_\R L_0\le n$ (Zung \cite[Theorem 1.2]{zung00}) where $L_0$ is the fiber containing the given orbit in real-analytic case,
\item[--] parabolic (the author and Martynchuk \cite{kud:mar21}) in the smooth $C^\infty$ case.
\end{itemize}
Furthermore, if the singular orbit is nondegenerate with Williamson type $(k_e,k_h,k_f)$, then this action extends to an effective (not locally-free) $F$-preserving Hamiltonian $(S^1)^{r+k_e+k_f}$-action (Zung \cite[Theorem 6.1]{zung96a}), moreover the system is fiberwise symplectomorphic to a linear model (\cite{ito91} for the real-analytic case, \cite{zung:mir04} for the $C^\infty$ case). 

For degenerate local singularities, only partial results were obtained about Hamiltonian torus actions (and more general group actions)  
\cite{bau:zung97, zung00, zung03, zung04}, classification \cite{dui, bro93, han} and structural stability \cite{ler87, kal, gia, vdm, var:giv82, gar04, osh-tuz, kud:mar21} of these singularities. For singularities of integrable systems with incomplete flows, only partial results were obtained \cite{kud:lep11, kud12}. 
Rigidity of some degenerate corank 1 singular points (i.e.\ in the ``very local'' setting), including parabolic ones, was proved under some assumptions only (Givental and Varchenko \cite {var:giv82}, Garay \cite {gar04}).

The purpose of this paper is to study the Hamiltonian torus actions near local singularities in more detail
and to apply them for studying and classifying degenerate singularities
of singular Lagrangian fibrations associated with integrable systems (including systems with incomplete flows), and for studying structural stability of such singularities.
More specifically, 
three natural circles of questions are addressed in this paper. The first circle of questions is as follows:
\begin{itemize}
\item[(1a)] Does there exist an effective Hamiltonian $(S^1)^{r+\varkappa_e}$-action on a neighbourhood of a {\em degenerate} orbit ${\cal O}$, that preserves the momentum map, 
where the $(S^1)^{r}$-subaction is locally-free, and the $(S^1)^{\varkappa_e}$-subaction leaves the orbit ${\cal O}$ fixed?
\item[(1b)] Can this torus action be extended to an effective Hamiltonian $(S^1)^{r+\varkappa_e+\varkappa_h}$-action on a small open complexification $U^{\mathbb C}$ of a {\em degenerate} (respectively, {\em nondegenerate}) orbit ${\cal O}$, that preserves the holomorphic extension $F^{\mathbb C}$ of the momentum map to $U^{\mathbb C}$, where the $(S^1)^{\varkappa_h}$-subaction leaves the orbit ${\cal O}$ fixed?
\item[(1c)] Is this torus action persistent under small real-analytic integrable perturbations of the system?
\end{itemize}
The second circle of questions is about symplectic normalization of a torus action:
\begin{itemize}
\item[(2a)] Can a Hamiltonian torus action be written in a ``simple'' normal form (so-called canonical model) w.r.t.\ some symplectic coordinates on a neighbourhood of an orbit ${\cal O}$?
\item[(2b)] Can a (more general) Hamiltonian torus action be written in a ``simple'' symplectic normal form, when some of $S^1$-subactions fixing the orbit ${\cal O}$ are allowed to be generated by real-analytic functions multiplied with $\sqrt{-1}$?
\item[(2c)]
Is this normalization of the torus action persistent (rigid) under small real-analytic integrable perturbations of the torus action?
\end{itemize}
The third circle of questions is on symmetry and structural stability of singularities:
\begin{itemize}
\item[(3a)] 
For a given singularity of an integrable system, describe its hidden toric symmetries, in particular compute all its resonances (elliptic, hyperbolic and twisting) introduced in this paper.
\item[(3b)] 
For a given (degenerate) local singularity, how to prove its structural stability under small integrable perturbations (in real-analytic case, Definition \ref {def:stab})?
\item[(3c)] 
Can the singular fibration and its toric symmetry be written in a ``simple'' form
(called preliminary normal form) 
w.r.t.\ some coordinates in a neighbourhood of the singularity, and what form has the symplectic structure in these coordinates?
\item[(3d)]
Is this preliminary normal form persistent under small integrable perturbations?
\end{itemize}

In this paper, we solve the above questions (1a)--(2c) affirmatively, under certain weak conditions (Theorems \ref {thm:kudr}, \ref {thm:kudr:} and \ref {thm:period}, \ref {thm:period:}). In particular, if the orbit is compact and nondegenerate, then $r+\varkappa_e+\varkappa_h=n$. For illustration, we answer the question (3a) for several infinite series of local singularities (Examples \ref {exa:res}, \ref {exa:res:}), and conjecture structural stability of all these singularities (Example \ref {exa:stab} (B)). As an application, we solve the questions (3b)--(3d) for parabolic orbits with resonances, which are degenerate local singularities with $r=n-1$, $\varkappa_e=\varkappa_h=0$ (Proposition \ref {pro:kal}). 
As a by-product, we prove an equivariant version of the Vey theorem \cite{vey} about symplectic local normal form for nondegenerate singularities (Lemmata \ref {lem:period} and \ref {lem:period:}).
Our proofs are analogous to the proofs of theorems about torus actions in \cite{zung96, zung96a, bau:zung97, BF, zung00, zung03}.

When using the term ``hidden'' for toric symmetries, we mean the following:
\begin{itemize}
\item This symmetry is a Hamiltonian $(S^1)^{r+\varkappa_e+\varkappa_h}$-action generated by some smooth functions (``actions'') depending on the first integrals of the system, but these ``action functions'' are often not given or not known in advance.
\item The $(S^1)^{r+\varkappa_e}$-subaction is defined only on a small neighbourhood $U$ of the given singular orbit, so this subaction does not necessarily extend globally to the whole phase space.
\item The $(S^1)^{\varkappa_h}$-subaction is generated by imaginary-valued ``action functions'', so this subaction is defined only on a small open complexification $U^{\mathbb C}$ of the neighbourhood $U$ of the singular orbit.
\end{itemize}
The above properties of being ``hidden'' for the $(S^1)^{r+\varkappa_e+\varkappa_h}$-action and its subactions are observed in many integrable mechanical systems, e.g.\ for regular orbits (so-called Liouville tori, $r=n$, $\varkappa_e=\varkappa_h=0$), nondegenerate singular orbits ($r<n$, $r+\varkappa_e+\varkappa_h=n$) and parabolic orbits \cite{kud:mar21} ($r=k_e=1$, $\varkappa_e=\varkappa_h=0$, $n=2$).

Remarks.

1) In the case of a locally-free $(S^1)^r$-action, the question (1a) was solved 
for nondegenerate singularities by Ito \cite{ito91} (for real-analyltic case), Fomenko \cite[Proposition 4]{fom}, \cite[Theorem 3.2]{BF} (for the twisted hyperbolic case when $n=2$ and $r=1$), Zung \cite {zung96, zung96a}; 
for not too degenerate singularities by Bau and Zung \cite[Theorem 2.1]{bau:zung97}; 
for finite type singularities by Zung \cite {zung03}, \cite [Theorem 3.7]{zung04} (see \S\ref{subsec:adj});
for degenerate corank-1 analytic singularities by Zung \cite[Theorem 1.2]{zung00}; 
for smooth parabolic singularities by the author and Martynchuk \cite{kud:mar21}.

2) For nondegenerate singularities of integrable systems, the questions (1a), (2a), (3a) and (3c) w.r.t.\ the $(S^1)^{r+\varkappa_e}$-subaction were solved by Ito \cite{ito91} (for real-analytic case), Miranda and Zung \cite{zung:mir04} (for equivariant $C^\infty$ case).
For nondegenerate singularities, the questions (3b) and (3d) w.r.t.\ the $(S^1)^{r+\varkappa_e}$-subaction were partially solved by Miranda and Zung \cite{zung:mir04} (for equivariant $C^\infty$ case), for a weaker notion of structural stability (resp.\ persistence), namely for structural stability (resp.\ persistence) under parametric families of integrable perturbations.

3) For (degenerate) rank-0 singularities, a solution to the questions (1a), and partially (1b), (2a), (2b) was described by Zung \cite {zung04} in terms of the Poincar\'e-Birkhoff normal form, by proving its convergence.

4) For an arbitrary smooth symplectic action of a compact Lie group on a neighbourhood of its fixed point, the question (2a) was partially solved by Weinstein \cite[Lecture 5]{wei}, \cite{del:mel} by linearizing the action in some symplectic coordinates.
For an arbitrary smooth symplectic action of a Lie group on a neighbourhood of its orbit, a solution to the questions (2a), (2b) was given 
by Marle (\cite{mar83}, \cite[Propositions 1.9 and 1.10]{mar}), Guillemin and Sternberg \cite{gui:ste} in terms of a linear model of the action. 
The question (2a) in a global setting was solved in some cases (cf.\ \cite{kar:tol14, kar:zil16} and references therein).

5) For some degenerate corank 1 singularities including parabolic ones, in the analytic case, the questions (3b)--(3d) were solved under some assumptions in the ``very local'' setting (in a neighbourhood of a singular point) by Givental and Varchenko \cite {var:giv82}, Garay \cite {gar04}. For parabolic orbits,
 % (without resonance), 
the question (3b) on structural stability was solved by Lerman and Umanskii \cite{ler87}, see \cite{kud:mar21} for smooth structural stability.
For parabolic orbits with resonances, the question (3b) was partially solved by Kalashnikov \cite {kal}, for a weaker notion of structural stability, namely for structural stability under $S^1$-symmetry-preserving integrable perturbations.
Infinitesimal stability (i.e.\ stability under infinitesimal integrable deformations of the system \cite[Definition 8]{gia}) was studied for 2-degrees of freedom integrable systems, namely: nondegenerate rank-0 and rank-1 singular points and a rank-1 parabolic singular point are infinitesimally stable \cite[Definition 9, Theorems 2 and 3]{gia}. This partially solves the question (3b) for these singularities, for a weaker notion of structural stability, namely for infinitesimal stability.

6) For a corank-2 singularity ``integrable Hamiltonian Hopf bifurcation'' of integrable Hamiltoninan systems with 3 degrees of freedom, the question (3c) was partially solved by van der Meer \cite{vdm} (without studying the symplectic structure).

7) For saddle-saddle fibers satisfying a ``non-splitting'' condition, structural stability under ``component-wise'' $C^\infty$ integrable perturbations was proved by Oshemkov and Tuzhilin \cite{osh-tuz}. This partially solves the question (3b) for such semilocal singularities, for a weaker notion of structural stability.

Our solutions to the questions (1a)--(3d) have the following advantages:
\begin{itemize}
\item Our solution to the questions (1a)--(1c) (in Theorems \ref {thm:kudr}, \ref {thm:kudr:}) 
can be applied to degenerate orbits of any corank, to 
torus actions that are not locally-free, and to integrable systems with incomplete flows.
\item 
Our solution to the questions (2a)--(2c) (in Theorems \ref {thm:period}, \ref {thm:period:}) can be applied to integrable systems with degenerate singularities; furthermore we show that our canonical model is not only linear but also has  a ``diagonal'' form (strengthening
 %in contrast to 
\cite {wei, del:mel, mar83, mar, gui:ste}).
\item 
In our solution to the question (3a) (in Examples \ref {exa:res}, \ref {exa:res:}), 
local singularities have resonances of different qualitative nature (so-called elliptic, hyperbolic and twisting resonances), and these resonances cannot be reduced or simplified.
In our solution to the questions (3b)--(3d) for parabolic orbits with resonances (in Proposition \ref {pro:kal}), we prove structural stability of a singularity (resp., persistence of a preliminary normal form) under arbitrary small integrable perturbations (Definition \ref {def:stab}). In particular, we do not assume that the $S^1$-action is preserved under the perturbation. We also do not assume that the perturbation is parametric, so our ``perturbed'' system is not necessarily included into a parametric family of integrable systems containing the ``unperturbed'' one. We also do not assume that the ``perturbed'' system has a singular orbit close to the ``unperturbed'' orbit.
\end{itemize}
We expect that our solutions to the questions (1a)--(2c), and (3a) for several infinite series of local singularities, as well as (3b)--(3d) for parabolic orbits with resonances, will be helpful for solving the questions (3b)--(3d) for other singularities, including those from Examples \ref {exa:res} and \ref {exa:res:}, as we conjectured in Example \ref {exa:stab} (B).

Our results on degenerate singularities and their structural stability can be eventually applied in Mechanics and Physics to study nearly integrable Hamiltonian systems. H.~Poincar\'e \cite{poi} considered such systems as the {\em basic problems of mechanics}. One of the most famous approaches to Poincar\'e's basic problem is the KAM (Kolmogorov-Arnold-Moser) theory.
The classical KAM theorem says that when the perturbation is small, most of the invariant tori of the unperturbed system will not be destroyed but only slightly deformed. This theorem is stated under the following nondegeneracy condition, also known as Kolmogorov's condition: $\det(\partial^2H_0/\partial I_i\partial I_j)\ne0$, where $H_0$ is the unperturbed Hamilton function,
and $(I_i)$ is a local system of action variables. In practice, this condition is not easy at all to verify directly, because the computation of the determinant often involves transcendental functions. However, one can often verify this condition easily using singularities (Kn\"orrer \cite{kno}, Bau and Zung \cite{bau:zung97}).

\section {A hidden toric symmetry: existence and persistence under small integrable perturbations} \label {subsec:period}

This section is devoted to solving the questions (1a)--(1c) from Introduction.

The following theorem tells us how to partially solve the questions (1a) and (1c) w.r.t.\ some of the $S^1$-subactions of the desired torus action. In detail: we should apply this theorem several times, in order to find several $S^1$-actions (some of them will be locally free, and the others will leave our orbit fixed). 
Such $S^1$-actions automatically pairwise commute, so all together they form a single $(S^1)^{r'+\varkappa_e'}$-action, that will be a subaction of the desired $(S^1)^{r+\varkappa_e}$-action, where $0\le r'\le r$, $0\le \varkappa_e'\le \varkappa_e$.
This theorem gives sufficient conditions for
\begin{itemize}
\item the existence of a (not necessarily locally-free) Hamiltonian 
$S^1$-action that preserves the momentum map, and
\item persistence of such an action under small integrable perturbations.
\end{itemize}
The case of a nondegenerate orbit (the existence part only) was treated in \cite {ito91,fom,BF,zung96,zung96a}.
The case of a locally-free action 
(again the existence part only) was treated in \cite[Theorem 2.1]{bau:zung97}, \cite[Theorem 1.2]{zung00} (the corank-1 case), \cite {zung03}.
Actually our proof is analogous to the proof of theorems about torus actions in \cite{zung96, zung96a, bau:zung97, BF, zung00, zung03}.

Denote by ${\cal O}_m$ the orbit of a point $m\in M$ under the (local) Hamiltonian action of $\R^n$ on $(M,\Omega)$ generated by the functions $f_1,\dots,f_n$.
Denote by $X_{f}$ the Hamiltonian vector field with the Hamilton function $f$.

In the following two theorems, the momentum map $F$ is not necessarily proper. 

\begin{Theorem} \label {thm:kudr}
Let $(M^{2n},\Omega,F)$ be a real-analytic integrable Hamiltonian system, $m_0\in M$ a singular point of the momentum map $F=(f_1,\dots,f_n)$, $F(m_0)=(0,\dots,0)$, and $L_0=F^{-1}(0,\dots,0)$ the singular fiber containing the point $m_0$.

Suppose there exists a point $m_1\in L_0$ satisfying the following conditions:
\begin{itemize}
\item[\rm(i)] $\rank \ddd F(m_1)=n$, i.e.\ $m_1$ is a regular point of the momentum map $F$,
\item[\rm(ii)] there exist a compact trajectory $\gamma_0$ (i.e.\ a closed trajectory or an equilibrium) and a continuous one-parameter family of $2\pi$-periodic trajectories $\gamma_u\subset {\cal O}_{m_1}$, $0<u\le 1$, of the vector field $X_{f_1}$ such that $m_1\in\gamma_1=:\gamma$ and $m_0\in\gamma_0\subset\overline{\bigcup\limits_{0<u\le 1}\gamma_u}=:C$.\footnote{In practice, one can verify the
condition (ii) of the theorems \ref {thm:kudr} and \ref {thm:kudr:} using Lemma \ref {lem:ii}.}
\end{itemize}
Then

{\rm (a)} There exist a neighbourhood $U$ of the set $C$ in $M$ (cf.\ {\rm(ii)}) and a unique $F$-preserving Hamiltonian $S^1$-action on $U$ generated by a function $I(f_1,\dots,f_n)$, where $I(z_1,\dots,z_n)$ is a real-analytic function on the neighbourhood $V=F(U)$ of the origin in $\R^n$ such that $I(z_1,\dots,z_n)=z_1+O(\sum\limits_{j=1}^n|z_j|^2)$.
The action function $I(z_1,\dots,z_n)$ can be computed by the Mineur-Arnold \cite {min} integral formula
\begin{equation} \label {eq:action:real:new}
I(z_1,\dots,z_n)=\frac1{2\pi}\oint\limits_{\gamma_{(z_1,\dots,z_n)}} \alpha+\const, \qquad \qquad (z_1,\dots,z_n)\in V.
\end{equation}
Here $\gamma_{(z_1,\dots,z_n)}\subset F^{-1}(z_1,\dots,z_n)$ denotes a closed curve depending continuousely on $(z_1,\dots,z_n)$ such that $\gamma_{(0,\dots,0)}=\gamma$, and $\alpha$ is any analytic 1-form on a neighbourhood of $\gamma$ such that $\Omega=\ddd\alpha$ (such a 1-form always exists, see (\ref {eq:alpha}) below).

{\rm (b)} This $S^1$-action is persistent under real-analytic integrable perturbations in the following sense.
Suppose we are given an integer $k\ge2$, a neighbourhood $U'\subset U$ of the set $C$ in $M$ (cf.\ {\rm(ii)}) and a neighbourhood $V'$ of the origin in $\R^n$ having compact closures $\overline{U'}\subset U$ and $\overline{V'}\subset V$.
Then there exists $\varepsilon>0$ such that, for any (``perturbed'') real-analytic integrable Hamiltonian system $(M^{2n},\tilde\Omega,\tilde F)$ that is $\varepsilon-$close to the initial system in $C^k$-norm, the following properties hold.
There exist a bigger neighbourhood $\tilde U\supset U'$ and a unique $\tilde F$-preserving Hamiltonian (w.r.t.\ the ``perturbed'' symplectic structure) $S^1$-action on $\tilde U$ generated by a function 
$\tilde I(\tilde f_1,\dots,\tilde f_n)$, where $\tilde I(z_1,\dots,z_n)$ is a real-analytic function on some bigger neighbourhood $\tilde V=\tilde F(\tilde U)\supset V'$ that is $O(\varepsilon)-$close to the function $I(z_1,\dots,z_n)$ in $C^k$-norm.
The ``perturbed'' action function $\tilde I(z_1,\dots,z_n)$ can be computed by the Mineur-Arnold integral formula
\begin{equation} \label {eq:action:real:new:}
\tilde I(z_1,\dots,z_n)=\frac1{2\pi}\oint\limits_{\tilde\gamma_{(z_1,\dots,z_n)}} \tilde\alpha+\const, \qquad \qquad (z_1,\dots,z_n)\in \tilde V,
\end{equation}
where $\tilde\gamma_{(z_1,\dots,z_n)}\subset \tilde F^{-1}(z_1,\dots,z_n)$ denotes a closed curve close to $\gamma$, and $\tilde\alpha$ is any analytic 1-form on a neighbourhood of $\gamma$ such that $\tilde\Omega=\ddd\tilde\alpha$ (such a 1-form always exists).
\end{Theorem}

In the following theorem, we show how to solve the questions (1b) and (1c) from Introduction, as well as the remaining part of the questions (1a) and (1c). Similarly to the previous theorem, we formulate our solution w.r.t.\ some of the $S^1$-subactions of the desired torus action. By using this theorem, one can obtain the remaining $(S^1)^{r-r'+\varkappa_e-\varkappa_e'+\varkappa_h}$-subaction of the desired $(S^1)^{r+\varkappa_e+\varkappa_h}$--action.

Denote by $M^{\mathbb C}$ a small open complexification of $M$, on which $\Omega^{\mathbb C}$ and $F^{\mathbb C}$ are defined.
Denote by ${\cal O}_m^{\mathbb C}$ the orbit of a point $m\in M^{\mathbb C}$ under the (local) Hamiltonian action of ${\mathbb C}^{n}$ on $M^{\mathbb C}$ generated by the functions $f_1^{\mathbb C},\dots,f_n^{\mathbb C}$.
For a holomorphic function $f$ on $M^{\mathbb C}$, denote by $X_{f}$ the Hamiltonian vector field on $M^{\mathbb C}$ with the Hamilton function $f$.
Denote by $\Sing(F^{\mathbb C})$ the set of singular points of $F^{\mathbb C}$.

\begin{Theorem} \label {thm:kudr:}
Let $(M^{2n},\Omega,F)$ be a real-analytic integrable Hamiltonian system, $m_0\in M$ a singular point of the momentum map $F=(f_1,\dots,f_n)$, $F(m_0)=(0,\dots,0)$.

Suppose there exists a point $m_1\in L_0^{\mathbb C}=(F^{\mathbb C})^{-1}(0,\dots,0)$ and $\lambda\in{\mathbb C}\setminus\{0\}$ such that
\begin{itemize}
\item[\rm(i)] $\rank \ddd F^{\mathbb C}(m_1)=n$, i.e.\ $m_1$ is a regular point of the map $F^{\mathbb C}$,
\item[\rm(ii)] there exist a compact trajectory $\gamma_0$ (i.e.\ a closed trajectory or an equilibrium) and a continuous one-parameter family of $2\pi$-periodic trajectories $\gamma_u\subset {\cal O}_{m_1}^{\mathbb C}$, $0<u\le 1$,
of the vector field $X_{\lambda f_1}$ such that $m_1\in\gamma_1=:\gamma$ and $m_0\in\gamma_0\subset\overline{\bigcup\limits_{0<u\le 1}\gamma_u}=:C$.$^1$
\end{itemize}
Then

{\rm (a)} There exist a neighbourhood $U$ of the set $C$ in $M^{\mathbb C}$ (cf.\ {\rm(ii)}) and a unique $F^{\mathbb C}$-preserving Hamiltonian $S^1$-action on $U$ generated by a function $I(f_1^{\mathbb C},\dots,f_n^{\mathbb C})$, where $I(z_1,\dots,z_n)$ is a holomorphic function on the neighbourhood $V=F^{\mathbb C}(U)$ of the origin in ${\mathbb C}^n$ such that $I(z_1,\dots,z_n)=\lambda z_1+O(\sum\limits_{j=1}^n|z_j|^2)$.
The action function $I(z_1,\dots,z_n)$ can be computed by the Mineur-Arnold \cite {min} integral formula
\begin{equation} \label {eq:action:real:C:new}
I(z_1,\dots,z_n)=\frac1{2\pi}\oint\limits_{\gamma_{(z_1,\dots,z_n)}} \alpha^{\mathbb C}+\const, \qquad \qquad (z_1,\dots,z_n)\in V,
\end{equation}
where $\gamma_{(z_1,\dots,z_n)}\subset (F^{\mathbb C})^{-1}(z_1,\dots,z_n)$ denotes a closed curve depending continuousely on $(z_1,\dots,z_n)$ such that $\gamma_{(0,\dots,0)}=\gamma$, and $\alpha^{\mathbb C}$ is any holomorphic 1-form on a neighbourhood of $\gamma$ such that $\Omega^{\mathbb C}=\ddd\alpha^{\mathbb C}$ (such a 1-form always exists, similarly to (\ref {eq:alpha})).

Let, in addition, the curve $\gamma$ in {\rm(ii)} be {\em homologically symmetric} in the following sense:
\begin{itemize}
\item[\rm(iii)] for each $\varepsilon>0$, there exists $a\in \R^n$, $|a|<\varepsilon$, such that the closed path $\gamma_a$ from {\rm(\ref {eq:action:real:C:new})} is homological in the fiber $(F^{\mathbb C})^{-1}(a)\setminus \Sing(F^{\mathbb C})$ to its $\mathbb C$-conjugated path $\overline{\gamma_a}$ (respectively, to the closed path obtained from $\overline{\gamma_a}$ by reversing orientation).
\end{itemize}
Then $\lambda\in \R$ (respectively $\lambda\in i\R$) and the ``normalized'' action function $\frac1\lambda I(z_1,\dots,z_n)$ is real-valued (and, hence, real-analytic) on the domain $V\cap\R^n$.

{\rm (b)} This $S^1$-action is persistent under real-analytic integrable perturbations in the following sense.
Suppose we are given $k\in{\mathbb Z}_+$, a neighbourhood $U'$ of the set $C$ in $M^{\mathbb C}$ and a neighbourhood $V'$ of the origin in ${\mathbb C}^n$ having compact closures $\overline{U'}\subset U$ and $\overline{V'}\subset V$.
Then there exists $\varepsilon>0$ such that, for any (``perturbed'') real-analytic integrable Hamiltonian system $(M,\tilde\Omega,\tilde F)$ whose holomorphic extension to $M^{\mathbb C}$ is $\varepsilon-$close to $(M^{\mathbb C},\Omega^{\mathbb C},F^{\mathbb C})$ in $C^0-$norm, the following properties hold. On some neighbourhood $\tilde U\supset U'$, there exists a unique $\tilde F^{\mathbb C}$-preserving Hamiltonian (w.r.t.\ the ``perturbed'' symplectic structure $\tilde\Omega$) $S^1$-action generated by a function $\tilde I(\tilde f_1^{\mathbb C},\dots,\tilde f_n^{\mathbb C})$, where $\tilde I(z_1,\dots,z_n)$ is a holomorphic function on some neighbourhood $\tilde V\supset V'$ that is $O(\varepsilon)-$close to $I(z_1,\dots,z_n)$ in $C^k-$norm. The ``perturbed'' action function $\tilde I(z_1,\dots,z_n)$ can be computed by the Mineur-Arnold integral formula
\begin{equation} \label {eq:action:real:C:new:}
\tilde I(z_1,\dots,z_n)=\frac1{2\pi}\oint\limits_{\tilde\gamma_{(z_1,\dots,z_n)}} \tilde\alpha^{\mathbb C}+\const, \qquad (z_1,\dots,z_n)\in \tilde V,
\end{equation}
where $\tilde\gamma_{(z_1,\dots,z_n)}\subset (\tilde F^{\mathbb C})^{-1}(z_1,\dots,z_n)$ is a closed curve close to $\gamma$,
$\tilde\alpha^{\mathbb C}$ is a holomorphic 1-form on a neighbourhood of $\gamma$ such that $\tilde\Omega^{\mathbb C}=\ddd\tilde\alpha^{\mathbb C}$ (such a 1-form always exists).

Let, in addition, the curve $\gamma$ in {\rm(ii)} be {\em homologically symmetric}, i.e.\ satisfy {\rm(iii)} from {\rm(a)}.
Then the ``perturbed'' ``normalized'' action function $\frac1\lambda \tilde I(z_1,\dots,z_n)$ is real-valued (and, hence, real-analytic) on the domain $\tilde V\cap\R^n$.
\end{Theorem}

Proofs of Theorems \ref {thm:kudr} and \ref {thm:kudr:} are given in \S \ref {sec:app}.

\subsection {Nondegenerate singularities, singular orbits of finite type} \label {subsec:adj}

A singular point $m_0$ of rank $0$ is called {\em nondegenerate} (cf. e.g. \cite{des90}) if the linearizations $A_j$ of the Hamiltonian vector fields $X_{f_j}$ at the point $m_0$ span a Cartan subalgebra of the Lie algebra of the Lie group $\Symp(T_{m_0}M,\Omega|_{m_0})\simeq\Symp(2n,\R)$, i.e. the operators $A_1,\dots,A_n$ span an $n$-dimensional commutative subalgebra and there exists a linear combination $A=\sum\limits_{j=1}^nc_jA_j$, $c_j\in\R$, having a simple spectrum: $|\spec A|=2n$. 
A singular point $m_0$ of rank $r$ is called {\em nondegenerate} (cf. e.g. \cite{des90}) if the corresponding rank-$0$ singular point of the corresponding reduced integrable Hamiltonian system with $n-r$ degrees of freedom (obtained by local symplectic reduction under the action of $f_1,\dots,f_r$ such that $df_1\wedge\dots\wedge df_r|_{m_0}\ne0$) is nondegenerate.

A singular orbit (respectively, fiber) is called {\em nondegenerate} if each singular point contained in this orbit (fiber) is nondegenerate.

A singular orbit ${\cal O}_{m_0}$ is called of {\em finite type} \cite[Definition 3.6]{zung04} if there is only a finite number of orbits of the infinitesimal action of ${\mathbb C}^n\approx\R^{2n}$ on the fiber $L_0^{\mathbb C}\subset M^{\mathbb C}$ containing $m_0$, and
 $L_0^{\mathbb C}$ contains a regular point of the map $F^{\mathbb C}$.

Due to the Vey theorem \cite{vey}, each nondegenerate orbit is of finite type.

\subsection {On topological conditions (ii), (iii) in Theorems \ref {thm:kudr}, \ref {thm:kudr:}} \label {subsec:ii:iii}

In practice, the periodicity condition (ii) in the theorems \ref {thm:kudr} and \ref {thm:kudr:} 
can be verified using Lemma \ref {lem:ii} below.

Suppose $m_0$ is a rank 0 singularity (the general rank $r$ case can be reduced to this special case by local symplectic reduction under the action of $f_1,\dots,f_r$ such that $df_1\wedge\dots\wedge df_r|_{m_0}\ne0$).
Consider the commuting vector fields $X_{f_1},\dots,X_{f_n}$ on the regular part $L_0\setminus \Sing(F)$ of the singular fiber.
Since they are linearly independent, they form a basis of $T_mL_0$ at each point $m\in L_0\setminus \Sing(F)$. Since they pairwise commute, they define a flat affine connection on $L_0\setminus \Sing(F)$, which
has a trivial holonomy. Roughly speaking, the periodicity condition (ii) from Theorem \ref {thm:kudr} simply means that 
there exists a closed geodesic $\gamma$ on $L_0\setminus \Sing(F)$ w.r.t.\ this 
flat affine connection, moreover 
this geodesic is ``close enough'' to the point $m_0$ (see Lemma below for a precise formulation).
Similarly, the periodicity condition (ii) from Theorem \ref {thm:kudr:} means that 
there exists a closed geodesic $\gamma$ on $L_0^{\mathbb C}\setminus \Sing(F^{\mathbb C})$ w.r.t.\ the similar flat affine connection on the complexified fiber $L_0^{\mathbb C}\setminus \Sing(F^{\mathbb C})$, moreover this geodesic is ``close enough'' to the point $m_0$.

Let us proceed with precise statements.
Denote by $\theta_1,\dots,\theta_n$ the 1-forms on $L_0\setminus \Sing(F)$ forming a dual basis of $T_m^*L_0$ to the basis $X_{f_1}|_m,\dots,X_{f_n}|_m$ of $T_mL_0$ at each point $m\in L_0\setminus \Sing(F)$, and by $\ddd s^2:=\sum\limits_{i=1}^n(\theta_i)^2$ the flat Riemannian metric on $L_0\setminus \Sing(F)$.

Let $B_{0,r}\subset\R^{2n}$ denote the open ball of radius $r$ centred at $0$.
Let $(U,\phi)$ be a canonical coordinate chart centred at $m_0$ with $\phi:U\to\phi(U)=B_{0,r}$, thus 
$\Omega=\phi^*\sum\limits_{i=1}^n\ddd x_{2i-1}\wedge\ddd x_{2i}$. 
Let $c>0$ be a real constant such that $|\ddd(f_i\circ\phi^{-1})(x)|\le c|x|$ for all $x\in B_{0,r}$, $1\le i\le n$.

\begin{Lemma} \label {lem:ii}
Suppose $m_0\in M$ is a rank $0$ singular point. Suppose there exists a closed curve $\hat\gamma$ (not necessarily a geodesic) in $L_0\setminus \Sing(F)$ such that   
$(a)$ 
$\oint\limits_{\hat\gamma}\theta_1\ne0$,
$(b)$ either the flows of all $X_{f_i}$ are complete in $L_0$ and $m_0\in \overline{{\mathcal O}_{m_1}}$ for $m_1\in\hat\gamma$, 
or there exists a continuous path $m_u\in(L_0\setminus \Sing(F))\cap\phi^{-1}(B_{0,re^{-nc\ell}})$, $0<u\le 1$, such that $m_1\in\hat\gamma$ and $m_0\in\overline{\{m_u\}_{0<u\le1}}$,
where $\ell:=\oint\limits_{\hat\gamma}\ddd s$. Then there exists a linear combination $f:=\sum\limits_{i=1}^n a_if_i$, $a_1\ne0$, such that the points $m_0,m_1$ and the momentum map $(f,f_2,\dots,f_n)$ satisfy the condition {\rm(ii)} from Theorem \ref {thm:kudr}.
\end{Lemma}

\begin{Remark} \label {rem:ii}  
The conditions $(a)$ and $(b)$ in Lemma \ref {lem:ii} can not be omitted, in general. To see necessity of $(b)$, consider the integrable system $({\mathbb C}^2,\re(\ddd z\wedge\ddd w),F)$ with the momentum map $F=(f_1,f_2)=(\re f(z,w),\im f(z,w))$ given by the hyperellyptic polynomial $f(z,w)=z^2+P_5(w)$, $P_5(w)=(w-a_1)(w-a_2)(w-a_3)w^2$, $a_1<a_2<a_3<0$. 
Topology and geometry of such Lagrangian fibrations are studied in \cite {kud:lep11, kud12}. The fibre $L_0=F^{-1}(0)$ is a punctured genus-2 surface with a pinch at the focus-focus point $m_0=(0,0)$ (Fig.~1) \cite[\S3--5]{kud:lep11}. Its regular part $L_0\setminus\{m_0\}$ consists of three families of periodic orbits and three unbounded trajectories $s_1=\{(\pm\sqrt{-P_5(w)},w)\mid w\in(-\infty,a_1]\}$, $s_2$ and $s_3$ of the vector field $X_{f_1}$ \cite[\S 6]{kud:lep11}. 
Orbits of two of the families approach the point $m_0$. The third family contains the periodic orbit $\gamma=\{(\pm\sqrt{-P_5(w)},w)\mid w\in[a_2,a_3]\}$.
Thus, the singular fibre $L_0$ in the chart $U={\mathbb C}^2$ contains the periodic orbit $\gamma$ that can not be connected to the singular point $m_0=(0,0)$ by a family of periodic orbits in this fibre.
To see necessity of $(a)$, consider a similar system with $P_5(w)$ replaced by $w^3$.
The circle $\hat\gamma=\{|v|=1\}$ in the fibre $L_0=\{(v^3,-v^2)\mid v\in{\mathbb C}\}$ satisfies $(b)$ of Lemma \ref {lem:ii} after compactifying fibres via attaching a point at infinity to each fibre \cite[Lemma 7]{kud:lep11}.
%being included into the family of concentric circles (as in the focus-focus case). 
But $\hat\gamma$ does not provide $S^1$-action near $m_0$,
since $\theta_1=\ddd\,\re(v^{-1})$, $\theta_2=-\ddd\,\im(v^{-1})$ are exact.
\end{Remark}

\begin{figure}[htbp] %\label {fig:genus:two}
\begin{center} 
\includegraphics[width=0.4\linewidth]{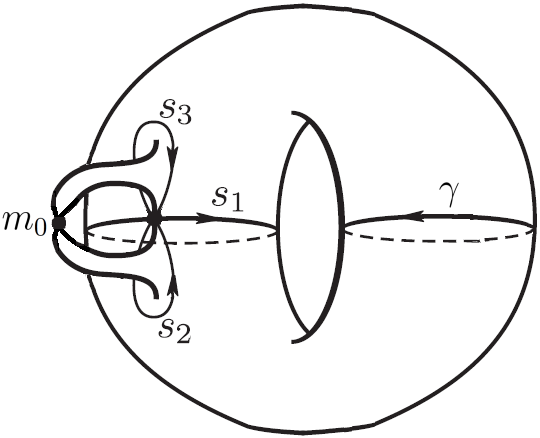} \qquad
\includegraphics[width=0.38\linewidth]{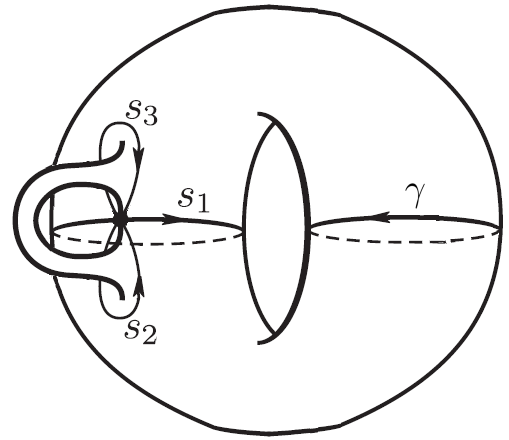}  
\caption{Trajectories of the vector field $X_{f_1}$ with incomplete flow 
on a singular fibre (left) and on a nearby regular fibre (right)}
\label{fig/standard_model}
\end{center}
\end{figure}

\begin{proof}  
Consider the loops in $L_0\setminus \Sing(F)$ based at $m_1$, which are homotopic to $\hat\gamma$ relatively basepoint. Let us show that this set of loops contains a shortest loop w.r.t.\ the Riemannian metric $\ddd s^2$. 
Since the loop $\hat\gamma$ has length $\ell$, it is sufficient to consider only piece-wise geodesic loops of lengths $\le\ell$.
All these loops lie in the image of the ball of radius $\ell/2$ under the exponential map at $m_1$ of the Riemannian metric $\ddd s^2$, provided that this map is well-defined on that ball.
Notice that any geodesic is a trajectory of the vector field $v_a=X_f$ for some constant vector $a=(a_1,\dots,a_n)\in\R^n$, $f:=\sum\limits_{i=1}^n a_if_i$, and the length of its velocity vector $v_a$ 
w.r.t.\ $\ddd s^2$ equals $|a|$. Denote by $\gamma_{m,a}=\gamma_{m,a}(t)$ the trajectory of $v_a$
starting at $m\in L_0\setminus \Sing(F)$. Consider two cases.

{\em Case 1:} the flows of all $X_{f_i}$ are complete in $L_0$. Then the Riemannian metric $\ddd s^2$ on $L_0\setminus \Sing(F)$ is geodesically complete, thus its exponential map is well-defined, which shows the existence of a shortest loop $\gamma$ in $L_0\setminus \Sing(F)$ 
based at $m_1$ homotopic to $\hat\gamma$ relatively basepoint.
Due to $(a)$, we have $\oint\limits_{\gamma}\theta_1=\oint\limits_{\hat\gamma}\theta_1\ne0$, thus $\gamma$ is not just a point.
By above, $\gamma$ is a $2\pi$-periodic trajectory $\gamma=\gamma_{m_1,a}$ of the vector field $v_a$ for some vector $a\in\R^n$, so $2\pi a_1=\oint\limits_{\gamma}\theta_1\ne0$. 
Therefore some neighbourhood of $\gamma$ in $L_0\setminus \Sing(F)$ is filled by $2\pi$-periodic trajectories of $v_a$.
Since $m_0\in \overline{{\mathcal O}_{m_1}}$, there exists a path $m_u\in {\mathcal O}_{m_1}$, $0<u\le 1$, 
such that $m_0\in\overline{\{m_u\}_{0<u\le 1}}$. 
By above arguments, the trajectories $\gamma_{m_u,a}$, $0< u\le 1$, of $v_a$ are $2\pi$-periodic and 
form a required family satisfying (ii).

{\em Case 2} (without completeness assumption). If $|a|=1$ and a trajectory $\gamma_{m,a}$ of the vector field $v_a=X_f$ lies in $U$, then 
$|\frac{\ddd}{\ddd t}\phi(\gamma_{m,a}(t))|
=|\ddd(f\circ\phi^{-1})(\phi(\gamma_{m,a}(t)))|
\le nc|\phi(\gamma_{m,a}(t))|$, since $\phi$ is a canonical chart. Hence,
\begin{equation} \label{eq:close}
\frac{\ddd}{\ddd t}|\phi(\gamma_{m,a}(t))|^2
= 2 \langle \phi(\gamma_{m,a}(t)), \frac{\ddd}{\ddd t}\phi(\gamma_{m,a}(t))\rangle
\le 2 nc |\phi(\gamma_{m,a}(t))|^2.
\end{equation}
Denote $r_1:=|\phi(m_1)|$.
Due to \eqref{eq:close}, each trajectory segment $\gamma_{m_1,a}|_{[0,\ell_1]}$ in $U$ with $|a|=1$ 
is contained in the closure of $\phi^{-1}(B_{0,r_1e^{nc\ell_1}})$. Due to $(b)$, we have $r_1<re^{-nc\ell}$, therefore the trajectory segment 
$\gamma_{m_1,a}|_{[0,\ell/2]}$ lies inside the ball $\phi^{-1}(B_{0,re^{-nc\ell/2}})$,
so it will not reach the boundary of $U$.
Thus the exponential map at $m_1$ is well-defined on the ball of radius $\ell/2$.
This shows the existence of a shortest loop $\gamma$ in $L_0\setminus \Sing(F)$ based at $m_1$ homotopic to $\hat\gamma$ relatively basepoint. The rest of the proof can be given similarly to Case 1 using 
the path $m_u$, $0<u\le 1$, from $(b)$.
\end{proof}

\begin{Ex} \label {exa:nondeg}
Let us verify the condition (iii) from Theorem \ref {thm:kudr:} (a) on homological symmetry for basic nondegenerate singularities: elliptic, hyperbolic and focus-focus.

(e) Consider an elliptic nondegenerate rank 0 singularity, given by $f_1=\frac12(p^2+q^2)$ and $\Omega=\ddd p\wedge \ddd q$ on $\R^2$ with coordinates $(p,q)$. Take a small $\varepsilon>0$ and a regular point $m_1=(\varepsilon,i\varepsilon)$. Since $f_1^{\mathbb C}(m_1)=0$, we have $m_1\in L_0^{\mathbb C}$. Then the Hamiltonian system has the form
$\frac{\ddd p}{\ddd t}=-\frac{\partial f_1}{\partial q}=-q$,
$\frac{\ddd q}{\ddd t}= \frac{\partial f_1}{\partial p}= p$.
Its solutions are $\gamma_{a,b}(t)=(ae^{-it}+be^{it},iae^{-it}-ibe^{it})$ with arbitrary constants $a,b\in\mathbb C$.
So we have a solution $\gamma(t)=\gamma_{\varepsilon,0}(t)=\varepsilon e^{-it}(1,i)$ with $\gamma(0)=:m_1\in L_0^{\mathbb C}$. Since this solution is $2\pi$-periodic, we have $\lambda=1\in\R$.
Take a real regular value $a=4\varepsilon^3\in\R$ close to $0$. Consider the closed path
$\gamma'(t)=\gamma_{\varepsilon,\varepsilon^2}(t)=\varepsilon e^{-it}(1+e^{2it}\varepsilon,i-ie^{2it}\varepsilon)$ in the (Milnor's) fiber $(F^{\mathbb C})^{-1}(a)=\{p^2+q^2=4\varepsilon^3\}$.
This path is a closed path obtained from $\gamma$ by a small deformation.
The $\mathbb C$-conjugated path is
$\overline{\gamma'(t)}=\varepsilon e^{it}(1+e^{-2it}\varepsilon,-i+ie^{-2it}\varepsilon)=\gamma_{\varepsilon^2,\varepsilon}(t)$.
Then $\gamma'(t)$ and $\overline{\gamma'(t)}$ are homological in the (Milnor's) fiber $(F^{\mathbb C})^{-1}(a)$.
Indeed, they are orbits of the Hamiltonian $S^1$-action generated by $f_1^{\mathbb C}$, thus they can be connected with each other by a 1-parameter family of such orbits in the regular (Milnor's) fiber $(F^{\mathbb C})^{-1}(a)$.
We have a real $\lambda=1\in\R$ and a real-analytic $2\pi$-periodic first integral $I(f_1)=f_1$, as Theorem \ref {thm:kudr:} (a) asserts.

(h) Consider a hyperbolic nondegenerate rank 0 singularity, given by $f_1=pq$ and $\Omega=\ddd p\wedge \ddd q$ on $\R^2$ with coordinates $(p,q)$.
Then the Hamiltonian system with the Hamilton function $if_1^{\mathbb C}$ has the form
$\frac{\ddd p}{\ddd t}=-i\frac{\partial f_1}{\partial q}=-ip$,
$\frac{\ddd q}{\ddd t}= i\frac{\partial f_1}{\partial p}= iq$.
Its solutions are $\gamma_{a,b}(t)=(ae^{-it},be^{it})$ with arbitrary constants $a,b\in\mathbb C$.
Take a small $\varepsilon>0$ and a regular point $m_1=(\varepsilon,0)\in L_0$.
So we have a solution $\gamma(t)=\gamma_{\varepsilon,0}(t)=(\varepsilon e^{-it},0)$ with $\gamma(0)=m_1\in L_0$. Since this solution is $2\pi$-periodic, we have $\lambda=i\in i\R$.
The $\mathbb C$-conjugated path is $\overline{\gamma(t)}=(\varepsilon e^{it},0)=\gamma(-t)$.
Thus $\gamma(t)$ and $\overline{\gamma(-t)}$ are homological in the regular part of the fiber $L_0^{\mathbb C}$, since they just coincide.
We have an imaginary $\lambda=i\in i\R$ and a holomorphic $2\pi$-periodic first integral $I(f_1)=if_1$, moreover $iI(f_1)=-f_1$ is real-valued,
as Theorem \ref {thm:kudr:} (a) asserts.

(f) Consider a focus-focus nondegenerate rank 0 singularity, given by $f_1=p_1q_2-p_2q_1$, $f_2=p_1q_1+p_2q_2$ and $\Omega=\ddd p_1\wedge \ddd q_1+\ddd p_2\wedge \ddd q_2$ on $\R^4$ with coordinates $(p,q)=(p_1,p_2,q_1,q_2)$. We have two commuting Hamiltonian $S^1$-actions on a small open complexification of the origin, namely those generated by $f_1^{\mathbb C}$ and $if_2^{\mathbb C}$.

The Hamiltonian $S^1$-action generated by $f_1^{\mathbb C}$ is given by the Hamiltonian system
$\frac{\ddd p_1}{\ddd t}=-\frac{\partial f_1}{\partial q_1}=p_2$,
$\frac{\ddd p_2}{\ddd t}=-\frac{\partial f_1}{\partial q_2}=-p_1$,
$\frac{\ddd q_1}{\ddd t}= \frac{\partial f_1}{\partial p_1}=q_2$,
$\frac{\ddd q_2}{\ddd t}= \frac{\partial f_1}{\partial p_2}=-q_1$.
Its trajectory $\gamma(t)=\varepsilon(\cos t,-\sin t,0,0)$ with $\gamma(0)=(\varepsilon,0,0,0)=:m_1$ lies on a regular part of $L_0$, and Theorem \ref {thm:kudr} can be applied to it.
We have a real $\lambda=1\in\R$ and a real-analytic $2\pi$-periodic first integral $I_1(f_1,f_2)=f_1$, as Theorem \ref {thm:kudr} (a) asserts.

The Hamiltonian $S^1$-action generated by $if_2^{\mathbb C}$ is given by the system
$\frac{\ddd p_1}{\ddd t}=-i\frac{\partial f_2}{\partial q_1}=-ip_1$,
$\frac{\ddd p_2}{\ddd t}=-i\frac{\partial f_2}{\partial q_2}=-ip_2$,
$\frac{\ddd q_1}{\ddd t}= i\frac{\partial f_2}{\partial p_1}=iq_1$,
$\frac{\ddd q_2}{\ddd t}= i\frac{\partial f_2}{\partial p_2}=iq_2$.
Its orbit $\gamma(t)=\varepsilon(e^{-it},0,0,0)$ with $\gamma(0)=(\varepsilon,0,0,0)=:m_1$
lies on a regular part of $L_0^{\mathbb C}$. This solution is $2\pi$-periodic, thus $\lambda=i\in i\R$.
The $\mathbb C$-conjugated path is
$\overline{\gamma(t)}=(\varepsilon e^{it},0,0,0)=\gamma(-t)$.
Thus $\gamma(t)$ and $\overline{\gamma(-t)}$ are homological in the regular part of the fiber $L_0^{\mathbb C}$, since they coincide.
We have an imaginary $\lambda=i\in i\R$ and a holomorphic $2\pi$-periodic first integral $I_2(f_1,f_2)=if_2$, moreover $iI_2(f_1,f_2)=-f_2$ is real-valued,
as Theorem \ref {thm:kudr:} (a) asserts. \qed
\end{Ex}

\section {Normalization of a torus action near a singular orbit. Elliptic, hyperbolic and twisting resonances} \label {subsec:mult}

In this section, we solve the questions (2a)--(2c) and (3a) from Introduction.

In particular, we describe any Hamiltonian torus action on a neighborhood of its orbit, and prove persistence of its canonical model under perturbations.

Such a torus action can be obtained e.g.\ from a Hamiltonian $\R^n$-action generated by the momentum map $F=(f_1,\dots,f_n)$ of an integrable Hamiltonian system, via either results of \cite{zung96, zung96a, bau:zung97, fom, BF, zung00, zung03} or our results of the previous section.
If we do so, we will obtain a Hamiltonian action of a torus generated by some functions of the form $I_j=I_j(f_1,\dots,f_n)$, $1\le j\le r+\varkappa$.
As we are mostly interested in the singular Lagrangian fibration (rather than specific commuting functions $f_1,\dots,f_n$), we allow ourselves to replace $f_1,\dots,f_n$ with $I_1,\dots,I_{r+\varkappa},f_{r+\varkappa+1},\dots,f_n$ where $\frac{\partial (I_1,\dots,I_{r+\varkappa})}{\partial (f_1,\dots,f_{r+\varkappa})}\ne 0$.
So, we can assume that some of the components of the momentum map generate a torus action.
In this section (except for Examples \ref {exa:res}, \ref {exa:res:}), we forget about other components of the momentum map (even about their existence). After that, one can study normal form of the system and its perturbations, see Examples \ref {exa:res}, \ref {exa:res:} below.

It is well known that a smooth action of a compact Lie group $G$ is linearizable on a neighbourhood of its fixed point (\cite {boc}, \cite[Sec 3.1.4]{cha}).
According to the Darboux-Weinstein theorem (which is an equivariant Darboux theorem), a smooth symplectic action of a compact Lie group $G$ is symplectically linearizable on a neighbourhood of its fixed point (\cite[Lecture 5]{wei}, \cite{del:mel}). This extends to arbitrary Lie groups and their arbitrary orbits under natural assumptions \cite{mar83, mar, gui:ste}.
Moreover, an integrable Hamiltonian system admitting a symplectic action of a compact Lie group $G$ preserving the momentum map of the system, is equivariantly fiberwise symplectomorphic to a linear model on a neighbourhood of an invariant compact nondegenerate singular orbit $\cal O$ \cite {zung:mir04}.

In this section, we formulate Theorems \ref {thm:period} and \ref {thm:period:} about
\begin{itemize}
\item reduction to a symplectic normal form, so-called canonical model, 
for a Hamiltonian torus action generated by smooth (real-analytic, respectively) functions (some of which are multiplied by $\sqrt{-1}$, respectively) on a neighbourhood of a (may be degenerate) singular orbit $\cal O$ in $M$ (in $M^{\mathbb C}$, resp.);
\item persistence of the canonical model of a torus action under small perturbations of the action in the class of Hamiltonian torus actions.
\end{itemize}
In particular, we will define discrete parameters that completely determine the canonical model up to symplectomorphism: so-called Williamson type of the orbit, $\varkappa$ tuples of integers and $r$ tuples of ratios called elliptic, hyperbolic and twisting resonances of the orbit, respectively.
Here $r$ and $r+\varkappa$ denote dimensions of the orbit and the torus, respectively.
Proofs of Theorems \ref {thm:period} and \ref {thm:period:} will be given in Sections \ref{sec:app:fixed}, \ref {sec:app:}.

\subsection {The linear model in elliptic case} \label {ref:subsec:ell}

In this subsection, we solve the questions (2a), (2c) and partially (3a) from Introduction.

Denote by $\lambda=(\lambda_1,\dots,\lambda_r)$ a linear coordinate system on a small ball $D^r$ of dimension $r$ centred at the origin, $\varphi=(\varphi_1,\dots,\varphi_r)$ a standard periodic coordinate system of the torus $(S^1)^r$, and $(x,y)=(x_1,y_1,\dots,x_{n-r},y_{n-r})$ a linear coordinate system on a small polydisc $(D^2)^{n-r}$ of dimension $2(n-r)$ centred at the origin. Consider the manifold
\begin{equation} \label {eq:mz:product}
V = D^r \times (S^1)^r \times (D^2)^{n-r},
\end{equation}
with the standard symplectic form $\sum\limits_{s=1}^r\ddd \lambda_s\wedge \ddd \varphi_s
+ \sum\limits_{j=1}^{n-r} \ddd x_{j}\wedge \ddd y_{j}$, and the following map:
\begin{equation} \label {eq:mz:H}
(\lambda,H) = (\lambda_1,\dots,\lambda_r,h_1,\dots,h_{k_e}):V\to\R^{r+k_e}
\end{equation}
where $r,k_e\in{\mathbb Z}_+$, $r+k_e\le n$,
\begin{equation} \label {eq:mz:hj}
h_j=\frac{x_{j}^2+y_{j}^2}2 \quad \mbox{for} \quad 1\le j\le k_e.
\end{equation}

Let $\Gamma$ be a group acting on the product $D^r \times (S^1)^r \times (D^2)^{n-r}$ by symplectomorphisms preserving the map $(\lambda,H)$.
We will say that the action of $\Gamma$ is {\em linear} (compare \cite[\S2.3]{zung:mir04}) if it satisfies the following property:

(L) {\em $\Gamma$ acts on the product $V = D^r \times (S^1)^r \times (D^2)^{n-r}$ componentwise; the action of $\Gamma$ on $D^r$ is trivial, its action on $(S^1)^r$ is by translations (w.r.t.\ the coordinate system $\varphi$), and its action on $(D^2)^{n-r}$ is linear w.r.t.\ the coordinate system $(x,y)$.}

Suppose now that $\Gamma$ is a finite group with a free symplectic action on $V$ that is linear (see (L) above) and preserves the map $(\lambda,H)$.
Then we can form the quotient symplectic manifold $V/\Gamma$, with a {\em $(S^1)^{r+\varkappa_e}$-action} on it generated by the following {\em momentum map}, all whose components are linear or quadratic functions:
\begin{equation} \label {eq:mz:F}
(\lambda, Q)=(\lambda_1,\dots,\lambda_r,Q_1\dots,Q_{\varkappa_e}) : V/\Gamma \to \R^{r+\varkappa_e}, \qquad
\mbox{where} \quad Q_\ell := \sum\limits_{j=1}^{k_e} {p_{j\ell}} h_{j},
\end{equation}
$1\le \ell\le \varkappa_e$, for some integers $\varkappa_e\in{\mathbb Z}_+$ and $p_{j\ell}\in{\mathbb Z}$. Suppose also that
\begin{itemize}
\item[(i)] $\rank\|p_{j\ell}\|=\varkappa_e$, thus $\varkappa_e\le k_e$,
\item[(ii)] $\Gamma$ acts on $(D^2)^{n-r}=
%(D^2)^{k_e}\times (D^2)^{n-r-k_e}
\underbrace{D^2 \times \dots \times D^2}_{k_e} \times \underbrace{D^2 \times \dots \times D^2}_{n-r-k_e}$ componentwise, and the induced action on $
(D^2)^{n-r-k_e}
%\underbrace{D^2 \times \dots \times D^2}_{n-r-k_e}
$ is by involutions.
\end{itemize}
The set
$$
{\cal O} := \{\lambda_s=x_j=y_j=0\}/\Gamma \subset V/\Gamma
$$
is a rank-$r$ orbit of the above $(S^1)^{r+\varkappa_e}$-action on $V/\Gamma$.

\begin{Definition} \label {def:lin:model}
Consider the above Hamiltonian $(S^1)^{r+\varkappa_e}$-action on $V/\Gamma$ generated by the momentum map (\ref {eq:mz:F}), satisfying the assumptions (i) and (ii) from above.
We will call this action the {\em linear $(S^1)^{r+\varkappa_e}$-action} (or {\em linear model}) of {\em rank} $r$, {\em Williamson type} $(k_e,0,0)$, {\em elliptic resonances}
\begin{equation} \label {eq:res:inf:ell}
(p_{1\ell}:\dots:p_{k_e,\ell})\in{\mathbb Q}P^{k_e-1}, \qquad 1\le \ell\le \varkappa_e,
\end{equation}
and {\em twisting group} $\Gamma$ (or, more precisely, {\em twisting linear action} of
$\Gamma$ on $V$), provided that the integer $k_e$ cannot be made smaller via a linear change of coordinates $(x,y)$ on $(D^2)^{n-r}$ (which is equivalent to the fact that, for each $j\in\{1,\dots,k_e\}$, either $\sum\limits_{\ell=1}^{\varkappa_e}|p_{j\ell}|>0$ or there exists $a\in\{1,\dots,r\}$ such that $2q_{\psi^a,j}\not\subset{\mathbb Z}$, see (\ref {eq:res:mon:ell}) below).
\end{Definition}

\begin{RemarkDefinition} \label {rem:period}
(A) Since $\Gamma$ freely acts on $V$ componentwise and the induced action on $(S^1)^r$ is by translations, we can regard $\Gamma$ as a subgroup of $(S^1)^r$.
One can show\footnote{Indeed, it follows from Theorem \ref {thm:period:} that the action of generators $\psi^a := p(\gamma^a)$ of $\Gamma$ on $V$ has the form
$$
(\lambda,\varphi,z_1,\dots,z_{n-r}) \mapsto (\lambda,\varphi+\psi^a,e^{2\pi i q_{\gamma^a,1}}z_1,\dots,e^{2\pi i q_{\gamma^a,n-r}}z_{n-r}), \quad 1\le a\le r,
$$
for some $q_{\gamma^a,j}\in{\mathbb Q}$. Here $\gamma^1,\dots,\gamma^r$ is a generating set of the lattice $p^{-1}(\Gamma)\subset\R^r$, $p:\R^r\to(S^1)^r$ denotes the projection. Since $p^{-1}(\Gamma)$ is a lattice in $\R^r$, there exists a unique linear map $\R^r\to\R^{n-r}$ sending $\gamma^a\mapsto2\pi(q_{\gamma^a,1},\dots,q_{\gamma^a,n-r})$, $1\le a\le r$. Clearly, this linear map has the form $\gamma\mapsto(\langle m_1,\gamma\rangle,\dots,\langle m_{n-r},\gamma\rangle)$, $\gamma\in\R^r$, for some $m_1,\dots,m_{n-r}\in\R^r$. From the short exact sequence $0\to2\pi{\mathbb Z}^r\to p^{-1}(\Gamma)\to\Gamma\to0$, we conclude that $\langle m_j,\gamma\rangle\in2\pi\mathbb Z$, provided that $\gamma\in2\pi{\mathbb Z}^r$. Therefore $m_j\in{\mathbb Z}^r$. This proves (\ref {eq:extension}).}
that a (free) twisting linear action of $\Gamma$ on $V$ has the form
\begin{equation} \label {eq:extension}
(\lambda,\varphi,z_1,\dots,z_{n-r}) \mapsto (\lambda,\varphi+\psi,e^{i\langle m_1,\psi\rangle}z_1,\dots, e^{i\langle m_{n-r},\psi\rangle}z_{n-r}), \quad 
\psi\in\Gamma\subset(S^1)^r,
\end{equation}
for some $m_1,\dots,m_{n-r}\in{\mathbb Z}^r$, for an appropriate choice of symplectic coordinates $(x,y)$ on $(D^2)^{n-r}$ satisfying (\ref {eq:mz:F}) (corresponding to a root decomposition of $\R^{2(n-r)}$ w.r.t.\ the commuting $(S^1)^{\varkappa_e}$-action and $\Gamma$-action).
Here we used the notation $m_j=:(m_{j1},\dots,m_{jr})$, $\psi:=(\psi_1,\dots,\psi_r)$,
$\langle m_j,\psi\rangle:=\sum\limits_{s=1}^r m_{js}\psi_s$, $z_j:=x_j+iy_j$, $1\le j\le n-r$.
In other words, a (free) twisting linear action of $\Gamma$ on $V$ can be extended to a (free) linear Hamiltonian action of $(S^1)^r\supset\Gamma$ 
preserving the momentum map (\ref {eq:mz:F}).

(B) We will call such a twisting linear action of $\Gamma$ on $V$ the {\em twisting linear action} with {\em twisting} {\em resonances}
\begin{equation} \label {eq:res:mon:ell}
(q_{\psi^a,1},\dots,q_{\psi^a,n-r})\in({\mathbb Q}/{\mathbb Z})^{n-r}, \qquad
\mbox{where} \quad
q_{\psi^a,j}:=\langle m_j,\frac{\psi^a}{2\pi}\rangle\mod 1\in{\mathbb Q}/{\mathbb Z},
\end{equation}
$1\le a\le r$, $1\le j\le n-r$, and $2q_{\psi^a,j}\subset{\mathbb Z}$ for $k_e+1\le j\le n-r$ (due to assumption (ii) from above). Here $\gamma^1,\dots,\gamma^r$ denote a basis of the homology group $H_1({\cal O})\simeq p^{-1}(\Gamma)\subset2\pi{\mathbb Q}^r\subset\R^r$, thus $\psi^a := p(\gamma^a)$, $1\le a\le r$, is a generating set of the group $\Gamma$, where $p:\R^r\to(S^1)^r$ is the projection.

(C) We always may assume that each basic cycle $\gamma^a$ has coordinates $\gamma^a_s=\frac{2\pi}{N_a}\delta^a_s$, $1\le a,s\le r$, where $N_a$ is a positive integer.
In this case, we have $q_{\psi^a,j} = \frac{m_{ja}}{N_a}\mod 1$.
\end{RemarkDefinition}

\begin{Definition} \label {def:period}
In notations of Definition \ref {def:lin:model}, Remark and Definition \ref {rem:period}, we will call the linear $(S^1)^{r+\varkappa_e}$-action on $V/\Gamma$ from Definition \ref {def:lin:model} the {\em linear $(S^1)^{r+\varkappa_e}$-action} with
\begin{itemize}
\item dimension $2n$, rank $r$, {\em Williamson type} $(k_e,0,0)$,
\item {\em elliptic resonances} (\ref {eq:res:inf:ell}) assigned to the basic cycles of the subtorus $(S^1)^{\varkappa_e}$ of $(S^1)^{r+\varkappa_e}$ (which is the isotropy subgroup of some and, hence, any point of $\cal O$),
\item {\em twisting} {\em resonances} (\ref {eq:res:mon:ell}) assigned to the basic cycles of the orbit ${\cal O} = \{\lambda_s=x_j=y_j=0\}/\Gamma \approx (S^1)^r/\Gamma$ relatively the action of the subtorus $(S^1)^r$ of $(S^1)^{r+\varkappa_e}$.
\end{itemize}
Clearly, dimension, rank, the Williamson type together with elliptic and twisting resonances completely determine the symplectic manifold $V/\Gamma$ with the linear $(S^1)^{r+\varkappa_e}$-action on it, up to symplectomorphism. See also Remark \ref {rem:eigen} about uniqueness of the Williamson type and the resonances.
\end{Definition}

Now we can formulate our result in the elliptic case, which is the symplectic normalization theorem for singular orbits of Hamiltonian torus actions:

\begin{Theorem} \label {thm:period}
Suppose we are given an effective Hamiltonian action of the $(r+\varkappa_e)$-torus $(S^1)^{r+\varkappa_e}$ generated by $C^\infty$-smooth functions $I_1,\dots,I_r$, $J_1,\dots,J_{\varkappa_e}$ on a $C^\infty$-smooth symplectic manifold $(M,\Omega)$.
Suppose a point $m_0\in M$ is fixed under the $(S^1)^{\varkappa_e}$-subaction, and its orbit ${\cal O}_{m_0}$ is $r$-dimensional.
Then:

{\rm (a)}
There exist an invariant neighbourhood $U$ of ${\cal O}_{m_0}$,
a finite group $\Gamma$, a linear $(S^1)^{r+\varkappa_e}$-action of rank $r$ on the symplectic manifold $V/\Gamma$ given by (\ref {eq:mz:product})--(\ref {eq:mz:F}), and a smooth action-preserving symplectomorphism $\phi$ from $U$ to $V/\Gamma$, that sends the orbit ${\cal O}_{m_0}$ to the torus ${\cal O}=\{\lambda_s=x_j=y_j=0\}/\Gamma$. If the system is real-analytic, $\phi$ is real-analytic too.

In particular, dimensions of the torus $(S^1)^{r+\varkappa_e}$ and its subtorus 
$(S^1)^{\varkappa_e}$ satisfy the estimates $r+\varkappa_e\le r+k_e\le n$. If $r+\varkappa_e=n$ then the orbit ${\cal O}_{m_0}$ is nondegenerate.

{\rm (b)} The symplectic normalization
(\ref {eq:mz:product})--(\ref {eq:mz:F}) for a torus action is persistent under $C^\infty$-smooth perturbations in the following sense.
For any integer $k\ge5$ and any neighbourhood $U'$ of the orbit ${\cal O}_{m_0}$ having a compact closure $\overline{U'}\subset U$, there exists $\varepsilon>0$ satisfying the following. Suppose we are given a $C^\infty$-smooth 
Hamiltonian $(S^1)^{r+\varkappa_e}$-action on $\tilde M$ (w.r.t.\ the ``perturbed'' symplectic structure $\tilde\Omega$ that is $\varepsilon$-close to $\Omega$ in $C^{k+n-r-2}-$norm) generated by a (``perturbed'') momentum map $\tilde F=(\tilde I_1,\dots,\tilde I_r,\tilde J_1,\dots,\tilde J_{\varkappa_e})$ that is $\varepsilon$-close to $F=(I_1,\dots,I_r,J_1,\dots,J_{\varkappa_e})$ in $C^k-$norm, where $U\subset\tilde M\subset M$.
Then there exist an invariant neighbourhood $\tilde U\supset U'$,
and a smooth action-preserving symplectomorphism $\tilde\phi$ from $\tilde U$ to $V/\Gamma$
that is $O(\varepsilon)-$close to $\phi$ in $C^{k-4}-$norm. 
If the systems are real-analytic, $\tilde\phi$ is real-analytic too.
If the system depends smoothly (resp., analytically) on a local parameter (i.e.\ we have a local family of systems), $\phi$ can also be chosen to depend smoothly (resp., analytically) on that parameter.
\end{Theorem}

\begin{Remark} \label {rem:eigen}
(A) Consider a linear $(S^1)^{r+\varkappa_e}$-action on $V/\Gamma$ (see Definition \ref {def:lin:model}, Remark and Definition \ref {rem:period}, Definition \ref {def:period}).
What is a dynamical meaning of the elliptic resonances (\ref {eq:res:inf:ell}) and the twisting resonances (\ref {eq:res:mon:ell})?
The eigenvalues of the linearized $\ell$-th infinitesimal generator of the $\varkappa_e$-subtorus action on $(D^2)^{k_e}$ are pure imaginary and are in (\ref {eq:res:inf:ell}) resonance, which is exactly the elliptic resonance assigned to this generator. Further, the multipliers of a cycle $\gamma\in H_1({\cal O})\approx p^{-1}(\Gamma)\subset2\pi{\mathbb Q}^r\subset\R^r$ of the torus ${\cal O}=\{\lambda_s=x_j=y_j=0\}/\Gamma$ under the $r$-subtorus action on $(D^2)^{k_e}$ equal $e^{\pm2\pi i q_{\psi,j}}$, where $\psi:=p(\gamma)\in\Gamma$ and $q_{\psi,j}\mod1:=\langle m_j,\frac\psi{2\pi}\rangle\in{\mathbb Q}/{\mathbb Z}$, $1\le j\le n-r$. So, they are in $(q_{\psi,1},\dots,q_{\psi,n-r})$ resonance, which is exactly the twisting resonance assigned to the cycle $\gamma$. Here $p:\R^r\to(S^1)^r$ denotes the projection.

(B) Clearly, Williamson type is well defined, i.e.\ completely determined by the (``hidden'') torus-symmetry of the orbit $\cal O$.
Are the elliptic resonances also well defined (perhaps, up to some natural transformations)?
Clearly, the subtorus $(S^1)^{\varkappa_e}$ of the torus $(S^1)^{r+\varkappa_e}$ trivially acts on $\cal O$, moreover this subtorus is the isotropy subgroup of any point of $\cal O$, thus this subtorus is well defined, i.e.\ does not depend on the choice of generators of the $(S^1)^{r+\varkappa_e}$-action.
It follows that the elliptic resonances (\ref {eq:res:inf:ell}), $1\le\ell\le \varkappa_e$, are well defined up to replacing the vectors $(p_{1\ell},\dots,p_{k_e,\ell})$, $1\le\ell\le \varkappa_e$, by their linear combinations with integer coefficients forming a nondegenerate $\varkappa_e\times \varkappa_e$-matrix.

(C) Are the twisting resonances also well defined (perhaps, up to some natural transformations)?
Suppose we allow ourselves to change the generators of the subtorus $(S^1)^{r}$ that acts locally-freely on our manifold $V/\Gamma$.
At the same time, suppose that we fixed basic cycles $\gamma^1,\dots,\gamma^r\in H_1({\cal O})\subset{\mathbb Q}^r\subset\R^r$ of the orbit ${\cal O}$, and we want that the cycle $\gamma^a$ has coordinates $\gamma^a_s=\frac{2\pi}{N_a}\delta^a_s$ w.r.t.\ to the generators of the subtorus $(S^1)^r$, both before and after the change, may be with different integers $N_a$, $1\le a,s\le r$ (see Remark and Definition \ref {rem:period} (C)).
Thus, the above change is equivalent to replacing the coordinates $\lambda_s, \varphi_s, z_j$ with
$$
\hat\lambda_s := \nu_s \lambda_s + \sum\limits_{\ell=1}^{\varkappa_e} n_{\ell s} Q_\ell, \qquad
\hat\varphi_s := \frac1{\nu_s}\varphi_s, \qquad
\hat z_j := e^{-i  \langle (pn)_j , \varphi \rangle}z_j
$$
for $1\le s\le r$ and $1\le j\le n-r$.
Here $\nu_s,n_{1 s},\dots,n_{\varkappa_e s}$ are coprime integers such that $\nu_s\ne 0$, and we denoted
$(pn)_{js}:=\sum\limits_{\ell=1}^{\varkappa_e} p_{j\ell}n_{\ell s}$ if $1\le j\le k_e$,
$(pn)_{js}:=0$ if $k_e+1\le j\le n-r$,
$(pn)_j:=((pn)_{j1},\dots,(pn)_{jr})$,
$\varphi:=(\varphi_1,\dots,\varphi_r)$,
$\langle (pn)_j , \psi\rangle := \sum\limits_{s=1}^r (pn)_{js}\psi_s$.
This replacement will lead to the following replacements: $\psi_j^a$ with $\hat\psi_j^a$ such that $\nu_j\hat\psi_j^a=\psi_j^a$, and
$m_j\in{\mathbb Z}^r$ with $\hat m_{js}:=\nu_s(m_{js}-(pn)_{js})$.
Thus, the twisting resonances (\ref {eq:res:mon:ell}) will be replaced with
$$
(\hat q_{\psi^a,1},\dots,\hat q_{\psi^a,n-r})\in({\mathbb Q}/{\mathbb Z})^{n-r}, \quad
\mbox{where} \quad
\hat q_{\psi^a,j} = q_{\psi^a,j} + \langle (pn)_j,\frac{\psi^a}{2\pi}\rangle
\mod 1\in{\mathbb Q}/{\mathbb Z},
$$
$1\le a\le r$, $1\le j\le n-r$.
If we recall that $\gamma^a_s=\frac{2\pi}{N_a}\delta^a_s$, $1\le a,s\le r$, then 
$\hat q_{\psi^a,j} = q_{\psi^a,j} + (pn)_{ja}/N_a$.
In other words, the twisting resonances (\ref {eq:res:mon:ell}) are well defined up to adding any linear combinations of ``extended'' elliptic resonances (\ref {eq:res:inf:ell}) with rational coefficients forming a $r\times \varkappa_e$-matrix. Here, by the $\ell$-th {\em extended elliptic resonance}, $1\le \ell\le \varkappa_e$, we mean the vector $(p_{1,\ell},\dots,p_{k_e,\ell},0,\dots,0)\in{\mathbb Z}^{n-r}$.
\end{Remark}

\begin{Ex} \label {exa:res}
Suppose ${\cal O}$ is a compact $r$-dimensional orbit of the momentum map of a real-analytic integrable Hamiltonian system with $n$ degrees of freedom. 
%In other words, ${\cal O}$ is a rank-$r$ local singularity.
Thus, we have a rank-$r$ local singularity at ${\cal O}$. 
Suppose this singularity has one of the following types:
\begin{itemize}
\item[\rm(a)] a parabolic orbit with resonance $\ell/s$ (see \cite{kal} or \S \ref {subsec:kal}), $s\ge5$, given by a momentum map $F=(H,I):(D^1_{(\lambda)}\times S^1_{(\varphi)}\times D^2_{(x,y)})/{\mathbb Z}_s \to\R^3$ where $I=\lambda$, $H=\re(z^s)+|z|^4+\lambda|z|^2$,
we denote $z=x+iy$; a generator of ${\mathbb Z}_s$ acts on $D^1\times S^1\times D^2$ by the transformation $(\lambda,\varphi,z)\mapsto(\lambda,\varphi+2\pi/s,e^{2\pi i\ell/s}z)$, $s,\ell\in\mathbb Z$, $0\le\ell<s$, $(s,\ell)=1$;
\item[\rm(b)] an integrable Hamiltonian Hopf bifurcation with resonance $p:q$ \cite{dui,han}, $p,q\in\mathbb Z$, $0<p<q$, $(p,q)=1$, $\frac pq\ne\frac13$, given by a momentum map $F=(H,I,J):D^1_{(\lambda)}\times S^1_{(\varphi)}\times D^4_{(z_1,z_2)}\to\R^3$ where $I=\lambda$, 
$J=p\frac{|z_1|^2}2+q\frac{|z_2|^2}2$, $H=\re(z_1^p\bar z_2^q)+a|z_2|^4+\lambda|z_2|^2$, one denotes $z_j=x_j+iy_j$, $a$ is a real parameter;
\item[\rm(c)] a normally-elliptic parabolic orbit 
with resonance $\ell/s$ (compare \cite{bro93}), $s\ge5$, 
given by a momentum map $F=(H,I,J):(D^1_{(\lambda)}\times S^1_{(\varphi)}\times D^4_{(z_1,z_2)})/{\mathbb Z}_s\to\R^3$ where $I=\lambda$, $J=\frac12|z_1|^2$, 
$H=\re(z_2^s)+|z_2|^4+(I\pm J)|z_2|^2$, one denotes $z_j=x_j+iy_j$, 
a generator of ${\mathbb Z}_s$ acts on $D^1\times S^1\times D^4$ by the transformation $(\lambda,\varphi,z_1,z_2)\mapsto(\lambda,\varphi+2\pi/s,z_1,e^{2\pi i\ell/s}z_2)$, $s,\ell\in\mathbb Z$, $0\le\ell<s$, $(s,\ell)=1$;
\item[\rm(d)] a Hamiltonian swallow-tail bifurcation with resonance $\ell/5$,
given by a momentum map $F=(H,I_1,I_2):(D^2_{(\lambda_1,\lambda_2)}\times (S^1)^2_{(\varphi_1,\varphi_2)}\times D^2_{(x,y)})/{\mathbb Z}_5 \to\R^3$ where $I_1=\lambda_1$, $I_2=\lambda_2$, 
$H=\re(z^5)+|z|^6+\lambda_2|z|^4+\lambda_1|z|^2$, one denotes $z=x+iy$,
a generator of ${\mathbb Z}_5$ acts on $D^2\times (S^1)^2\times D^2$ by the transformation $(\lambda_1,\lambda_2,\varphi_1,\varphi_2,z)\mapsto(\lambda_1,\lambda_2,\varphi_1+2\pi/5,\varphi_2,e^{2\pi i\ell/5}z)$, $\ell\in\mathbb Z$, $1\le\ell<5$.
\end{itemize}
Clearly, this orbit is of finite type \cite {zung03} (see \S\ref{subsec:adj}) and degenerate.
Due to the Zung result \cite {zung03}, there exists a locally-free $F$-preserving Hamiltonian $(S^1)^r$-action on some neighbourhood of ${\cal O}$. It can be shown from Theorem \ref {thm:kudr} or \ref {thm:kudr:} that this action extends to a ``hidden'' $(S^1)^{r+\varkappa_e}$-symmetry, i.e.\ to an effective (not locally-free) $F$-preserving Hamiltonian $(S^1)^{r+\varkappa_e}$-action, with $r+\varkappa_e=n-1$. Moreover, this action is persistent under real-analytic integrable perturbations. By Theorem \ref {thm:period}, the latter action is symplectomorphic to a linear $(S^1)^{r+\varkappa_e}$-action, that is also persistent under real-analytic integrable perturbations. In Table, we show Williamson type $(k_e,0,0)$, elliptic and twisting resonances (Definition \ref {def:period}) of this linear $(S^1)^{r+\varkappa_e}$-action, for each type of the singularity ${\cal O}$ from above.
\end{Ex}

\begin{tabular}{|c|c|c|c|c|c|c|}
\hline
Case &$n$& \multicolumn{2}{c|}{Subtori dim's} & Williamson    & \multicolumn{2}{c|}{Resonances} \\
\cline{3-4} \cline{6-7}
&  & \ \ \, $r$ \ \ & $\varkappa_e$ & type $(k_e,0,0)$ & elliptic & twisting \\
\hline
\hline
(a) & 2 & 1 &0& $(1,0,0)$ & no & $\ell/s\mod1\in{\mathbb Q}/{\mathbb Z}$ \\
\hline
(b) & 3 & 1 &1& $(2,0,0)$ & $(p:q)\in{\mathbb Q}P^{1}$ & $(0,0)\in({\mathbb Q}/{\mathbb Z})^2$ \\
\hline
(c) & 3 & 1 &1& $(2,0,0)$ & $(0:1)\in{\mathbb Q}P^{1}$ & $(\ell/s\mod1,0)\in({\mathbb Q}/{\mathbb Z})^2$ \\
\hline
(d) & 3 & 2 &0& $(1,0,0)$ & no & $\ell/5\mod1,\ 0\mod1\in{\mathbb Q}/{\mathbb Z}$ \\
\hline
\end{tabular}

\subsection {The linear model in general case}

In this subsection, we solve the questions (2b), (2c) and (3a) from Introduction.

As above, consider the manifold $V = D^r \times (S^1)^r \times (D^2)^{n-r}$ as in (\ref {eq:mz:product}),
with the standard symplectic form $\sum\limits_{s=1}^r\ddd \lambda_s\wedge \ddd \varphi_s
+ \sum\limits_{j=1}^{n-r} \ddd x_{j}\wedge \ddd y_{j}$, and the following map:
\begin{equation} \label {eq:mz:H:}
(\lambda,H) = (\lambda_1,\dots,\lambda_r,h_1,\dots,h_{2k_f+k_e+k_h}):V\to\R^{r+2k_f+k_e+k_h}
\end{equation}
where $r,k_e,k_h,k_f\in{\mathbb Z}_+$, $r+k_e+k_h+2k_f\le n$,
\begin{equation} \label {eq:mz:hj:}
\begin{array} {lll}
h_{2j-1}=\frac{x_{2j-1}^2+y_{2j-1}^2}2 - \frac{x_{2j}^2+y_{2j}^2}2 && \mbox{and}  \\
h_{2j}=x_{2j-1}y_{2j} + x_{2j}y_{2j-1} & & \mbox{for} \quad 1\le j\le k_f, \\
h_j=\frac{x_{j}^2+y_{j}^2}2 && \mbox{for} \quad 2k_f+1\le j\le 2k_f+k_e, \\
h_j=x_{j}y_{j} && \mbox{for} \quad 2k_f+k_e+1\le j\le 2k_f+k_e+k_h.
\end{array}
\end{equation}

Let $\Gamma$ be a group acting on the product $D^r \times (S^1)^r \times (D^2)^{n-r}$ by symplectomorphisms preserving the map $(\lambda,H)$ given by (\ref {eq:mz:H:}), (\ref {eq:mz:hj:}).
Suppose the group $\Gamma$ is finite, and its action on $V$ is free and linear (see Property (L) in \S \ref {ref:subsec:ell}).
Then we can form the quotient symplectic real-analytic manifold $V/\Gamma$, with a real-analytic {\em momentum map}
\begin{equation} \label {eq:mz:F:}
(\lambda,Q) = (\lambda_1,\dots,\lambda_r,Q_1\dots,Q_\varkappa) : V/\Gamma \to \R^{r+\varkappa}, \quad
\mbox{where} \ Q_\ell := \sum\limits_{j=1}^{k_e+k_h+2k_f} {p_{j\ell}} h_{j},
\end{equation}
$1\le \ell\le \varkappa=\varkappa_e+\varkappa_h$,
and a {\em $(S^1)^{r+\varkappa}$-action} on its small open complexification $(V/\Gamma)^{\mathbb C}$ generated by the map
\begin{equation} \label {eq:mz:F::}
(\lambda_1,\dots,\lambda_r, Q_1\dots,Q_{\varkappa_e}, iQ_{\varkappa_e+1},\dots,iQ_{\varkappa_e+\varkappa_h}) : V/\Gamma \to \R^{r+\varkappa_e}\times(i\R)^{\varkappa_h},
\end{equation}
for some integers $\varkappa_e,\varkappa_h\in{\mathbb Z}_+$ and $p_{j\ell}\in{\mathbb Z}$ such that each component $Q_\ell$ (respectively $iQ_\ell$) of the map (\ref {eq:mz:F::}) is a linear combination of elliptic (respectively hyperbolic) $h_j$, i.e.\ $p_{j\ell}=0$ if at least one of the following conditions (i) and (ii) holds:
\begin{itemize}
\item[(i)] $1\le\ell\le\varkappa_e$ and 
$j\in\{2,4,6,\dots,2k_f\}\cup\{2k_f+k_e+1,\dots,2k_f+k_e+k_h\}$,
\item[(ii)] $\varkappa_e<\ell\le \varkappa_e+\varkappa_h$ and 
$j\in\{1,3,5,\dots,2k_f-1\}\cup\{2k_f+1,\dots,2k_f+k_e\}$.
\end{itemize}
Suppose also that
\begin{itemize}
\item[(iii)] $\rank\|p_{j\ell}\|=\varkappa_e+\varkappa_h$, thus $\varkappa_e\le k_e+k_f$, $\varkappa_h\le k_h+k_f$, and $\varkappa_e+\varkappa_h\le k_e+k_h+2k_f$,
\item[(iv)] $\Gamma$ acts on $(D^2)^{n-r}=
\underbrace{D^2 \times \dots \times D^2}_{k_e+k_h+2k_f} \times \underbrace{D^2 \times \dots \times D^2}_{n-r-k_e-k_h-2k_f}$ componentwise, and the induced action on $
(D^2)^{n-r-k_e-k_h-2k_f}
%\underbrace{D^2 \times \dots \times D^2}_{n-r-k_e-k_h-2k_f}
$ is by involutions.
%(D^2)^{k_e+k_h+2k_f}\times (D^2)^{n-r-k_e-k_h-2k_f}$ componentwise, and the induced action on $(D^2)^{n-r-k_e-k_h-2k_f}$ is by involutions.
\end{itemize}
The set
$$
{\cal O}:=\{\lambda_s=x_j=y_j=0\}/\Gamma \subset V/\Gamma
$$
is a rank-$r$ orbit of the above $(S^1)^{r+\varkappa_e+\varkappa_h}$-action on $(V/\Gamma)^{\mathbb C}$.

\begin{Definition} \label {def:lin:model:}
Consider the above Hamiltonian $(S^1)^{r+\varkappa_e+\varkappa_h}$-action on $(V/\Gamma)^{\mathbb C}$ generated by the map (\ref {eq:mz:F::}), satisfying the assumptions (i)--(iv) from above. 
We will call this action the {\em linear $(S^1)^{r+\varkappa_e+\varkappa_h}$-action} (or {\em linear model}) of {\em rank} $r$, {\em Williamson type} $(k_e,k_h,k_f)$, {\em elliptic resonances}
\begin{equation} \label {eq:res:inf:e}
(p_{1,\ell}:-p_{1,\ell}:p_{3,\ell}:-p_{3,\ell}:\dots:p_{2k_f-1,\ell}:-p_{2k_f-1,\ell}:
p_{2k_f+1,\ell}:\dots:p_{2k_f+k_e,\ell})\in{\mathbb Q}P^{2k_f+k_e-1},
\end{equation}
$1\le\ell\le\varkappa_e$, {\em hyperbolic resonances}
\begin{equation} \label {eq:res:inf:h}
(p_{2,\ell}:p_{2,\ell}:p_{4,\ell}:p_{4,\ell}:\dots:p_{2k_f,\ell}:p_{2k_f,\ell}:
p_{2k_f+k_e+1,\ell}:\dots:p_{2k_f+k_e+k_h,\ell})\in{\mathbb Q}P^{2k_f+k_h-1},
\end{equation}
$\varkappa_e+1\le \ell\le \varkappa_e+\varkappa_h$, and {\em twisting group} $\Gamma$ (or, more precisely, {\em twisting linear action} of $\Gamma$ on $V$), provided that the integer $k_e+k_h+2k_f$ cannot be made smaller via a linear change of coordinates $(x,y)$ on $(D^2)^{n-r}$ (which is equivalent to the fact that, for each $j\in\{1,3,5,\dots,2k_f-1\}\cup\{2k_f+1,\dots,2k_f+k_e\}$, 
either $\sum\limits_{\ell=1}^{\varkappa_e}|p_{j\ell}|>0$ or there exists $a\in\{1,\dots,r\}$ such that $2q_{\psi^a,j}\not\subset{\mathbb Z}$, see (\ref {eq:res:mon}) below;
furthermore for each $j\in\{2,4,6,\dots,2k_f\}\cup\{2k_f+k_e+1,\dots,2k_f+k_e+k_h\}$, we have $\sum\limits_{\ell=\varkappa_e+1}^{\varkappa_e+\varkappa_h}|p_{j\ell}|>0$). 
\end{Definition}

\begin{RemarkDefinition} \label {rem:period:}
(A) Since $\Gamma$ freely acts on $V$ componentwise and the induced action on $(S^1)^r$ is by translations, we can regard $\Gamma$ as a subgroup of $(S^1)^r$.
One can show similarly to (\ref {eq:extension}) that a (free) twisting linear action of $\Gamma$ on $V$ has the form
$$
(\lambda,\varphi,z_1,\dots,z_{n-r}) \mapsto
(\lambda,\varphi+\psi, e^{i\langle m_1,\psi\rangle}z_1, e^{-i\langle m_1,\psi\rangle}z_2,\dots, e^{i\langle m_{k_f},\psi\rangle}z_{2k_f-1}, e^{-i\langle m_{k_f},\psi\rangle}z_{2k_f},
$$
$$
e^{i\langle m_{k_f+1},\psi\rangle}z_{2k_f+1},\dots,
e^{i\langle m_{k_f+k_e},\psi\rangle}z_{2k_f+k_e},
\chi_{1}(\psi)z_{2k_f+k_e+1},\dots,
\chi_{k_h}(\psi)z_{2k_f+k_e+k_h},
$$
$$
e^{i\langle m_{k_f+k_e+k_h+1},\psi\rangle}z_{2k_f+k_e+k_h+1},\dots,
e^{i\langle m_{n-r-k_f},\psi\rangle}z_{n-r}), \quad \psi\in\Gamma\subset(S^1)^r,
$$
for some integer vectors $m_1,\dots,m_{k_f+k_e}$, $m_{k_f+k_e+k_h+1},\dots,m_{n-r-k_f}\in{\mathbb Z}^r$ and characters $\chi_1,\dots,\chi_{k_h}:\Gamma\to\{1,-1\}$, for an appropriate choice of coordinates $(x,y)$ on $(D^2)^{n-r}$ satisfying (\ref {eq:mz:F:}) (corresponding to a root decomposition of $\R^{2(n-r)}$ w.r.t.\ the commuting $(S^1)^{\varkappa_e+\varkappa_h}$-action and $\Gamma$-action).
Here we used the notation $m_j=:(m_{j1},\dots,m_{jr})$, $\psi=(\psi_1,\dots,\psi_r)\in(\R/2\pi\mathbb Z)^r$,
$\langle m_j,\psi\rangle:=\sum\limits_{a=1}^r m_{ja}\psi_a$, $z_j:=x_j+iy_j$.

In other words, the twisting group $\Gamma$ freely acts on the product $V=D^r \times (S^1)^r \times (D^2)^{2k_f+k_e}\times (D^2)^{k_h}\times (D^2)^{n-r-2k_f-k_e-k_h}$ componentwise, where its action on the ``hyperbolic'' component $(D^2)^{k_h}$ has the form
$$
(z_{2k_f+k_e+1},\dots, z_{2k_f+k_e+k_h}) \mapsto
(\chi_{1}(\psi)z_{2k_f+k_e+1},\dots,\chi_{k_h}(\psi)z_{2k_f+k_e+k_h}),
\qquad \psi \in \Gamma,
$$
while its action on the product $D^r \times (S^1)^r \times (D^2)^{2k_f+k_e}\times (D^2)^{n-r-2k_f-k_e-k_h}$ of the remaining components can be extended to a (free) linear Hamiltonian action of $(S^1)^r\supset\Gamma$ on $V$
preserving the momentum map (\ref {eq:mz:F:}).

(B) We will call such a twisting linear action of $\Gamma$ on $V$ the {\em twisting linear action} with {\em twisting} {\em resonances}
\begin{equation} \label {eq:res:mon}
(q_{\psi^a,1},-q_{\psi^a,1},\dots,q_{\psi^a,k_f},-q_{\psi^a,k_f}, \
q_{\psi^a,k_f+1},\dots,q_{\psi^a,n-r-k_f})
\in({\mathbb Q}/{\mathbb Z})^{n-r},
\end{equation}
$1\le a\le r$, where
$q_{\psi^a,j}:=\langle m_j,\frac{\psi^a}{2\pi}\rangle\mod 1\in{\mathbb Q}/{\mathbb Z}$ for $j\in\{1,\dots,k_f+k_e\}\cup\{k_f+k_e+k_h+1,\dots,n-r-k_f\}$,
$q_{\psi^a,j}:=\frac{1-\chi_{j}(\psi^a)}4$ for $k_f+k_e+1\le j\le k_f+k_e+k_h$, and
$2q_{\psi^a,j}\subset{\mathbb Z}$ for $k_f+k_e+k_h+1\le j\le n-r-k_f$ (due to assumption (iv) from above).
Here $\gamma^1,\dots,\gamma^r$ denote a basis of the homology group $H_1({\cal O}_{m_0})\simeq p^{-1}(\Gamma)\subset2\pi{\mathbb Q}^r\subset\R^r$, thus $\psi^a := p(\gamma^a)$, $1\le a\le r$, is a generating set of $\Gamma$, $p:\R^r\to(S^1)^r$ is the projection.

(C) Similarly to the elliptic case (see Remark and Definition \ref {rem:period} (C)), we
always may assume that $\gamma^a_s=\frac{2\pi}{N_a}\delta^a_s$, $1\le a,s\le r$, where
$N_a$ is a positive integer.
In this case, we have $q_{\psi^a,j} = \frac{m_{ja}}{N_a}\mod 1$.
\end{RemarkDefinition}

\begin{Definition} \label {def:period:}
We will call the linear $(S^1)^{r+\varkappa_e+\varkappa_h}$-action on $V/\Gamma$ from Definition \ref {def:lin:model:} the {\em linear $(S^1)^{r+\varkappa_e+\varkappa_h}$-action} with
\begin{itemize}
\item dimension $2n$, rank $r$, {\em Williamson type} $(k_e,k_h,k_f)$,
\item {\em elliptic resonances} (\ref {eq:res:inf:e}) and {\em hyperbolic resonances} (\ref {eq:res:inf:h}) assigned to the basic cycles of the subtori $(S^1)^{\varkappa_e}$ and $(S^1)^{\varkappa_h}$ of $(S^1)^{r+\varkappa_e+\varkappa_h}$, respectively,
\item {\em twisting} {\em resonances} (\ref {eq:res:mon}) assigned to the basic cycles of the orbit ${\cal O}=\{\lambda_s=x_j=y_j=0\}/\Gamma\approx(S^1)^r/\Gamma$ 
relatively the action of the subtorus $(S^1)^r$ of $(S^1)^{r+\varkappa_e+\varkappa_h}$.
\end{itemize}
Clearly, dimension, rank, the Williamson type, the elliptic, hyperbolic and twisting resonances completely determine the symplectic manifold $V/\Gamma$ with the linear $(S^1)^{r+\varkappa_e+\varkappa_h}$-action on $(V/\Gamma)^{\mathbb C}$, up to symplectomorphism. See also Remark \ref {rem:eigen:} about uniqueness of the Williamson type and the resonances.
\end{Definition}

Now we can formulate our result in the general real-analytic case, which is the symplectic normalization theorem for singular orbits of Hamiltonian torus actions:

\begin{Theorem} \label {thm:period:}
Suppose we are given real-analytic functions $I_1,\dots,I_r$, $J_1,\dots,J_{\varkappa_e+\varkappa_h}$ on a real-analytic symplectic manifold $(M,\Omega)$.
Suppose the functions $I_1^{\mathbb C},\dots,I_r^{\mathbb C}$, $J_1^{\mathbb C},\dots,J_{\varkappa_e}^{\mathbb C}$, $iJ_{\varkappa_e+1}^{\mathbb C},\dots,iJ_{\varkappa_e+\varkappa_h}^{\mathbb C}$ generate an effective Hamiltonian $(S^1)^{r+\varkappa_e+\varkappa_h}$-action on $(M^{\mathbb C},\Omega^{\mathbb C})$.
Suppose a point $m_0\in M$ is fixed under the $(S^1)^{\varkappa_e+\varkappa_h}$-subaction, and its orbit ${\cal O}_{m_0}$ is $r$-dimensional.
Then:

{\rm (a)} There exist an $(S^1)^{r+\varkappa_e+\varkappa_h}$-invariant neighbourhood $U$ of ${\cal O}_{m_0}$ in $M^{\mathbb C}$,
a finite group $\Gamma$, a linear $(S^1)^{r+\varkappa_e+\varkappa_h}$-action of rank $r$ on the symplectic manifold $(V/\Gamma)^{\mathbb C}$ given by (\ref {eq:mz:product}), (\ref {eq:mz:H:})--(\ref {eq:mz:F:}), and a real-analytic 
symplectomorphism $\phi$ from $U\cap M$ to $V/\Gamma$ having an action-preserving holomorphic extension to $U$ and 
sending the orbit ${\cal O}_{m_0}$ to the torus ${\cal O}=\{\lambda_s=x_j=y_j=0\}/\Gamma$.

In particular, dimensions of the torus $(S^1)^{r+\varkappa_e+\varkappa_h}$ and its subtori 
$(S^1)^{\varkappa_e}$ and $(S^1)^{\varkappa_h}$ satisfy the estimates $\varkappa_e\le k_e+k_f$, $\varkappa_h\le k_h+k_f$, and $r+\varkappa_e+\varkappa_h\le r+k_e+k_h+2k_f\le n$. If $r+\varkappa_e+\varkappa_h=n$ then the orbit ${\cal O}_{m_0}$ is nondegenerate.

{\rm (b)} The symplectic normalization 
(\ref {eq:mz:product}), (\ref {eq:mz:H:})--(\ref {eq:mz:F:}) for a torus action is persistent under analytic perturbations in the following sense. 
For any $k\in{\mathbb Z}_+$ and any neighbourhood $U'$ of the orbit ${\cal O}_{m_0}$ in $M^{\mathbb C}$ having a compact closure $\overline{U'}\subset U$, there exists $\varepsilon>0$ satisfying the following.
Suppose $\tilde F=(\tilde I_1,\dots,\tilde I_r,\tilde J_1,\dots,\tilde J_{\varkappa_e+\varkappa_h})$ and $\tilde\Omega$ are analytic (``perturbed'') momentum map and symplectic structure, whose holomorphic extensions to $M^{\mathbb C}$ are $\varepsilon$-close to $F^{\mathbb C}
=(I_1^{\mathbb C},\dots,I_r^{\mathbb C},J_1^{\mathbb C},\dots,J_{\varkappa_e+\varkappa_h}^{\mathbb C})$ and $\Omega^{\mathbb C}$, respectively, in $C^0-$norm.
Suppose the functions $\tilde I_1^{\mathbb C},\dots,\tilde I_r^{\mathbb C}$, $\tilde J_1^{\mathbb C},\dots,\tilde J_{\varkappa_e}^{\mathbb C}$, $i\tilde J_{\varkappa_e+1}^{\mathbb C},\dots,i\tilde J_{\varkappa_e+\varkappa_h}^{\mathbb C}$ generate a Hamiltonian $(S^1)^{r+\varkappa_e+\varkappa_h}$-action on $(\tilde M^{\mathbb C},\tilde\Omega^{\mathbb C})$, where $U\subset \tilde M^{\mathbb C}\subset M^{\mathbb C}$.
Then there exist an invariant neighbourhood $\tilde U\supset U'$, and an analytic 
symplectomorphism $\tilde\phi$ from $\tilde U\cap M$ to $V/\Gamma$,
whose holomorphic extension to $\tilde U$ is action-preserving and 
$O(\varepsilon)-$close to $\phi^{\mathbb C}$ in $C^k-$norm.
If the system depends on a local parameter (i.e.\ we have a local family of systems), 
moreover its holomorphic extension to $M^{\mathbb C}$ depends smoothly (resp., analytically) on that parameter, then
$\phi$ can also be chosen to depend smoothly (resp., analytically) on that parameter.
\end{Theorem}

\begin{Remark} \label {rem:eigen:}
(A) Consider a linear $(S^1)^{r+\varkappa_e+\varkappa_h}$-action on $(V/\Gamma)^{\mathbb C}$ (see Definition \ref {def:lin:model:}, Remark and Definition \ref {rem:period:}, Definition \ref {def:period:}).
What is a dynamical meaning of the elliptic resonances (\ref {eq:res:inf:e}), hyperbolic resonances (\ref {eq:res:inf:h}) and the twisting resonances (\ref {eq:res:mon})?
The eigenvalues of the linearized $\ell$-th infinitesimal generator of the 
%$\varkappa_e$-dimensional 
$\varkappa_e$-subtorus action on $(D^2)^{k_e+k_h+2k_f}$ are pure imaginary and in (\ref {eq:res:inf:e}) resonance, which is exactly the elliptic resonance assigned to this generator, $1\le \ell\le \varkappa_e$.
Further, the eigenvalues of the linearized $\ell$-th infinitesimal generator of the 
%$\varkappa_h$-dimensional 
$\varkappa_h$-subtorus action on $((D^2)^{k_e+k_h+2k_f})^{\mathbb C}$ are real and in (\ref {eq:res:inf:h}) resonance, which is exactly the hyperbolic resonance assigned to this generator, $\varkappa_e+1\le \ell\le \varkappa_e+\varkappa_h$.
Furthermore, the multipliers of a cycle $\gamma\in H_1({\cal O}) \approx p^{-1}(\Gamma)
\subset2\pi{\mathbb Q}^r\subset\R^r$ of the torus ${\cal O}=\{\lambda_s=x_j=y_j=0\}/\Gamma$ under the $r$-subtorus action on $(D^2)^{n-r}$ equal
$e^{\pm2\pi i q_{\psi,j}},e^{\pm2\pi i q_{\psi,j}}$ for $1\le j\le k_f$,
$e^{\pm2\pi i q_{\psi,j}}$ for $k_f+1\le j\le n-r-k_f$.
So, they are in $(q_{\psi,1},-q_{\psi,1},\dots,q_{\psi,k_f},-q_{\psi,k_f}$,
$q_{\psi,k_f+1},\dots,q_{\psi,n-r-k_f})\in({\mathbb Q}/{\mathbb Z})^{n-r}$ resonance, which is exactly the twisting resonance assigned to the cycle $\gamma$.
Here $p:\R^r\to(S^1)^r$ denotes the projection, $\psi := p(\gamma)\in\Gamma$.

(B) Clearly, the Williamson type is well defined, i.e.\ completely determined by the
(``hidden'') torus-symmetry of the orbit $\cal O$.
Are the elliptic and hyperbolic resonances also well defined (perhaps, up to some natural
transformations)? 
At first, we recall that the Hamiltonian action of the subtorus $(S^1)^{r+\varkappa_e}$ of the torus $(S^1)^{r+\varkappa_e+\varkappa_h}$ is generated by real-valued functions, while 
the Hamiltonian action of the subtorus $(S^1)^{\varkappa_h}$ is generated by imaginary-valued functions. Thus, the subtori $(S^1)^{r+\varkappa_e}$ and $(S^1)^{\varkappa_h}$ of the torus $(S^1)^{r+\varkappa_e+\varkappa_h}$ are well defined, i.e.\ they do not depend on the choice of generators of the $(S^1)^{r+\varkappa_e+\varkappa_h}$-action.
At second, the subtorus $(S^1)^{\varkappa_e}$ of the subtorus $(S^1)^{r+\varkappa_e}$ trivially acts on $\cal O$, moreover this subtorus is the isotropy subgroup of any point of $\cal O$ under the $(S^1)^{r+\varkappa_e}$-action, thus the subtorus $(S^1)^{\varkappa_e}$ is also well defined.
It follows that the elliptic resonances (\ref {eq:res:inf:e}) are well defined up to replacing the vectors $(p_{1\ell},p_{3\ell},\dots,p_{2k_f-1,\ell},\
p_{2k_f+1,\ell},\dots,p_{2k_f+k_e,\ell})$, $1\le\ell\le\varkappa_e$, by their linear
combinations with integer coefficients forming a nondegenerate
$\varkappa_e\times\varkappa_e$-matrix. Similarly, the hyperbolic resonances
(\ref {eq:res:inf:h}) are well defined up to replacing the
vectors $(p_{2\ell},p_{2\ell},\dots,p_{2k_f,\ell},\
p_{2k_f+k_e+1,\ell},\dots,p_{2k_f+k_e+k_h,\ell})$, $\varkappa_e+1\le\ell\le \varkappa_e+\varkappa_h$, by their
linear combinations with integer coefficients forming a nondegenerate
$\varkappa_h\times\varkappa_h$-matrix.

(C) Are the twisting resonances also well defined (perhaps, up to some natural
transformations)?
Suppose we allow ourselves to change the generators of the subtorus $(S^1)^r$ that acts locally-freely on our manifold $(V/\Gamma)^{\mathbb C}$.
At the same time, suppose that we fixed basic cycles $\gamma^1,\dots,\gamma^r\in H_1({\cal O})$ of the orbit ${\cal O}$, and we want that each cycle $\gamma^a$ has coordinates
$\gamma^a_s=\frac{2\pi}{N_a}\delta^a_s$, $1\le a,s\le r$ (see Remark and Definition \ref {rem:period:} (C)).
Thus, the above change is equivalent to replacing the variables $\lambda_s, \varphi_s, z_j$ with
$$
\hat\lambda_s := \nu_s \lambda_s + \sum\limits_{\ell=1}^{\varkappa_e} n_{\ell s} Q_\ell, \qquad
\hat\varphi_s := \frac1{\nu_s}\varphi_s, \qquad
\hat z_j := e^{-i  \langle (pn)_j , \psi \rangle}z_j
$$
for $1\le s\le r$ and $1\le j\le n-r$.
Here $\nu_s,n_{1 s},\dots,n_{\varkappa_e s}$ are coprime integers such that $\nu_s\ne 0$, and we denoted
$(pn)_{2j-1,s}=-(pn)_{2j,s}:=\sum\limits_{\ell=1}^{\varkappa_e} p_{2j-1,\ell}n_{\ell s}$ if $1\le j\le k_f$,
$(pn)_{js}:=\sum\limits_{\ell=1}^{\varkappa_e} p_{k_f+j,\ell}n_{\ell s}$ if $2k_f+1\le j\le 2k_f+k_e+k_h$,
$(pn)_{js}:=0$ if $2k_f+k_e+k_h+1\le j\le n-r$,
$(pn)_j:=((pn)_{j1},\dots,(pn)_{jr})$,
$\psi:=(\psi_1,\dots,\psi_r)$,
$\langle (pn)_j , \varphi \rangle := \sum\limits_{s=1}^r (pn)_{js}\varphi_s$.
This replacement will lead to the following replacements:
$\psi_a=(\psi^1_a,\dots,\psi^r_a) \to \hat\psi_a=(\hat\psi^1_a,\dots,\hat\psi^r_a):=(\frac{\psi^1_a}{\nu_1},\dots,\frac{\psi^r_a}{\nu_r})$, and
$m_j=(m_{j1},\dots,m_{jr}) \to \hat m_{js}$ where 
$\hat m_{js}:=\nu_s(m_{js}-(pn)_{2j-1,s})$ if $1\le j\le k_f$,
$\hat m_{js}:=\nu_s(m_{js}-(pn)_{k_f+j,s})$ if $k_f+1\le j\le n-r-k_f$.
Thus, the twisting resonances (\ref {eq:res:mon}) will be replaced with
$$
(\hat q_{\psi^a,1},-\hat q_{\psi^a,1},\dots,
\hat q_{\psi^a,k_f},-\hat q_{\psi^a,k_f}, \hat q_{\psi^a,k_f+1},\dots,
\hat q_{\psi^a,n-r-k_f})\in({\mathbb Q}/{\mathbb Z})^{n-r}, 
$$
$1\le a\le r$, where $\hat q_{\psi^a,j}
= q_{\psi^a,j} + \langle (pn)_{2j-1},\frac{\psi^a}{2\pi}\rangle
\mod 1\in{\mathbb Q}/{\mathbb Z}$ if $1\le j\le k_f$,
$\hat q_{\psi^a,j}= q_{\psi^a,j} + \langle (pn)_{k_f+j},\frac{\psi^a}{2\pi}\rangle
\mod 1\in{\mathbb Q}/{\mathbb Z}$ if $k_f+1\le j\le n-r-k_f$.
If we recall that $\psi^a_s=\frac{2\pi}{N_a}\delta^a_s$, $1\le a,s\le r$,
then $\hat q_{\psi^a,j} = q_{\psi^a,j} + (pn)_{2j-1,a}/N_a$ if $1\le j\le k_f$,
$\hat q_{\psi^a,j} = q_{\psi^a,j} + (pn)_{k_f+j,a}/N_a$ if $k_f+1\le j\le n-r-k_f$.
In other words, the twisting resonances (\ref {eq:res:mon}) are well defined up to adding any linear combinations of the ``extended'' elliptic resonances (\ref {eq:res:inf:e})
with rational coefficients forming a $r\times\varkappa_e$-matrix.
Here, by the $\ell$-th {\em extended elliptic resonance}, $1\le \ell\le \varkappa_e$, we mean the vector
$(p_{1,\ell},-p_{1,\ell},p_{3,\ell},-p_{3,\ell},\dots,p_{2k_f-1,\ell},-p_{2k_f-1,\ell},p_{2k_f+1,\ell},\dots,p_{2k_f+k_e,\ell},0,\dots,0)\in{\mathbb Z}^{n-r}$.
\end{Remark}

\begin{Ex} \label {exa:res:}
In Table below, we give hyperbolic analogues of the singularities (b) and (c) from Example \ref {exa:res}.
In detail, suppose that ${\cal O}$ is a compact $r$-dimensional orbit of the momentum map of a real-analytic integrable Hamiltonian system with $n$ degrees of freedom. 
Suppose the local
%this 
singularity at ${\cal O}$ has one of the following types:
\begin{itemize}
\item[\rm(b)] a hyperbolic integrable Hamiltonian Hopf bifurcation
\cite [\S 2]{ler00}, given by a momentum map $F=(H,I,J):D^1_{(\lambda)}\times S^1_{(\varphi)}\times D^4_{(x_1,y_1,x_2,y_2)}\to\R^3$ where $I=\lambda$, $J=x_1y_1+x_2y_2$, $H=x_1y_2+(x_2y_1)^2+\lambda x_2y_1$;
\item[\rm(c)] a normally-hyperbolic parabolic orbit with resonance $\ell/s$ ($s\ge5$), given by a momentum map $F=(H,I,J):(D^1_{(\lambda)}\times S^1_{(\varphi)}\times D^4_{(z_1,z_2)})/{\mathbb Z}_s\to\R^3$ where $I=\lambda$, $J=x_1y_1$,
$H=\re(z_2^s)+|z_2|^4+(I\pm J)|z_2|^2$, one denotes $z_j=x_j+iy_j$, 
a generator of ${\mathbb Z}_s$ acts on $D^1\times S^1\times D^4$ by the transformation $(\lambda,\varphi,z_1,z_2)\mapsto(\lambda,\varphi+2\pi/s,z_1,e^{2\pi i\ell/s}z_2)$, $s,\ell\in\mathbb Z$, $0\le\ell<s$, $(s,\ell)=1$.
\end{itemize}
One shows from Theorem \ref {thm:kudr:} that this singularity possesses a ``hidden'' $(S^1)^{r+\varkappa_e+\varkappa_h}$-symmetry, i.e.\ an effective (not locally-free) $F^{\mathbb C}$-preserving Hamiltonian $(S^1)^{r+\varkappa_e+\varkappa_h}$-action, with $r+\varkappa_e+\varkappa_h=n-1$. Moreover, this action is persistent under real-analytic integrable perturbations.
By Theorem \ref {thm:period:}, the latter action is symplectomorphic to a linear $(S^1)^{r+\varkappa_e+\varkappa_h}$-action, that is also persistent under real-analytic integrable perturbations. In Table, we show the Williamson type $(k_e,k_h,k_f)$, the elliptic, hyperbolic and twisting resonances (Definition \ref {def:period:}) of this linear $(S^1)^{r+\varkappa_e+\varkappa_h}$-action, for each type of the singularity at ${\cal O}$ from above.
\end{Ex}

\begin{tabular}{|c|c|c|c|c|c|c|c|c|c|}
\hline
Case &$n$& \multicolumn{3}{c|}{Subtori dim's} & Williamson & \multicolumn{3}{c|}{Resonances} 
\\
\cline{3-5}\cline{7-9}
&& \ $r$ \ &  $\varkappa_e$ & $\varkappa_h$ & type $(k_e,k_h,k_f)$ & elliptic & hyperbolic & twisting \\
\hline
\hline
(b) & 3 & 1 &0&1& $(0,2,0)$ &no& $(1:1)\in{\mathbb Q}P^{1}$ & $(0,0)\in({\mathbb Q}/{\mathbb Z})^2$ \\
\hline
(c) & 3 & 1 &0&1& $(1,1,0)$ &no& $1\in{\mathbb Q}P^{0}$ & $(\ell/s\mod1,0)\in({\mathbb Q}/{\mathbb Z})^2$ \\
\hline
\end{tabular}

\section {Application to structural stability of singularities} \label {subsec:stab}

In this section, we solve the questions (3b)--(3d) from Introduction, for parabolic orbits with resonances.

Two singularities will be called {\em equivalent} if there exists a fiberwise homeomorphism of fibration germs at these singularities.

Our central object will be {\em structurally stable} singularities (Definition \ref{def:stab} below), which are those singularities whose equivalence classes are {\em open} in the topology described below.
Such singularities occur
%are met 
in typical integrable systems, so they cannot disappear after small integrable perturbations.
We will assume that the manifold $M$, the symplectic structure $\Omega$ and the momentum map $F$ are {\em real-analytic}. Notice that, due to the Cauchy theorem or the Weierstrass theorem \cite[Ch.~I, \S2, Theorem 8]{shabat}, all compact-open $C^k-$topologies on the space of holomorphic pairs $(\Omega^{\mathbb C},F^{\mathbb C})$ on $M^{\mathbb C}$ coincide for all $k\in{\mathbb Z}_+$.
Here $M^{\mathbb C}$ denotes a small open complexification of $M$, while $\Omega^{\mathbb C},F^{\mathbb C}$ are holomorphic extensions of $\Omega,F$ to $M^{\mathbb C}$. Thus, below we can take $k=0$.

\begin{Definition} \label {def:stab}
A singularity at $\cal O$ of a singular Liouville fibration $(M,\Omega,F)$ will be called {\em structurally stable under real-analytic integrable perturbations}, or simply {\em structurally stable} if $\cal O$ has a neighbourhood $U_0$ such that, for any smaller neighbourhood $U_1$ with a compact closure $\overline{U_1}\subset U_0$, there exist $\varepsilon>0$ and a (small) open complexification $U_0^{\mathbb C}$ of $U_0$
satisfying the following.
For any real-analytic integrable perturbation $(U_0,\tilde\Omega,\tilde F)$ of $(U_0,\Omega|_{U_0},F|_{U_0})$ such that $\|\tilde\Omega^{\mathbb C}-\Omega^{\mathbb C}\|_{C^{k}}+\|\tilde F^{\mathbb C}-F^{\mathbb C}\|_{C^k}<\varepsilon$, the singular Liouville fibrations $(U,\Omega|_U,F|_U)$ and $(\tilde U,\tilde\Omega|_{\tilde U},\tilde F|_{\tilde U})$ are equivalent (i.e., fiberwise homeomorphic) for some neighbourhoods $U,\tilde U\subset U_0$ each of which contains $U_1$.
In a similar way, one defines {\em structural stability under integrable perturbations of some class}, e.g.\ the classes of {\em $C^\infty$ perturbations}, {\em symmetry-preserving perturbations} etc.
\end{Definition}

\begin{Ex} \label {exa:stab}
(A) Consider a nondegenerate singular orbit $\cal O$ of a real-analytic integrable system.
One can show, using the Vey-Eliasson theorem \cite {vey} (cf. \cite {eli} for $C^\infty$ case), that the singularity at every point $m\in\cal O$ is structurally stable under integrable perturbations (both in real-analytic and $C^\infty$ cases).
Let us show that the singularity at the orbit $\cal O$ is structurally stable (under real-analytic integrable perturbations), provided that $\cal O$ is compact.
Due to the result by Ito \cite {ito91}, our singular Lagrangian fibration on a neighbourhood $U$ of $\cal O$ can be defined by a momentum map $F=(I_1,\dots,I_r,J_1,\dots,J_{n-r})$ having a standard form (called a symplectic normalization, see \S\ref {subsec:mult} with $r+\varkappa_e+\varkappa_h=n$) w.r.t.\ some real-analytic coordinate system. In particular, on some small open complexification $U^{\mathbb C}$ of $U$, we have a $F$-preserving linear $(S^1)^n$-action generated by $I_1,\dots,I_r,J_1,\dots,J_{\varkappa_e}$, $iJ_{\varkappa_e+1},\dots,iJ_{n-r}$ (Definitions \ref {def:lin:model}, \ref {def:lin:model:}).
One directly checks (see e.g.\ Lemma \ref {lem:ii} and Example \ref {exa:nondeg}) that, for each $S^1$-subaction, there exists a point $m_1\in (F^{\mathbb C})^{-1}(F({\cal O}))$ satisfying 
the conditions (i)--(iii) of Theorem \ref {thm:kudr:} (a). Hence, by 
Theorem \ref {thm:kudr:} (b), for any ``perturbed'' real-analytic integrable system $(U,\widetilde\Omega,\widetilde F)$, there exists a $\widetilde F$-preserving ``perturbed'' $(S^1)^n$-action generated by $\widetilde I_1,\dots,\widetilde I_r,\widetilde J_1,\dots,\widetilde J_{\varkappa_e}$ and $i\widetilde J_{\varkappa_e+1},\dots,i\widetilde J_{n-r}$, for some real-analytic functions $\widetilde I_s,\widetilde J_j$ close to $I_s,J_j$.
By Theorem \ref {thm:period:} (b), the latter ``perturbed'' $(S^1)^n$-action is also linear, moreover there exists an $(S^1)^n$-action-preserving symplectomorphism (close to the identity) between the ``unperturbed'' and the ``perturbed'' fibrations.

(B) We conjecture that all degenerate local singularities from Examples \ref {exa:res} and \ref {exa:res:} are structurally stable. We expect that, similarly to (A), one can derive this from Theorems \ref{thm:kudr:} and \ref {thm:period:}.
In Proposition \ref {pro:kal} below, we do this for parabolic orbits with resonances, that were briefly described in Example \ref {exa:res} (a).
\end{Ex}

\subsection {Structural stability and preliminary normal form for Kalashnikov's parabolic orbits with resonances} \label {subsec:kal}

Consider integrable systems with 2 degrees of freedom.
Such a system is defined by a pair $F=(H,K)$ of Poisson-commuting functions on a symplectic 4-manifold $(M^4,\Omega)$.

An important property of parabolic orbits is their structural stability under small integrable perturbations (see Lerman and Umanskii \cite{ler87}).
This is one of the reasons why such orbits can be observed in many examples of integrable Hamiltonian systems: Kovalevskaya top \cite{bol:fom:ric}, other integrable cases in rigid body dynamics including Steklov case, Clebsch case, Goryachev--Chaplygin--Sretenskii case, Zhukovskii case, Rubanovskii case and Manakov top on $so(4)$ \cite{BF}, as well as systems invariant w.r.t. rotations \cite{kant, koz:osh, kudr:osh}, see also examples discussed in \cite{efs:gia}, \cite{dul:iva2}. Unlike nondegenerate singularities, however, in the literature on topology and singularities of integrable systems there are only few papers devoted to degenerate singularites including parabolic ones. We refer, first of all, to  the following six ---  L. Lerman, Ya. Umanskii \cite{ler:uma2}, V. Kalashnikov \cite{kal}, N. T. Zung \cite{zung00}, H. Dullin, A. Ivanov \cite{dul:iva2}, K. Efstathiou, A. Giacobbe \cite{efs:gia} and Y. Colin de Verdi\`ere \cite{CdV} --- which we consider to be very important in the context of general classification programme for bifurcations occurring in integrable systems.
Rigidity of Hamiltonian group actions was applied for obtaining structural stability results for nondegenerate singularities \cite{Mir14} and $S^1$-invariant degenerate singularities \cite{MM20} of integrable systems.

Parabolic orbits with resonances were discovered by Kalashnikov \cite {kal}, who proved that they form a complete list of typical degenerate rank-1 singularities, and are structurally stable under $S^1$-symmetry-preserving integrable perturbations (in $C^\infty$ case). As we will show in Proposition \ref {pro:kal} below, they are structurally stable (in real-analytic case) even in the following stronger sense: structurally stable under all integrable perturbations, not necessarily preserving the $S^1$-action.
Parabolic singularity with $\frac12$ resonance and a plus sign (known as ``elliptic period-doubling bifurcation'', cf.\ (\ref {eq:kal:f}) with $s=2$ and a plus sign) can be observed in the Sretenskii system; its ${\mathbb Z}_2$-symmetric 2-fold cover (known as ``elliptic pitchfork'') appear
%are met 
in the problems by Kovalevskaya, Steklov, Neumann, Clebsch \cite{bol:fom:ric}. 
Parabolic singularity with $\frac12$ resonance and a minus sign (known as ``hyperbolic period-doubling bifurcation'', cf.\ (\ref {eq:kal:f}) with $s=2$ and a minus sign) is another local singularity.
The corresponding semilocal singularity appears
%is met 
in the problems by Kovalevskaya and Sretenskii; 
its two topologically different ${\mathbb Z}_2$-symmetric 2-fold covers (known as ``hyperbolic pitchforks'') appear
%are met 
in the Kovalevskaya problem \cite{bol:fom:ric}.

It would be interesting to find examples of mechanical integrable systems having a parabolic orbit with ``higher order'' resonances (i.e., resonances different from $0$ and $1/2 \mod1$). More generally, to find examples of mechanical integrable systems having a rank-$r$ local singularity with ``higher order'' twisting resonances (\ref {eq:res:mon}), i.e.\ resonances $(q_{a,1},\dots,q_{a,n-r})\in({\mathbb Q}/{\mathbb Z})^{n-r}$, $1\le a\le r$, where at least some $q_{a,j}$ is different from $0$ and $1/2\mod1$ (for any choice of generators of an $(S^1)^r$-subaction, see Remark \ref {rem:eigen:} (C)).

It is known that from the smooth point of view, all parabolic orbits without resonances are equivalent, i.e., any two parabolic orbits admit fiberwise diffeomorphic neighborhoods 
(Lerman-Umanskii \cite{ler87, ler:uma2}, Kalashnikov \cite{kal}, the author and Martynchuk 
\cite{kud:mar21}). The same is true for cuspidal tori \cite{efs:gia, kud:mar21}. A symplectic classification of parabolic orbits in real-analytic case is studied in \cite{bgk}.

Below we describe the structure of the singular Lagrangian fibration in a neighborhood of a parabolic orbit with resonance in real-analytic case. As we are mostly interested in this fibration (rather than specific commuting functions $H$ and $K$), we allow ourselves to replace 
$H$ and $K$ with $H_1 = H_1(H,K)$ and $K_1 = K_1(H,K)$ where $\frac{\partial (H_1,K_1)}{\partial (H,K)}\ne 0$.

A model (called ``preliminary normal form'') for a parabolic singularity with resonance is as follows. Denote by $D^1$ a small interval centred at $0$ with coordinate $\lambda$, by $S^1$ a circle with a standard periodic angle coordinate $\varphi \mod 2\pi$, and by $D^2$ a small open disk centred at the origin with coordinates $x,y$. Consider the manifold 
$$
V=D^1\times S^1 \times D^2.
$$
Let $s$ be a positive integer (called the {\em resonance order}).
Consider the free action $\rho$ of the group $\Gamma={\mathbb Z}_s$ on $V$ generated by the map
\begin{equation} \label {eq:kal:action}
(\lambda,\varphi,z)\mapsto(\lambda,\varphi+\frac{2\pi}s,e^{2\pi i\frac\ell s}z), \qquad
 \mbox{where} \quad z=x+iy,
\end{equation}
$\ell\in\mathbb Z$ is an integer coprime with $s$ (we may assume that $0\le\ell<s$).
Consider the following 1-parameter family of ${\mathbb Z}_s$-invariant functions $f_{a,\lambda,s}(z)$ on $D^2$ with parameter $\lambda\in\R$ (where $a\in\R\setminus\{1,-1\}$ is an additional parameter appearing for $s=4$ only):
\begin{equation} \label {eq:kal:f}
\begin{array} {ll}
f_{\lambda,1}(x,y)=x^2+y^3+\lambda y, & \\
f_{\lambda,2}(x,y)=x^2\pm y^4+\lambda y^2, & \\
f_{\lambda,3}(z)=\re(z^3)+\lambda|z|^2, & \\
f_{a,\lambda,4}(z)=\re(z^4)+a|z|^4+\lambda|z|^2, & \qquad a^2\ne1, \\
f_{\lambda,s}(z)=\re(z^s)+|z|^4+\lambda|z|^2, & \qquad s\ge5.\
\end{array}
\end{equation}

We endow $V$ with the ${\mathbb Z}_s$-invariant symplectic 2-form
\begin{equation} \label {eq:kal:omega}
\Omega = \ddd\alpha(\lambda) \wedge \ddd\varphi + \pi^*\omega ,
\end{equation}
where 
$\pi:\R^3 \times S^1\to\R^3$ is the projection;
$\alpha(\lambda)$ is any real-analytic function such that $\alpha'(\lambda)>0$, moreover $\alpha(\lambda)\equiv\lambda$ if $s<5$;
$\omega$ is an arbitrary ${\mathbb Z}_s$-invariant closed 2-form on $D^2\times D^1$ such that $\omega\wedge \ddd\lambda$ nowhere vanish.
The latter condition on $\omega$ is equivalent to the following: 
\begin{equation} \label {eq:kal:PQR}
\omega=R(x,y,\lambda)\ddd x\wedge \ddd y+\ddd\lambda\wedge(Q(x,y,\lambda)\ddd x-P(x,y,\lambda)\ddd y),
\end{equation}
for some ${\mathbb Z}_s$-invariant real-analytic divergence-free vector field $(P,Q,R)$ on $D^2\times D^1$ having a nowhere vanishing ``vertical component'' $R(x,y,\lambda)$.
In a simplest case, we have $R(x,y,\lambda)\equiv1$ and $\omega=\ddd x\wedge \ddd y$. However, in general, we cannot assume that the coordinates $\alpha(\lambda), \varphi,x,y$ are canonical, since $\omega$ can depend on $\lambda$ in an essential way.

Now we can form the quotient symplectic manifold $V/{\mathbb Z}_s$, with an integrable system on it given by the following two functions:
\begin{equation} \label {eq:kal:HF}
H = f_{a(\lambda),\lambda,s}(x,y) \quad\mbox{and} \quad K = \lambda,
\end{equation}
where $a(\lambda)$ is any real-analytic function such that $a(\lambda)\ne\pm1$ (such a function appears if $s=4$ only).
The functions (\ref {eq:kal:HF}) Poisson-commute w.r.t.\ the symplectic 2-form (\ref {eq:kal:omega}).

Due to a result by V. Kalashnikov \cite{kal}, the curve $\gamma(t)=(0, 0, 0, t)$ is a {\em parabolic orbit with $\ell/s$ resonance} of an integrable Hamiltonian system defined by commuting functions $H$ and $K$ on $(V/{\mathbb Z}_s,\Omega)$.
A formal definition of a parabolic orbit with $\ell/s$ resonance is given in \cite{kal} (see also \cite{bgk} for $s=1$).

In according to the definitions \ref {def:lin:model} and \ref {def:period}, $\gamma(t)$ is a rank-1 orbit admitting a linear $S^1$-action (see Example \ref {exa:res} (a) for properties of this $S^1$-action). 

\begin{Proposition} \label {pro:kal}
Consider the rank-1 orbit ${\cal O}=\{\gamma(t)\}$ (called a parabolic orbit with $\ell/s$ resonance) of the real-analytic singular Lagrangian fibration $(V/{\mathbb Z}_s,\Omega,F)$, whose symplectic 2-form $\Omega$ and momentum map $F=(H,K)$ are in a ``preliminary normal form'' (\ref {eq:kal:f})--(\ref {eq:kal:HF}) on $V/{\mathbb Z}_s$, and the ${\mathbb Z}_s$-action on $V$ has the form (\ref {eq:kal:action}). Then:

{\rm (a)} The singularity at the orbit $\cal O$ is structurally stable under real-analytic integrable perturbations.

{\rm (b)} Moreover, the preliminary normal form (\ref {eq:kal:action})--(\ref {eq:kal:HF}) is persistent under real-analytic integrable perturbations (up to some left-right change of variables) in the following sense. For any $k\in{\mathbb Z}_+$ and any neighbourhood $U'$ of $\cal O$ having a compact closure $\overline{U'}\subset V/{\mathbb Z}_s$, there exists $\varepsilon>0$ satisfying the following.
For any (``perturbed'') real-analytic integrable system $(V/{\mathbb Z}_s,\widetilde\Omega,\widetilde F)$ whose holomorphic extension to $(V/{\mathbb Z}_s)^{\mathbb C}$ is $\varepsilon$-close to $((V/{\mathbb Z}_s)^{\mathbb C},\Omega^{\mathbb C},F^{\mathbb C})$ in $C^0-$norm, with $\widetilde F=(\widetilde H,\widetilde K)$, one can choose a neighbourhood $\widetilde U\supset U'$ in $V/{\mathbb Z}_s$, real-analytic coordinate changes $\widetilde \phi:\widetilde U\to V/{\mathbb Z}_s$ and
$\widetilde\chi:(\tilde H,\tilde K)\mapsto(\tilde H_1(\tilde H,\tilde K),\tilde K_1(\tilde H,\tilde K))$ 
that are $O(\varepsilon)-$close to the identities in $C^k-$norm and bring the ``perturbed'' singular Lagrangian fibration on $\widetilde U$ to a preliminary normal form
\begin{equation} \label {eq:kal:can}
(\widetilde\phi(\widetilde U),\ 
(\widetilde\phi^{-1})^*\widetilde\Omega=\ddd\widetilde\alpha(\lambda)\wedge\ddd\varphi + \pi^*\widetilde\omega,\ 
\widetilde\chi\circ\widetilde F\circ\widetilde\phi^{-1}=(f_{\lambda,\tilde a(\lambda),s}(z),\lambda)),
\end{equation}
for some analytic (``perturbed'') ${\mathbb Z}_s$-invariant closed 2-form $\widetilde\omega$ and functions $\tilde a(\lambda)$, $\widetilde\alpha(\lambda)$ that are $O(\varepsilon)-$close to the 2-form $\omega$ and the functions $a(\lambda)$, $\alpha(\lambda)$ in $C^k-$norm, $\widetilde\alpha(\lambda)\equiv\lambda$ if $s<5$ (the function $\tilde a(\lambda)$ appears only for resonance order $s=4$).

In particular, the ``unperturbed'' and the ``perturbed'' singular Lagrangian fibrations are locally (i.e., on some neighbourhoods $U,\widetilde U\supset U'$ of $\cal O$) fiberwise homeomorphic. Moreover, they are fiberwise diffeomorphic if $s\ne4$.
\end{Proposition}

\begin{Comment}
We stress that the 2-forms $\omega$ and $\widetilde\omega$, that appear in the preliminary normal form for the ``unperturbed'' and the ``perturbed'' fibrations, may be different.
Similarly, the functions $\alpha(\lambda)$ and $\widetilde\alpha(\lambda)$ may be different (if $s\ge5$).
However they affect the symplectic structure only, so they do not affect the topology of our singular Lagrangian fibration.
Similarly, the functions $a(\lambda)$ and $\widetilde a(\lambda)$ (if $s=4$) may be different. However they affect only the smooth structure (rather than the topology) of our singular Lagrangian fibrations, provided that $|a(\lambda)|-1$ and $|\widetilde a(\lambda)|-1$ have the same sign. If the signs are different, the singular Lagrangian fibrations corresponding to these two functions have different topology.
\end{Comment}

\begin{proof}
{\em Step 1.} Consider the linear $S^1$-action on $V$ by shifts along the $\varphi$-axis. Observe that it is a Hamiltonian $S^1$-action generated by 
the function $I(H,K)=\alpha(\lambda)$, w.r.t.\ the symplectic structure $\Omega$ in (\ref {eq:kal:omega}). 
Therefore, a function on $V$ Poisson commutes with $K=\lambda$ if and only if it is $S^1$-invariant, i.e.\ does not depend on $\varphi$.

{\em Step 2.} Suppose we are given a ``perturbed'' real-analytic integrable system $(U,\widetilde\Omega,\widetilde F)$, where $\widetilde F=(\widetilde H,\widetilde K)$.
Let us show that it admits a real-analytic $2\pi$-periodic first integral $\tilde I(\widetilde H,\widetilde K)$ on $\tilde U\supset U'$, where $\tilde I(z_1,z_2)$ is close to $I(z_1,z_2)=\alpha(z_2)$. 
One directly checks that there exists a point $m_1\in (F^{\mathbb C})^{-1}(F({\cal O}))$ satisfying either the conditions (i) and (ii) of Theorem \ref {thm:kudr} or the conditions (i)--(iii) of Theorem \ref {thm:kudr:}, moreover the paths $\gamma_a$ and $\overline{\gamma_a}$ in (iii) are homological to each other in the fiber.
Hence, by Theorem \ref {thm:kudr} (b) or Theorem \ref {thm:kudr:} (b), the ``perturbed'' system has a $2\pi$-periodic first integral $\tilde I(\tilde H,\tilde K)$, for some real-analytic function $\tilde I(z_1,z_2)$ close to $I(z_1,z_2)=\alpha(z_2)$. 
Define the following change of first integrals:
$$
\chi_1:(\tilde H,\tilde K)\mapsto(\tilde H,\tilde I(\tilde H,\tilde K)).
$$

{\em Step 3.} 
Consider the 1-parameter family of (``unperturbed'') symplectic 2-forms 
$\omega_\lambda:=R(x,y,\lambda)\ddd x\wedge \ddd y$ with real parameter $\lambda$.
By assumption, this 2-form is invariant under the action $\rho'$ of the group $\Gamma:={\mathbb Z}_s$ on $D^2$, where $\rho'(1\mod s)(z)=e^{2\pi i/s}z$.

Due to Theorem \ref {thm:period} (b), this $\Gamma$-action $\rho'$ is symplectically normalizable, i.e.\ there exists a real-analytic $\Gamma$-equivariant change of variables $\phi_\lambda:(x,y)\mapsto(x_1(x,y,\lambda),y_1(x,y,\lambda))$ that analytically depends on $\lambda$ and brings the symplectic form $\omega_\lambda$ to the standard form $\ddd x_1\wedge\ddd y_1$, so $\omega_\lambda=\phi_\lambda^*\ddd x_1\wedge\ddd y_1$.

Define the 
($\Gamma$-equivariant and $S^1$-equivariant) 
change $\phi_0:(\lambda,\varphi,x,y)\mapsto(\lambda_1=\alpha(\lambda),\varphi,x_1,y_1)$. Then $\Omega=\phi_0^* (\ddd \lambda_1\wedge(\ddd\varphi + \beta) + \ddd x_1\wedge\ddd y_1)$ for some 1-form $\beta$. Since $\ddd\Omega=0$, it follows that $\ddd \lambda_1\wedge\beta=\ddd \lambda_1\wedge\ddd g$ for some real-analytic function $g=g(x,y,\lambda)$, that is $\Gamma$-invariant and $S^1$-invariant. Define $\varphi_1=\varphi+g(x,y,\lambda)$.
Then, with respect to the coordinate system $(\lambda_1,\varphi_1,x_1,y_1)$, the symplectic 2-form $\Omega$ has the standard form, the $S^1$-action is linear, and the free action $\rho$ of $\Gamma$ is also linear.

Clearly, the change 
$\phi_1:(\lambda,\varphi,x,y)\mapsto(\lambda_1,\varphi_1,x_1,y_1)$ is $(\Gamma\times S^1)$-equivariant.

{\em Step 4.} By Step 3, the (``unperturbed'') $S^1$-action is linear with respect to the coordinate system $(\lambda_1,\varphi_1,x_1,y_1)$ on $V$.
By Step 2, there exists a ``perturbed'' Hamiltonian $S^1$-action generated by a function $\tilde I(\tilde H,\tilde K)$ close to $I(H,K)=\alpha(K)$.
Hence, by Theorem \ref {thm:period} (b), the latter ``perturbed'' $S^1$-action is linear too, moreover there exists a ``perturbed'' coordinate change 
$\widetilde\phi_1:(\lambda,\varphi,x,y)\mapsto(\tilde\lambda_1,\tilde\varphi_1,\tilde x_1,\tilde y_1)$ such that, with respect to the ``perturbed'' coordinate system $(\tilde\lambda_1,\tilde\varphi_1,\tilde x_1,\tilde y_1)$, the ``perturbed'' symplectic 2-form $\widetilde\Omega$ has the standard form, moreover the ``perturbed'' Hamiltonian $S^1$-action generated by $\tilde I(\tilde H,\tilde K)$ and the ``perturbed'' free $\Gamma$-action $\widetilde\rho$ are linear.
This implies that $\widetilde I\circ(\widetilde H,\widetilde K)\circ\widetilde\phi_1^{-1}=\tilde\lambda_1+\const$. Without loss of generality, the latter constant is 0.

How do $\phi_1$ and $\widetilde\phi_1$ transform the Hamilton functions $H$ and $\widetilde H$ (respectively)?

The ``unperturbed'' coordinate change $\phi_1$ brings the ``unperturbed'' function $H=f_{a(\lambda),\lambda,s}(z)$ to a function $H\circ\phi_1^{-1}=:G=g_{\lambda_1}(z_1)$.
The ``perturbed'' coordinate change $\widetilde\phi_1$ brings the ``perturbed'' function $\widetilde H$ to a function $\widetilde H\circ\widetilde\phi_1^{-1}=:\widetilde G=\widetilde g_{\tilde\lambda_1}(\tilde z_1)$.
Indeed, both functions $G$ and $\widetilde G$ are $S^1$-invariant (since $H$ and $I(H,K)$ Poisson commute, and $\widetilde H$ and $\widetilde I(\widetilde H,\widetilde K)$ Poisson commute), hence they don't dependent on $\varphi_1$ and $\tilde\varphi_1$, respectively. 

The functions $G$ and $\widetilde G$ are invariant under the free $\Gamma$-action $\rho$, since $\phi_1$ and $\widetilde\phi_1$ are $\Gamma$-equivariant, 
and $H$ and $\widetilde H$ are invariant under the free $\Gamma$-actions $\rho$ and $\widetilde\rho$, respectively.

{\em Step 5.} Let us compose these two functions $G=g_{\lambda_1}(z_1)$ and $\widetilde G=\widetilde g_{\tilde\lambda_1}(\tilde z_1)$ with the ($\Gamma$-equivariant and $S^1$-equivariant, see Step 3) coordinate change $\phi_1$. We will obtain two functions on $V$:
$$
G\circ\phi_1=H=f_{a(\lambda),\lambda,s}(z)
\quad \mbox{and} \quad
\widetilde G\circ\phi_1=\widetilde H\circ(\widetilde\phi_1^{-1}\circ\phi_1)=:\widetilde f_\lambda(z),
$$
that both are also $\Gamma$-invariant and $S^1$-invariant (since $G$ and $\widetilde G$ are, by Step 4). 
We can and will regard the ``perturbed'' function $\widetilde f_\lambda(z)$ as a 1-parameter family of functions in variables $z=(x,y)$ with parameter $\lambda$.

The diffeomorphism $\widetilde\phi_1^{-1}\circ\phi_1$ is $\Gamma$-equivariant, since both $\phi_1$ and $\widetilde\phi_1$ are.
Since $(\phi_1^{-1})^*\Omega$ is standard (by Step 3) and $(\widetilde\phi_1^{-1})^*\widetilde\Omega$ is standard (by Step 4), they coincide. Hence 
\begin{equation} \label {eq:Omega}
(\widetilde\phi_1^{-1}\circ\phi_1)^*\widetilde\Omega=\Omega.
\end{equation}
By Step 4 and Step 3,
\begin{equation} \label {eq:tildeI}
\widetilde I\circ(\widetilde H,\widetilde K)\circ \widetilde\phi_1^{-1} \circ \phi_1
= \tilde\lambda_1\circ\phi_1 = \alpha(\lambda).
\end{equation}

{\em Step 6.} Define the 1-parameter family of functions on $D^2$ with parameter $\lambda$:
\begin{equation} 
\widehat f_{a,\lambda,s}(z):= \left\{
\begin{array} {ll}
f_{\lambda,s}(z) , & \qquad 1\le s\le 3, \\
f_{a,\lambda,4}(z), & \qquad a^2\ne1,\ s=4, \\
\re(z^s)+a|z|^4+\lambda|z|^2 , & \qquad a>0,\ s\ge5,
\end{array}
\right.
\end{equation}
cf.\ (\ref {eq:kal:f}), where $a$ is an additional real parameter which is inessential if $s\le3$ (in particular, $f_{\lambda,s}(z)=\widehat f_{1,\lambda,s}(z)$ for $s\ge5$).

It follows from a result by Wassermann \cite {wass} that there exist a $\Gamma$-equivariant change of variables $(\lambda,z)\to(\lambda,\widehat z(z,\lambda))$ close to the identity and smooth (real-analytic in our case) functions $\widetilde a(\lambda),\widetilde b(\lambda),\widetilde c(\lambda)$ close to $a(\lambda),\lambda,0$, respectively, such that 
$$
\widetilde f_\lambda(z)
= \widehat f_{\widetilde a(\lambda),\widetilde b(\lambda),s}(\widehat z(z,\lambda))+\widetilde c(\lambda).
$$

Consider the change of variables
$$
\widehat\phi:(\lambda,\varphi,z)\mapsto(\lambda,\varphi,\widehat z).
$$
By above, it is close to the identity and transforms our ``perturbed'' functions $\widetilde f_\lambda(z)=\widetilde H\circ (\widetilde\phi_1^{-1}\circ\phi_1)$ and 
(\ref {eq:tildeI}) to 
\begin{equation} \label {eq:f:right:2D}
\widetilde H\circ \widetilde\phi_1^{-1}\circ\phi_1 \circ \widehat\phi^{-1}
=\widehat f_{\widetilde a(\lambda),\widetilde b(\lambda),s}(\widehat z)+\widetilde c(\lambda) \quad \mbox{and} \quad
\widetilde I\circ(\widetilde H,\widetilde K)\circ \widetilde\phi_1^{-1} \circ \phi_1 \circ \widehat\phi^{-1} = \alpha(\lambda).
\end{equation}
Clearly, $\widehat\phi$ is $\Gamma$-equivariant and brings the symplectic structure (\ref {eq:Omega}) to a form
\begin{equation} \label {eq:Omega:}
(\widetilde\phi_1^{-1}\circ\phi_1 \circ \widehat\phi^{-1})^*\widetilde\Omega
= (\widehat\phi^{-1})^* \Omega = d\alpha(\lambda) \wedge d\varphi + \pi^*\widehat\omega 
\end{equation}
similar to (\ref {eq:kal:omega}), with some $\Gamma$-invariant 2-form $\widehat\omega$ close to $\omega$.

{\em Step 7.} 
Recall that we are mostly interested in our singular Lagrangian fibration (rather than specific commuting functions $\widetilde H$ and $\widetilde K$), so we allow ourselves to replace the ``perturbed'' momentum map $\widetilde F=(\widetilde H,\widetilde K)$ with a composition of it from the left with some diffeomorphism.
Our goal is to show that some of the functions $\widetilde a(\lambda),\widetilde b(\lambda),\widetilde c(\lambda)$ in (\ref {eq:f:right:2D}) can be simplified by an appropriate $\Gamma$-equivariant change of variables close to the identity, both in the source and in the target (so-called left-right change of variables), perhaps at the expense of spoiling the function $\alpha(\lambda)$ in (\ref {eq:Omega:}) if $s\ge5$.

Firstly, by a change of first integrals
$\widehat\chi:(\widetilde H,\widetilde I)\mapsto
(\widetilde H-\widetilde c\circ\alpha^{-1}(\widetilde I),\alpha^{-1}(\widetilde I))$,
we can kill the function $\widetilde c(\lambda)$ and reduce $\widetilde I$ to $\lambda$.

Secondly, the function $\widetilde b(\lambda)$ can be reduced to $\lambda$ (at the expense of replacing $\widetilde a(\lambda)$ with $\widehat a(\lambda)=\widetilde a(\lambda)b_1(\lambda)^{(4-s)/(s-2)}$ if $s\ge5$, but preserving $\lambda$ and $\widetilde I$) 
via the following $\Gamma$-equivariant change of variables and first integrals:
$$
\phi_2\times\chi_2:
(\lambda,\varphi,\widehat x,\widehat y;\widetilde H,\widetilde I)\mapsto
(\lambda,\varphi,u(\lambda)\widehat x,v(\lambda)\widehat y;w(\widetilde I)\widetilde H, \widetilde I),
$$ 
where
\begin{itemize}
\item $u(\lambda)=b_1(\lambda)^{-1/2-s/4}$, $v(\lambda)=b_1(\lambda)^{-1/2}$, 
$w=b_1(\lambda)^{-1-s/2}$ if $s=1,2$;
\item $u(\lambda)=v(\lambda)=b_1(\lambda)^{-1/(s-2)}$, $w=b_1(\lambda)^{-s/(s-2)}$ if $s\ge3$,
\end{itemize}
and $b_1(\lambda)>0$ is a continuous function determined by the condition $\widetilde b(\lambda)=\lambda b_1(\lambda)$.

Finally, if $s\ge5$, we can reduce the function $\widehat a(\lambda)>0$ to the constant $1$ (at the expense of replacing $\widetilde I$ with a composition of $\widetilde I$ from the left with some diffeomorphism) by the following $\Gamma$-equivariant change of variables and first integrals:
$$
\phi_3\times\chi_3:
(\lambda,\varphi,\widehat z;\widetilde H,\widetilde I) \mapsto 
(t(\lambda)^{s-2}\lambda,\varphi,t(\lambda)\widehat z;t(\widetilde I)^s\widetilde H,t(\widetilde I)^{s-2}\widetilde I), \quad 
t(\lambda):=\widehat a(\lambda)^{1/(4-s)}>0.
$$

After applying the changes 
$\widetilde\phi=\phi_3\circ\phi_2\circ\widehat\phi\circ\phi_1^{-1}\circ\widetilde\phi_1$ and 
$\widetilde\chi=\chi_3\circ\chi_2\circ\widehat\chi\circ\chi_1$, the ``perturbed'' symplectic structure $\widetilde\Omega$ and the momentum map $\widetilde F=(\widetilde H,\widetilde K)$ will be transformed to the desired form (\ref {eq:kal:can}).
\end{proof}

\section{Proof of Theorems \ref {thm:kudr} and \ref {thm:kudr:} on a ``hidden'' toric symmetry near a singular orbit} \label {sec:app}

\subsection {Proof of Theorem \ref {thm:kudr}}

(a) On a small neighbourhood $U(m_1)$ of the point $m_1$, we can extend the functions $f_1,\dots,f_n$ to canonical coordinates $f_i,g_i$ such that $\Omega|_{U(m_1)}=\sum\limits_{i=1}^n\ddd f_i\wedge \ddd g_i$ (Darboux coordinates). This is possible because the point $m_1$ is regular by (i). Without loss of generality, $g_i(m_1)=0$.
Consider the time-$2\pi$ map $h=\phi_{f_1}^{2\pi}:U(m_1)\to M$ of the Hamiltonian flow generated by the function $f_1$. Here $t\mapsto \phi_H^t(m)$ denotes the trajectory of the vector field $X_H$ with initial condition $\phi_H^0(m)=m$. Since $h(m_1)=m_1$, $h^*\Omega=\Omega$ and $F\circ h=F$, it follows that the map $h$ w.r.t.\ the local coordinates $f_i,g_i$ on
$$
U'(m_1):=U(m_1)\cap h^{-1}(U(m_1))
$$
has the form $h(f_1,\dots,f_n,g_1,\dots,g_n)=(f_1,\dots,f_n,g_1+h_1(f_1,\dots,f_n),\dots,g_n+h_n(f_1,\dots,f_n))$ for some functions $h_i(f_1,\dots,f_n)$ such that $h_i(0,\dots,0)=(0,\dots,0)$. We can and will assume that each common level set $U'(m_1)\cap F^{-1}(a)$ is connected (this can be achieved by choosing a smaller neighbourhood $U(m_1)$ of the point $m_1$ in $M$).

Since the map $h$ preserves the symplectic structure $\Omega|_{U'(m_1)}=\sum\limits_{i=1}^n\ddd f_i\wedge \ddd g_i$, we conclude that the functions $g_i+h_i(f_1,\dots,f_n)$ pairwise Poisson commute. Hence $\{g_i+h_i(f_1,\dots,f_n),g_j+h_j(f_1,\dots,f_n)\} = \partial h_j/\partial f_i - \partial h_i/\partial f_j$ equals $0$. It follows from the Theory of PDE's that $h_i(z_1,\dots,z_n)=\partial S(z_1,\dots,z_n)/\partial z_i$, for some function $S=S(z_1,\dots,z_n)$ on the neighbourhood $F(U'(m_1))$ of the origin, i.e.\ $h=\phi_{S\circ F}^1$ on $U'(m_1)$. Here we used the connectedness of the sets $U'(m_1)\cap F^{-1}(a)$. Due to $h_i(0,\dots,0)=(0,\dots,0)$, we have $\ddd S(0,\dots,0)=0$. We can and will assume that $S(0,\dots,0)=0$. 
The function $S$ is real-analytic, since the functions $h_i$ are. Let us define a function
$$
I(z_1,\dots,z_n):=z_1-\frac1{2\pi}S(z_1,\dots,z_n)
$$
on $F(U'(m_1))$. Let us show that it has the required properties.

From the properties of the function $S$, we have $I(z_1,\dots,z_n)=z_1+O(\sum\limits_{j=1}^n|z_j|^2)$.
From the equalities $\phi_{f_1}^{2\pi}=h=\phi_{S\circ F}^1=\phi_{\frac1{2\pi}S\circ F}^{2\pi}$, we obtain $\phi_{f_1 - \frac1{2\pi}S\circ F}^{2\pi}=\Id$ on $U'(m_1)$. That is, the function $f_1 - \frac1{2\pi}S\circ F=I\circ F$ is a $2\pi$-periodic first integral on a neighbourhood of $\gamma_1$ containing $U'(m_1)$.

The function $I\circ F$ is defined and is real-analytic on the neighbourhood $U(L_0)=F^{-1}(F(U'(m_1)))$ of the singular fiber $L_0:=F^{-1}(0,\dots,0)$. The Hamiltonian flow $\phi_{I\circ F}^{2\pi}$ of the function $I\circ F$ in time $2\pi$ is defined on some neighbourhood $U_1\subseteq U(L_0)$ of the set $C$, due to the condition (ii) and the equality $X_{I\circ F}=X_{f_1}$ on $F^{-1}(0,\dots,0)$. We can and will assume that $U_1':=U_1\cap (\phi_{I\circ F}^{2\pi})^{-1}(U_1)$ is connected, since $C$ is. Let $U$ be the union of all trajectories through points of $U_1'$ of the Hamiltonian vector field $X_{I\circ F}$, which all are $2\pi$-periodic by construction. Let $V:=F(U)$.

Thus, the map $\phi_{I\circ F}^{2\pi}$ is real-analytic, is defined on the connected open domain $U_1'$, and its restriction to the sub-domain $U'(m_1)\cap U\ne\varnothing$ is the identity. Therefore, by the uniqueness of analytic continuation, this map is the identity on the whole $U_1'$, therefore $U$ is filled by $2\pi$-periodic trajectories of $X_{I\circ F}$.

Let us prove the integral formula (\ref {eq:action:real:new}) for the action function $I$. 
Similarly to the beginning of the proof, we can extend the functions $I\circ F,f_2,\dots,f_n$ to canonical coordinates $I\circ F,f_2,\dots,f_n,u_1,\dots,u_n$ on a smaller neighbourhood of $m_1$ (we can and will denote this neighbourhood again by $U_1'$) such that $\Omega|_{U_1'}=\ddd(I\circ F)\wedge \ddd u_1+\sum\limits_{i=2}^n\ddd f_i\wedge \ddd u_i$ (Darboux coordinates).
Recall that $U$ denotes the union of all trajectories through points of $U_1'$ of the Hamiltonian vector field $X_{I\circ F}$.
Observe that the functions $u_2,\dots,u_n$ and the 1-form $\ddd u_1$ on $U_1'$ are preserved by the Hamiltonian flow generated by the action function $I\circ F$.
We extend these functions and 1-form to $U$ by making them invariant under the Hamiltonian $S^1$-action generated by $I\circ F$. Define the 1-form
\begin{equation} \label {eq:alpha}
\alpha:=(I\circ F)\ddd u_1+\sum\limits_{i=2}^n f_i \ddd u_i
\end{equation}
on $U$, then $\ddd\alpha=\Omega|_U$, thus the symplectic structure is exact on $U$.
By the Stokes formula, in order to prove (\ref {eq:action:real:new}), it suffices to prove the integral formula
\begin{equation} \label {eq:action:real}
I(z_1,\dots,z_n)
 = \frac1{2\pi}\int\limits_{C_{\gamma,\gamma_{(z_1,\dots,z_n)}}} \Omega ,
\qquad \qquad (z_1,\dots,z_n)\in V,
\end{equation}
where $C_{\gamma,\gamma_{(z_1,\dots,z_n)}}$ denotes the cylinder (with boundary $\partial C_{\gamma,\gamma_{(z_1,\dots,z_n)}}=\gamma_{(z_1,\dots,z_n)}-\gamma$) formed by the closed curves $\gamma_{(tz_1,\dots,tz_n)}$, $0\le t\le 1$.
The integral formula (\ref {eq:action:real}) follows from the following facts:
\begin{itemize}
\item the Hamiltonian flow generated by $I\circ F$ is $2\pi$-periodic on a neighbourhood $U(\gamma)\subset U$ of $\gamma=\gamma_1$, and
\item its closed orbits on each fiber $U(\gamma)\cap F^{-1}{(z_1,\dots,z_n)}$ are homological in this fiber to the closed curve $\gamma_{(z_1,\dots,z_n)}$, ${(z_1,\dots,z_n)}\in V\cap F(U(\gamma))$.
\end{itemize}
Indeed, the right-hand side of (\ref {eq:action:real}) will not change if we replace the closed curves $\gamma_{(z_1,\dots,z_n)}$ with $2\pi$-periodic trajectories $\hat\gamma_{(z_1,\dots,z_n)}\subset F^{-1}(z_1,\dots,z_n)$ of the Hamiltonian flow generated by $I\circ F$. If we compute the right-hand side of (\ref {eq:action:real}) in this way, via the Fubini formula and taking into account that $\Omega(\cdot,\frac{\ddd}{\ddd t}\hat\gamma_{(z_1,\dots,z_n)})=\ddd(I\circ F)|_{\gamma_{(z_1,\dots,z_n)}}$, then the resulting value equals $I(z_1,\dots,z_n)-I(0,\dots,0)=I(z_1,\dots,z_n)$, which is the left-hand side of (\ref {eq:action:real}), as required.

(b) Let us fix a smaller neighbourhood $U'$ of the set $C$ and a smaller neighbourhood $V'$ of the origin in $\R^n$ having compact closures $\overline{U'}\subset U$ and $\overline{V'}\subset V$.

We can extend the functions $\tilde f_1,\dots,\tilde f_n$ on $U(m_1)$ to canonical coordinates $\tilde f_i,\tilde g_i$ such that $\tilde\Omega|_{U(m_1)}=\sum\limits_{i=1}^n\ddd\tilde f_i\wedge \ddd\tilde g_i$ (``perturbed'' Darboux coordinates). 
If one follows an explicit construction of such coordinates, one can manage to have
$\|\tilde g_i-g_i\|_{C^{k-2}}=O(\varepsilon)$.

Let us consider the ``perturbed'' map $\tilde h=\phi_{\tilde f_1}^{2\pi}:U(m_1)\to M$, where $t\mapsto \phi_{\tilde H}^t(m)$ is the trajectory of the Hamiltonian vector field $X_{\tilde H}$ with the Hamilton function $\tilde H$, the symplectic structure $\tilde\Omega$ and initial condition $\phi_{\tilde H}^0(m)=m$.
Since $\tilde h^*\tilde\Omega=\tilde\Omega$ and $\tilde F\circ \tilde h=\tilde F$, it follows that the map $\tilde h$ w.r.t.\ the local coordinates $\tilde f_i,\tilde g_i$ on
$$
\tilde U(m_1):=U'(m_1)\cap \tilde h^{-1}(U'(m_1))
$$
has the form $\tilde h(\tilde f_1,\dots,\tilde f_n,\tilde g_1,\dots,\tilde g_n)=(\tilde f_1,\dots,\tilde f_n,\tilde g_1+\tilde h_1(\tilde f_1,\dots,\tilde f_n),\dots,\tilde g_n+\tilde h_n(\tilde f_1,\dots,\tilde f_n))$, for some functions $\tilde h_i(z_1,\dots,z_n)$ on $\tilde V:=\tilde F(\tilde U(m_1))$ $O(\varepsilon)-$close to $h_i(z_1,\dots,z_n)$ in $C^{k-2}-$norm.

Since the map $\tilde h$ preserves the symplectic structure $\tilde\Omega|_{\tilde U(m_1)}$, we have similarly to (a) that $\tilde h_i(z_1,\dots,z_n)=\partial \tilde S(z_1,\dots,z_n)/\partial z_i$ for some real-analytic function $\tilde S=\tilde S(z_1,\dots,z_n)$ $O(\varepsilon)-$close to the function $S$ in $C^{k-1}-$norm on some (a bit perturbed) neighbourhood $\tilde V_1$ of the origin. We have $\tilde V_1\supset V'$ if the perturbation is small enough.
Let us define on $\tilde V_1$ the real-analytic function
$$
\tilde I(z_1,\dots,z_n):=z_1-\frac1{2\pi}\tilde S(z_1,\dots,z_n).
$$

Observe that the ``perturbed'' function $\tilde I\circ \tilde F$ is $O(\varepsilon)-$close in $C^{k-1}-$norm to the ``unperturbed'' function $I\circ F$, moreover the ``unperturbed'' map $\phi_{I\circ F}^{2\pi}$ is defined and is real-analytic on the connected domain $U\supset U'$ (cf.\ (a)).
This implies that the ``perturbed'' map $\phi_{\tilde I\circ \tilde F}^{2\pi}$ is defined and is real-analytic on (a bit smaller) connected domain $\tilde U\supset U'$. But the ``perturbed'' map $\phi_{\tilde I\circ \tilde F}^{2\pi}$ is the identity on the sub-domain $\tilde U(m_1)\cap \tilde U\ne\varnothing$ (by construction of the function $\tilde I$). Therefore, by the uniqueness of analytic continuation, it equals the identity on the whole of $\tilde U$. Thus, $\tilde I\circ \tilde F$ is a $2\pi$-periodic first integral of the perturbed system on $\tilde U$.

The integral formula (\ref {eq:action:real:new:}) for the ``perturbed'' action function $\tilde I$ follows by the same arguments as for the ``unperturbed'' case. 
The estimate $\|\tilde\alpha-\alpha\|_{C^{k-3}}=O(\varepsilon)$ follows from construction.
The estimate $\|\tilde I-I\|_{C^k}=O(\varepsilon)$ follows from (\ref {eq:action:real}) and its analogue for the ``perturbed'' system.
\qed

\subsection {Proof of Theorem \ref {thm:kudr:}}

Similarly to the proof of Theorem \ref {thm:kudr} (a) (respectively, (b)), we can construct a complex-valued $2\pi$-periodic first integral $I\circ F^{\mathbb C}$ (respectively, $\tilde I\circ\tilde F^{\mathbb C}$).
Such a construction can be performed on a small open complexification $U$ (respectively, $\tilde U\supset U'$) of a small neighbourhood of the set $C$, since, by assumption, $C$ is connected and invariant under the Hamiltonian flow generated by the function $\lambda f_1^{\mathbb C}$. As a result, we obtain a holomorphic function $I=I(z_1,\dots,z_n)$ on a neighbourhood $V$ (respectively, a holomorphic function $\tilde I=\tilde I(z_1,\dots,z_n)$ on a smaller neighbourhood $\tilde V\supset V'$) of the origin in ${\mathbb C}^n$ such that
\begin{itemize}
\item $I(z_1,\dots,z_n)=\lambda z_1+O(\sum\limits_{j=1}^n|z_j|^2)$ (respectively, $\tilde I$ is $O(\varepsilon)-$close to $I$ and $\tilde I(0,\dots,0)=0$),
\item the set $U$ (respectively, $\tilde U$) is invariant under the Hamiltonian flow generated by $I\circ F^{\mathbb C}$ (respectively, by $\tilde I\circ\tilde F^{\mathbb C}$), and this flow is $2\pi$-periodic on $U$ (respectively, on $\tilde U$),
\item the integral formula (\ref {eq:action:real:C:new}) for $I(z_1,\dots,z_n)$ (respectively, (\ref {eq:action:real:C:new:}) for $\tilde I(z_1,\dots,z_n)$) holds.
\end{itemize}
Due to the Cauchy integral formula for holomorphic functions and the first property from above, on a smaller neighbourhood, the function $\tilde I$ is $O(\varepsilon)-$close to $I$ in $C^k-$norm, for any $k\in{\mathbb Z}_+$.
We want to prove that $\frac1\lambda I$ (resp., $\frac1\lambda\tilde I$)
is real-valued on $V\cap \R^n$ (resp., on $\tilde V\cap\R^n$), moreover $\lambda\in\R\cup i\R$, provided that the corresponding condition (iii) for the circle $\gamma$ of being homologically symmetric holds. Consider two cases.

{\em Case 1:} the trajectory $\gamma\ni m_1$ of the Hamiltonian flow generated by the function $\lambda f_1^{\mathbb C}$ (and, hence, by $I\circ F^{\mathbb C}$) is homological to its conjugate in the following sense.
There exists $m_1'\in U(m_1)$ such that $a_1:=F^{\mathbb C}(m_1')\in\R^n$ and the closed path $\gamma_{a_1}$ in the fiber $L_{a_1}^{\mathbb C}:= (F^{\mathbb C})^{-1}(a_1)$ is homological in the fiber $L_{a_1}^{\mathbb C}\setminus \Sing(F^{\mathbb C})$ to the closed path $\overline{\gamma_{a_1}}$ (obtained from $\gamma_{a_1}$ by $\mathbb C$-conjugation).

(a) Choose a small neighbourhood $U(m_1')$ of $m_1'$ in $U(m_1)$. Take any point $m_2$ in $U^\R(m_1'):=U(m_1')\cap (F^{\mathbb C})^{-1}(\R^n)$, thus $a_2:=F^{\mathbb C}(m_2)\in\R^n$.
Consider the closed path $\gamma_{a_2}$ in the fiber $L_{a_2}^{\mathbb C} :=
(F^{\mathbb C})^{-1}(a_2)$.
Since the momentum map $F=(f_1,\dots,f_n)$ is real-analytic and $a_2\in\R^n$, it follows that the closed path 
$\overline{\gamma_{a_2}}$ is contained in the same fiber $L_{a_2}^{\mathbb C}$, moreover the paths $\gamma_{a_2}$ and $\overline{\gamma_{a_2}}$ are homological to each other in the fiber $L_{a_2}^{\mathbb C}\setminus \Sing(F^{\mathbb C})$, due to the homological symmetry condition.
On the other hand, we have the integral formula (\ref {eq:action:real:C:new}) for $I(z_1,\dots,z_n)$, which by Stokes' formula reads:
$$
I\circ F^{\mathbb C}(m_2)-I\circ F^{\mathbb C}(m_1')
 = I(a_2) - I(a_1)
 = \frac1{2\pi}\int\limits_{C_{\gamma_{a_1},\gamma_{a_2}}} \Omega^{\mathbb C}, \qquad m_2\in U^\R(m_1'),
$$
where $C_{\gamma_{a_1},\gamma_{a_2}}$ is a cylinder in $U$ with boundary $\partial C_{\gamma_{a_1},\gamma_{a_2}}=\gamma_{a_2}-\gamma_{a_1}$. Since (the real and imaginary parts of) $\Omega^{\mathbb C}$ is closed, the integral does not depend on the choice of the cylinder in the given homotopy class relatively boundary.
We want to show that the latter integral is in fact real.
By changing coordinates under the integral, we obtain
\begin{equation} \label {eq:action}
\int\limits_{C_{\gamma_{a_1},\gamma_{a_2}}} \Omega^{\mathbb C}
 =
\int\limits_{C_{\overline{\gamma_{a_1}},\overline{\gamma_{a_2}}}} \overline{\Omega^{\mathbb C}}
 =
\overline{\int\limits_{C_{\overline{\gamma_{a_1}},\overline{\gamma_{a_2}}}} \Omega^{\mathbb C}}.
\end{equation}
Since the symplectic structure vanishes on each fiber, the resulting integral in (\ref {eq:action}) does not depend on the choice of closed paths in given homology classes, so this integral will not change if we replace $\overline{\gamma_{a_1}},\overline{\gamma_{a_2}}$ with their conjugates. Thus
$$
\int\limits_{C_{\gamma_{a_1},\gamma_{a_2}}} \Omega^{\mathbb C}
= \overline{\int\limits_{C_{\overline{\gamma_{a_1}},\overline{\gamma_{a_2}}}} \Omega^{\mathbb C}}
= \overline{\int\limits_{C_{\gamma_{a_1},\gamma_{a_2}}} \Omega^{\mathbb C}},
$$
which shows that the integral is in fact real.
We conclude that
$$
I\circ F^{\mathbb C}(m_2)-I\circ F^{\mathbb C}(m_1') \in\R, \qquad m_2\in U^\R(m_1'),
$$
thus $\im(I\circ F^{\mathbb C})$ is constant on $U^\R(m_1')$.
Since $m_1'\in U(m_1)$ and $U(m_1)$ consists of regular points of $F^{\mathbb C}$, the set $V(a_1):=F(U^\R(m_1'))$ is open in $\R^n$, so it is a neighbourhood of $a_1$ in $\R^n$.
By above, $\im I({z_1},\dots,{z_n})$ is constant on this neighbourhood. Hence $I(z_1,\dots,z_n)$ up to an additive constant is real-analytic on $V(a_1)\subset V\cap \R^n$.
This implies, by the uniqueness of analytic continuation, that $I(z_1,\dots,z_n)$ is real-analytic on the entire neighbourhood $V\cap \R^n$ of $(0,\dots,0)$, up to an additive constant. The latter additive constant is in fact real, since $I(0,\dots,0)=0$ by properties of the function $I=I(z_1,\dots,z_n)$.

It remains to show that $\lambda\in\R$. On one hand, $\lambda=\frac{\partial I}{\partial z_1}(0,\dots,0)$. On the other hand, $I(z_1,\dots,z_n)\in\R$ for any $(z_1,\dots,z_n)\in V\cap \R^n$. Hence $\lambda\in\R$.

(b) For a perturbed system, the proof follows the same arguments.
In more detail, we take any two points $\tilde m_1,\tilde m_2$ in $\tilde U^\R(m_1'):=U(m_1')\cap (\tilde F^{\mathbb C})^{-1}(\R^n)$, thus $\tilde a_j:=\tilde F^{\mathbb C}(\tilde m_j)\in\R^n$ ($j=1,2$).
Consider the closed path $\tilde\gamma_{\tilde a_j}$ in the fiber $\tilde L_{\tilde a_j}:=(\tilde F^{\mathbb C})^{-1}(\tilde a_j)$ ($j=1,2$).
Since the ``perturbed'' momentum map $\tilde F=(\tilde f_1,\dots,\tilde f_n)$ is real-analytic and $\tilde a_j\in\R^n$, it follows that the closed path $\overline{\tilde\gamma_{\tilde a_j}}$ is contained in the same fiber $\tilde L_{\tilde a_j}$, moreover $\tilde\gamma_{\tilde a_j}$ and $\overline{\tilde\gamma_{\tilde a_j}}$ are homological to each other in the fiber $\tilde L_{\tilde a_j}\setminus \Sing(\tilde F^{\mathbb C})$, due to the homological symmetry condition ($j=1,2$).
On the other hand, we have the integral formula (\ref {eq:action:real:C:new:}) for $\tilde I(z_1,\dots,z_n)$, which by Stokes' formula reads:
$$
\tilde I\circ\tilde F^{\mathbb C}(\tilde m_2)-\tilde I\circ\tilde F^{\mathbb C}(\tilde m_1)
 = \tilde I({\tilde a_2}) - \tilde I({\tilde a_1})
 = \frac1{2\pi}\int\limits_{C_{\tilde \gamma_{\tilde a_1},\tilde\gamma_{\tilde a_2}}} \tilde\Omega^{\mathbb C},\qquad \tilde m_1,\tilde m_2\in \tilde U^\R(m_1'),
$$
where $C_{\tilde\gamma_{\tilde a_1},\tilde\gamma_{\tilde a_2}}$ is a cylinder with boundary $\partial C_{\tilde\gamma_{\tilde a_1},\tilde\gamma_{\tilde a_2}}=\tilde\gamma_{\tilde a_2}-\tilde\gamma_{\tilde a_1}$.
Similarly to (\ref {eq:action}), we compute
\begin{equation} \label {eq:tilde:action}
\int\limits_{C_{\tilde\gamma_{\tilde a_1},\tilde\gamma_{\tilde a_2}}} \tilde\Omega^{\mathbb C}
 =
\int\limits_{C_{\overline{\tilde\gamma_{\tilde a_1}},\overline{\tilde\gamma_{\tilde a_2}}}} \overline{\tilde\Omega^{\mathbb C}}
 =
\overline{\int\limits_{C_{\overline{\tilde\gamma_{\tilde a_1}},\overline{\tilde\gamma_{\tilde a_2}}}} \tilde\Omega^{\mathbb C}}.
\end{equation}
Similarly to the ``unperturbed'' system case, the resulting integral in (\ref {eq:tilde:action}) will not change if we replace the closed curves $\overline{\tilde\gamma_1},\overline{\tilde\gamma_2}$ with their conjugates, due to the homological symmetry condition. Thus
$$
\int\limits_{C_{\tilde\gamma_{\tilde a_1},\tilde\gamma_{\tilde a_2}}} \tilde\Omega^{\mathbb C}
 =
\overline{\int\limits_{C_{\overline{\tilde\gamma_{\tilde a_1}},\overline{\tilde\gamma_{\tilde a_2}}}} \tilde\Omega^{\mathbb C}}
=
\overline{\int\limits_{C_{\tilde\gamma_{\tilde a_1},\tilde\gamma_{\tilde a_2}}} \tilde\Omega^{\mathbb C}},
$$
thus the integral is in fact real, thus
$$
\tilde I\circ \tilde F^{\mathbb C}(\tilde m_2)-\tilde I\circ \tilde F^{\mathbb C}(\tilde m_1) \in\R,\qquad \tilde m_1,\tilde m_2\in \tilde U^\R(m_1').
$$
Similarly to the ``unperturbed'' system case, we conclude that
$\im\tilde I({z_1},\dots,{z_n})$ is
constant on the neighbourhood $\tilde V(a_1):=\tilde F(\tilde U^\R(m_1'))$ of $a_1$ in $\R^n$.
Hence $\tilde I(z_1,\dots,z_n)$ up to an additive constant is real-analytic on $\tilde V(a_1)\subset\tilde V\cap\R^n$.
This implies that $\tilde I(z_1,\dots,z_n)$ is real-analytic on the entire neighbourhood $\tilde V\cap \R^n$ of $(0,\dots,0)$, up to an additive constant. The latter additive constant is in fact real, since $\tilde I(0,\dots,0)=0$ by construction (see properties of the function $\tilde I$).

{\em Case 2:} the trajectory $\gamma\ni m_1$ of the Hamiltonian flow generated by the function $\lambda f_1^{\mathbb C}$ (and, hence, by $I\circ F^{\mathbb C}$) is reverse-homological to its conjugate in the following sense.
There exists $m_1'\in U(m_1)$ such that $a_1:=F^{\mathbb C}(m_1')\in\R^n$ and the closed path $\gamma_{a_1}$ in the fiber $L_{a_1}^{\mathbb C}:= (F^{\mathbb C})^{-1}(a_1)$ is homological in the fiber $L_{a_1}^{\mathbb C}\setminus \Sing(F^{\mathbb C})$ to the closed path obtained from the ($\mathbb C$-conjugated to $\gamma_{a_1}$) path $\overline{\gamma_{a_1}}$ by reversing orientation.

Similarly to the proof in the case $\lambda\in\R$, we conclude that each resulting integral in (\ref {eq:action}) and (\ref {eq:tilde:action}) will change to the opposite if we replace the closed curves $\overline{\gamma_{a_1}},\overline{\gamma_{a_2}}$ and $\overline{\tilde\gamma_{\tilde a_1}},\overline{\tilde\gamma_{\tilde a_2}}$ with their conjugates.
This immediately shows that each of these integrals is in fact purely imaginary.
We conclude that
$$
I\circ F^{\mathbb C}(m_2)-I\circ F^{\mathbb C}(m_1') \in i\R, \qquad
\tilde I\circ \tilde F^{\mathbb C}(\tilde m_2)-\tilde I\circ \tilde F^{\mathbb C}(\tilde m_1) \in i\R.
$$
Thus $\re I({z_1},\dots,{z_n})$ is constant on $V(a_1)$, and
$\re\tilde I({z_1},\dots,{z_n})$ is constant on $\tilde V(a_1)$.
Hence $I(z_1,\dots,z_n)$ and $\tilde I(z_1,\dots,z_n)$ up to additive constants are imaginary-valued on open subsets $V(a_1)\subset V\cap\R^n$ and $\tilde V(a_2)\subset \tilde V\cap \R^n$, respectively.
Therefore $iI(z_1,\dots,z_n)$ and $i\tilde I(z_1,\dots,z_n)$ are real-analytic on the entire neighbourhoods $V\cap\R^n$ and $\tilde V\cap \R^n$ of the origin in $\R^n$, up to additive constants. The latter additive constants are in fact real, since $I(0,\dots,0)=\tilde I(0,\dots,0)=0$ by construction.

It remains to show that $\lambda\in i\R$. On one hand, $\lambda=\frac{\partial I}{\partial z_1}(0,\dots,0)$. On the other hand, $I(z_1,\dots,z_n)\in i\R$ for any $(z_1,\dots,z_n)\in V\cap \R^n$. Hence $\lambda\in i\R$.
\qed

\section{Equivariant symplectic normalization of a torus action near a fixed point} \label {sec:app:fixed}

In this section, we study equivariant version of the theorems \ref {thm:period} and \ref {thm:period:} in the case when the orbit ${\cal O}_{m_0}$ is a point, and we prove that the torus action can be normalized symplectically in an equivariant way.

For proving Theorem \ref {thm:period}, we will need the following lemma about an equivariant normal form of a Hamiltonian torus action near a fixed point.

\begin{Lemma} \label {lem:period}
Suppose we are given a Hamiltonian $(S^1)^{\varkappa}$-action generated by $C^\infty$-smooth functions $J_1,\dots,J_\varkappa$ 
on a $C^\infty$-smooth symplectic manifold $(M,\Omega)$, where the momentum map $F=(J_1,\dots,J_\varkappa)$ is not necessarily proper.
Suppose a point $m_0\in M$ is fixed under this $(S^1)^{\varkappa}$-action.
Suppose a finite Abelian group $\Gamma$ acts on a neighbourhood of $m_0$ by $F$-preserving symplectomorphisms $\rho(\psi)$, $\psi\in\Gamma$, fixing the point $m_0$. Here $\rho:\Gamma\to\Symp(M,\Omega)$ is a homomorphism from $\Gamma$ to the group of (local) symplectomorphisms.
Then:

{\rm (a)}
There exist a $\Gamma\times(S^1)^{\varkappa}$-invariant neighbourhood $U$ of $m_0$ and smooth local coordinates $x_1,y_1,\dots,x_{n},y_{n}$ on $U$ such that $x_j(m_0)=y_j(m_0)=0$, $1\le j\le n$, and
\begin{equation} \label {eq:linear:coord} 
\Omega|_U =
\sum\limits_{j=1}^{n} \ddd x_{j}\wedge \ddd y_{j} , \quad
J_\ell|_U = c_\ell+\sum\limits_{j=1}^{k_e}p_{j\ell}\frac{x_{j}^2+y_{j}^2}2, \quad
1\le \ell\le \varkappa,
\end{equation}
for some constants $c_\ell\in\R$ and integers $k_e\in{\mathbb Z}_+$, $p_{j\ell}\in\mathbb Z$. Moreover, the $\rho(\Gamma)$-action on $U$ has the form
\begin{equation} \label {eq:twist}
\rho(\psi) : (z_1,\dots,z_{n}) \mapsto (e^{2\pi iq_{\psi,1}}z_1,\dots,e^{2\pi iq_{\psi,n}}z_{n}), \qquad \psi \in \Gamma,
\end{equation}
for some $q_{\psi,1},\dots,q_{\psi,n}\in{\mathbb Q}$ depending on $\psi$ such that $q_{\psi,j}\in\frac12\mathbb Z$ for $k_e+1\le j\le n$. Here we used the notation $z_j:=x_j+iy_j$, $1\le j\le n$. 
If the system is analytic, the coordinates $x_j,y_j$ are analytic too.

{\rm (b)} The equivariant symplectic normalization (\ref {eq:linear:coord})--(\ref {eq:twist}) for a torus action is persistent under $C^\infty$-smooth perturbations in the following sense.
For any integer $k\ge4$ and any neighbourhood $U'$ of $m_0$ having a compact closure $\overline{U'}\subset U$, there exists $\varepsilon>0$ satisfying the following.
Suppose $\tilde F=(\tilde J_1,\dots,\tilde J_\varkappa)$ is a (``perturbed'') momentum map $\varepsilon$-close to $F$ in $C^k-$norm, and $\tilde\Omega$ is a (``perturbed'') symplectic structure $\varepsilon$-close to $\Omega$ in $C^{k+n-2}-$norm.
Suppose the functions $\tilde J_1,\dots,\tilde J_\varkappa$ generate a Hamiltonian $(S^1)^{\varkappa}$-action on $(\tilde M,\tilde\Omega)$, where $U\subset\tilde M\subset M$.
Suppose we are given a ``perturbed'' action of the group $\Gamma$ on $\tilde M$ by $\tilde F$-preserving symplectomorphisms $\tilde\rho(\psi)$, $\psi\in\Gamma$, where $\tilde\rho:\Gamma\to\Symp(\tilde M,\tilde\Omega)$ is a homomorphism $\varepsilon-$close to $\rho$ in $C^{k-1}-$norm.
Then there exist a point $\tilde m_0\in U'$ fixed under the ``perturbed'' $\Gamma\times(S^1)^\varkappa$-action, an invariant (under the ``perturbed'' $\Gamma\times(S^1)^\varkappa$-action) neighbourhood $\tilde U\supset U'$, and smooth local coordinates
$\tilde x_1,\tilde y_1,\dots,\tilde x_{n},\tilde y_{n}$ on $\tilde U$ $O(\varepsilon)-$close to
$x_1,y_1,\dots,x_{n},y_{n}$ in $C^{k-4}-$norm such that $\tilde x_j(\tilde m_0)=\tilde y_j(\tilde m_0)=0$, $1\le j\le n$, and an analogue of (\ref {eq:linear:coord}) holds for $\tilde\Omega|_{\tilde U}$ and $\tilde J_\ell|_{\tilde U}$ w.r.t.\ the ``perturbed'' coordinates, with the same integers $p_{j\ell}$ as in {\rm(a)} and with some constants $\tilde c_\ell\in\R$ close to $c_\ell$.
Moreover, an analogue of {\rm (\ref  {eq:twist})} holds for the ``perturbed'' action $\tilde\rho(\psi)$ of each element $\psi \in \Gamma$ on $\tilde U$, with the same rational numbers $q_{\psi,j}$ as in {\rm(a)}.
If the systems are analytic, the coordinates $\tilde x_j,\tilde y_j$ are analytic too. 
If the 
system $(M,\Omega,F)$ and the action $\rho(\Gamma)$ depend smoothly (resp., analytically) on a local parameter, i.e.\ we have a local family of systems with actions, then the local coordinates can also be chosen to depend smoothly (resp., analytically) on that parameter.
\end{Lemma}

The following lemma generalizes the real-analytic part of the previous lemma to the case when the torus acts holomorphically on a small open complexification of the manifold, and the generating functions of this action are real-analytic functions some of which are multiplied by $\sqrt{-1}$. This lemma will be used for proving Theorem \ref {thm:period:}.

\begin{Lemma} \label {lem:period:}
Suppose we are given real-analytic functions $J_1,\dots,J_{\varkappa_e+\varkappa_h}$ 
on a real-analytic symplectic manifold $(M,\Omega)$, where the map $F=(J_1,\dots,J_{\varkappa_e+\varkappa_h})$ is not necessarily proper.
Suppose the functions $J_1^{\mathbb C},\dots,J_{\varkappa_e}^{\mathbb C}$, $iJ_{\varkappa_e+1}^{\mathbb C},\dots,iJ_{\varkappa_e+\varkappa_h}^{\mathbb C}$ generate a Hamiltonian $(S^1)^{\varkappa_e+\varkappa_h}$-action on $(M^{\mathbb C},\Omega^{\mathbb C})$, where $\varkappa_e,\varkappa_h\in{\mathbb Z}_+$.
Suppose a point $m_0\in M$ is fixed under this $(S^1)^{\varkappa_e+\varkappa_h}$-action.
Suppose a finite Abelian group $\Gamma$ acts on a neighbourhood of $m_0$ by $F$-preserving real-analytic symplectomorphisms $\rho(\psi)$, $\psi\in\Gamma$, fixing the point $m_0$. Here $\rho:\Gamma\to\Symp(M,\Omega)$ is a homomorphism from $\Gamma$ to the group of (local) symplectomorphisms.
Then:

{\rm (a)}
There exist a $\Gamma\times(S^1)^{\varkappa_e}$-invariant neighbourhood $U$ of $m_0$ in $M$ and real-analytic local coordinates $x_1,y_1,\dots,x_{n},y_{n}$ on $U$ such that $x_j(m_0)=y_j(m_0)=0$, $1\le j\le n$, and
\begin{equation} \label {eq:linear:coord:omega}
\Omega|_{U} = \sum\limits_{j=1}^{n} \ddd x_{j}\wedge \ddd y_{j} ,
\end{equation}
\begin{equation} \label {eq:linear:coord:ell}
J_\ell|_{U}=c_\ell
+\sum\limits_{j=1}^{k_f}             p_{j\ell} h_{2j-1}
+\sum\limits_{j=1}^{k_e} p_{k_f+j,\ell} h_{2k_f+j}
, \quad 1\le \ell\le \varkappa_e,
\end{equation}
\begin{equation} \label {eq:linear:coord:hyp}
J_\ell|_{U}=c_\ell
+\sum\limits_{j=1}^{k_f}
p_{j,\ell} h_{2j}
+\sum\limits_{j=1}^{k_h} p_{k_f+j,\ell} h_{2k_f+k_e+j}, \quad \varkappa_e+1\le \ell\le \varkappa_e+\varkappa_h
\end{equation}
(cf.\ (\ref {eq:mz:hj:})), for some real constants $c_\ell\in\R$ and integers $p_{j\ell}\in\mathbb Z$ and
$k_f,k_e,k_h\in{\mathbb Z}_+$ such that $2k_f+k_e+k_h\le n$. Moreover, the $\rho(\Gamma)$-action on $U$ has the form
\begin{equation} \label {eq:twist:}
\rho(\psi) : (z_1,\dots,z_{n}) \mapsto (e^{2\pi iq_{\psi,1}}z_1,\dots,e^{2\pi iq_{\psi,n}}z_{n}), \qquad \psi \in \Gamma,
\end{equation}
for some $q_{\psi,1},\dots,q_{\psi,n}\in{\mathbb Q}$ depending on $\psi$ such that $q_{\psi,2j-1}+q_{\psi,2j}\in{\mathbb Z}$ for $1\le j\le k_f$, and $q_{\psi,j}\in\frac12{\mathbb Z}$ for $2k_f+k_e+1\le j\le n$.

{\rm (b)} The equivariant symplectic normalization (\ref {eq:linear:coord:omega})--(\ref {eq:twist:}) for a torus action is persistent under analytic perturbations in the following sense.
For any $k\in{\mathbb Z}_+$ and any neighbourhood $U'$ of the point $m_0$ in $M$ having a compact closure $\overline{U'}\subset U$, there exists $\varepsilon>0$ satisfying the following.
Suppose $\tilde F=(\tilde J_1,\dots,\tilde J_{\varkappa_e+\varkappa_h})$ and $\tilde\Omega$ are analytic (``perturbed'') momentum map and symplectic structure 
whose holomorphic extensions to $M^{\mathbb C}$ are $\varepsilon$-close to $F^{\mathbb C}$ and $\Omega^{\mathbb C}$, respectively, in $C^0-$norm.
Suppose the functions $\tilde J_1^{\mathbb C},\dots,\tilde J_{\varkappa_e}^{\mathbb C}$, $i\tilde J_{\varkappa_e+1}^{\mathbb C},\dots,i\tilde J_{\varkappa_e+\varkappa_h}^{\mathbb C}$ generate a Hamiltonian $(S^1)^{\varkappa_e+\varkappa_h}$-action on $(\tilde M^{\mathbb C},\tilde\Omega^{\mathbb C})$, where $U\subset\tilde M^{\mathbb C}\subset M^{\mathbb C}$.
Suppose we are given a ``perturbed'' action of the group $\Gamma$ on $\tilde M=\tilde M^{\mathbb C}\cap M$ by $\tilde F$-preserving analytic symplectomorphisms $\tilde\rho(\psi)$, $\psi\in\Gamma$, where $\tilde\rho:\Gamma\to\Symp(\tilde M,\tilde\Omega)$ is a homomorphism whose holomorphic extension to $\tilde M^{\mathbb C}$ is $\varepsilon-$close to $\rho^{\mathbb C}$ in $C^0-$norm.
Then there exist a point $\tilde m_0\in U'$ fixed under the ``perturbed'' $\Gamma\times(S^1)^{\varkappa_e}$-action, an invariant (under the ``perturbed'' $\Gamma\times(S^1)^{\varkappa_e}$-action) neighbourhood $\tilde U\supset U'$, and analytic  (``perturbed'') local coordinates
$\tilde x_1,\tilde y_1,\dots,\tilde x_{n},\tilde y_{n}$ on $\tilde U$ $O(\varepsilon)-$close to $x_1,y_1,\dots,x_{n},y_{n}$ in $C^{k}-$norm, 
such that $\tilde x_j(\tilde m_0)=\tilde y_j(\tilde m_0)=0$, $1\le j\le n$, and analogues of (\ref {eq:linear:coord:omega})--(\ref {eq:linear:coord:hyp}) hold for $\tilde\Omega|_{\tilde U}$ and $\tilde J_\ell|_{\tilde U}$ w.r.t.\ the ``perturbed'' coordinates,
with the same integers $p_{j\ell}$ and $k_f,k_e,k_h$ as in {\rm(a)} and with some constants $\tilde c_\ell\in\R$.
Moreover, an analogue of {\rm (\ref  {eq:twist:})} holds for the ``perturbed'' action $\tilde\rho(\psi)$ of each element $\psi \in \Gamma$ on $\tilde U$, with the same rational numbers $q_{\psi,j}$ as in {\rm(a)}.
If the system $(M,\Omega,F)$ and the action $\rho(\Gamma)$ depend on a local parameter, i.e.\ we have a local family of systems with actions, 
moreover the holomorphic extensions of the system and the action to $M^{\mathbb C}$ depend smoothly (resp., analytically) on that parameter,
then the local coordinates can also be chosen to depend smoothly (resp., analytically) on that parameter.
\end{Lemma}

\subsection{Proof of part (a) of Lemmata \ref {lem:period} and \ref {lem:period:}} \label {subsec:proof:lem:period}

We will give a proof of Lemma \ref {lem:period:} (a) in the real-analytic case. If $\varkappa_h=0$ (i.e., all functions $J_\ell$ generate $2\pi$-periodic flows), all our arguments and constructions literally work both in the real-analytic and $C^\infty$ cases. This gives a proof of Lemma \ref {lem:period} (a) too.

{\em Step 1.}
On the vector space $V:=T_{m_0}M$, consider the linear operators
\begin{itemize}
\item $A_\ell$, $1\le\ell\le \varkappa_e+\varkappa_h$, the linearizations of the vector fields $X_{J_\ell}$ at their common equilibrium point $m_0$,
\item $M_a$, $1\le a\le r$, the linearizations of the symplectomorphisms $\rho(\psi^a)$ at their common fixed point $m_0$, where $\psi^a\in\Gamma$ are generators of the group $\Gamma$.
\end{itemize}
Observe that the operators $A_\ell$ ($1\le\ell\le \varkappa_e+\varkappa_h$) and $M_a$ ($1\le a\le r$) are
\begin{itemize}
\item Hamiltonian and symplectic (respectively) w.r.t.\ the symplectic form $\Omega|_{m_0}$,
\item semisimple (i.e.\ diagonalizable over $\mathbb C$) and pairwise commute,
\item the operators $A_1,\dots,A_{\varkappa_e}$ and $M_a$ are
elliptic, while the operators $A_{\varkappa_e+1},\dots,A_{\varkappa_e+\varkappa_h}$ are hyperbolic. In other words, all eigenvalues of $A_1,\dots,A_{\varkappa_e}$ belong to the unit circle in $\mathbb C$, all eigenvalues of $M_a$ are purely imaginary (i.e.\ belong to $i\R$), all eigenvalues of $A_{\varkappa_e+1},\dots,A_{\varkappa_e+\varkappa_h}$ are real.
\end{itemize}
It follows from a standard result of Linear Algebra that the linear hull of the operators $A_\ell$ is contained in a maximal Abelian subalgebra consisting of semisimple elements (called a {\em Cartan subalgebra}) of the Lie algebra $sp(2n,\R)$ of the Lie group $Sp(2n,\R)$ of real linear symplectomorphisms. Moreover, each symplectic operator $M_a$ is the exponent of some element of this subalgebra.

Furthermore, it is known that each Cartan subalgebra of $sp(2n,\R)$ is uniquely determined (up to conjugation) by a triple $(\hat k_e,\hat k_h,\hat k_f)$ of non-negative integers such that $\hat k_e+\hat k_h+2\hat k_f=n$.
If we decompose the symplectic vector space $(\R^{2n}_{(x_1,\dots,x_n,y_1,\dots,y_n)},\sum\limits_{j=1}^n\ddd x_j\wedge \ddd y_j)$ into the direct product of $\hat k_e+\hat k_h$ symplectic subspaces $\R^2_{(x_j,y_j)}$, $1\le j\le \hat k_e+\hat k_h$, and $\hat k_f$ symplectic subspaces $\R^4_{(x_{2j-1},x_{2j},y_{2j-1},y_{2j})}$, $1\le j\le \hat k_f$, then the corresponding Cartan subalgebra (as a vector space) is the direct product of $\hat k_e$ copies of the ``elliptic'' Cartan subalgebra of $sp(2,\R)$, $\hat k_h$ copies of the ``hyperbolic'' Cartan subalgebra of $sp(2,\R)$, and $\hat k_f$ copies of the ``focus-focus'' Cartan subalgebra of $sp(4,\R)$.
Here the ``elliptic'' (respectively, ``hyperbolic'') Cartan subalgebra of $sp(\R^2_{(x_j,y_j)})$ is one-dimensional and is spanned by the Hamiltonian operator with the quadratic Hamilton function $(x_j^2+y_j^2)/2$ (respectively, $x_jy_j$), while the ``focus-focus'' Cartan subalgebra of $sp(\R^4_{(x_{2j-1},x_{2j},y_{2j-1},y_{2j})})$ is two-dimensional and is spanned by two Hamiltonian operators with the quadratic Hamilton functions
$x_{2j-1}y_{2j}-y_{2j-1}x_{2j}$ and $x_{2j-1}y_{2j-1}+x_{2j}y_{2j}$ (cf.\ (\ref {eq:linear:coord:ell}), (\ref {eq:linear:coord:hyp})).

Using the above facts and taking into account that the flows of $X_{J_\ell}$, $1\le\ell\le \varkappa_e$, and $iX_{J_\ell}$, $\varkappa_e+1\le\ell\le \varkappa_e+\varkappa_h$, are $2\pi$-periodic, one can show the existence of non-negative integers $k_e\le\hat k_e$, $k_h\le\hat k_h$, $k_f\le\hat k_f$ and coordinates $\hat x_j,\hat y_j$ on $V$, in which the following formulae (\ref {eq:hat:omega:})--(\ref {eq:twist:::}) corresponding to (\ref {eq:linear:coord:omega})--(\ref {eq:twist:}) hold:
\begin{equation} \label {eq:hat:omega:}
\Omega|_{m_0}=\sum\limits_{j=1}^n \ddd\hat x_j\wedge \ddd\hat y_j,
\end{equation}
\begin{equation} \label {eq:d2J:ell}
\ddd^2J_\ell|_{m_0}
=\sum\limits_{j=1}^{k_f}     p_{j\ell} (\ddd\hat x_{2j-1}^2+\ddd\hat y_{2j-1}^2-\ddd\hat x_{2j}^2-\ddd\hat y_{2j}^2)
+\sum\limits_{j=2k_f+1}^{2k_f+k_e} p_{j-k_f,\ell} (\ddd\hat x_{j}^2+\ddd\hat y_{j}^2), \ \ell\le \varkappa_e,
\end{equation}
\begin{equation} \label {eq:d2J:hyp}
\ddd^2J_\ell|_{m_0}
=\sum\limits_{j=1}^{k_f}
p_{j\ell} (\ddd\hat x_{2j-1}\ddd\hat y_{2j}+\ddd\hat x_{2j}\ddd\hat y_{2j-1})
+\sum\limits_{j=2k_f+k_e+1}^{2k_f+k_e+k_h} p_{j-k_f-k_e,\ell} \ddd\hat x_{j} \ddd\hat y_{j}
, \ \ell>\varkappa_e,
\end{equation}
$1\le\ell\le \varkappa_e+\varkappa_h$, for some integers $p_{j\ell}$,
\begin{equation} \label {eq:twist::}
M_a : (\hat z_1,\dots,\hat z_{n}) \mapsto (e^{2\pi iq_{a,1}}\hat z_1,\dots,e^{2\pi iq_{a,n}}\hat z_{n}), \qquad 1\le a\le r,
\end{equation}
for some $q_{a,1},\dots,q_{a,n}\in{\mathbb Q}$ such that
\begin{equation} \label {eq:twist:::}
q_{a,2j-1}+q_{a,2j}\in{\mathbb Z} \quad \mbox{for} \quad 1\le j\le k_f, \qquad
q_{a,2k_f+k_e+j}\in\frac12{\mathbb Z} \quad \mbox{for} \quad 1\le j\le k_h.
\end{equation}
Here we denoted $\hat z_j:=\hat x_j+\hat y_j$.

Let us explicitely describe (in Substeps 1a--1d below) a construction of such coordinates $\hat x_j,\hat y_j$ on $V$ satisfying (\ref {eq:hat:omega:})--(\ref {eq:twist:::}).

{\em Substep 1a.} Since the vector fields $X_{J_\ell}$ pairwise commute and are $\rho(\Gamma)_*$-invariant, moreover $\Gamma$ is commutative, we conclude that the operators $A_\ell$ and $M_a$ pairwise commute too. From a standard assertion of Linear Algebra, there exists a unique (up to a permutation of terms) decomposition
$$
V=\bigoplus\limits_{s}V_s
$$
such that each $V_s$ is invariant under each $A_\ell$ and $M_a$, moreover $\spec(A_\ell|_{V_s})=\{\pm\lambda_{\ell s},\pm\overline{\lambda_{\ell s}}\}$, $\spec(M_a|_{V_s})=\{\mu_{as}^{\pm1},\overline{\mu_{as}^{\pm1}}\}$, and for any $s\ne s'$ there exists either $\ell\in\{1,\dots,\varkappa_e+\varkappa_h\}$ such that $\lambda_{\ell s'}\not\in\{\pm\lambda_{\ell s},\pm\overline{\lambda_{\ell s}}\}$ or $a\in\{1,\dots,r\}$ such that $\mu_{as'}\not\in\{\mu_{as}^{\pm1},\overline{\mu_{as}^{\pm1}}\}$.

Consider the symplectic form $\Omega|_{m_0}$ on $V$, and denote it by $\Omega$. Since the operators $A_\ell$ and $M_a$ are Hamiltonian and symplectic respectively, moreover they pairwise commute, it follows that the subspaces $V_s$ are symplectic and pairwise skew-orthogonal \cite[Proposition 3.1]{gol:ste}.
On the other hand, these subspaces are invariant under each operator $A_\ell$. Therefore $V_s$ are pairwise ``orthogonal'' w.r.t.\ the second differential of each function $J_\ell$, that is given by the symmetric bilinear form $\ddd^2J_\ell|_{m_0}(\xi_1,\xi_2)=\Omega(\xi_1,A_\ell\xi_2)$ on $V$. So, it remains to compute the restrictions of all $\ddd^2J_\ell|_{m_0}$ and $M_a$ to each subspace $V_s$.

Let us fix $\ell\in\{1,\dots,\varkappa_e+\varkappa_h\}$. By assumption, either the flow of $X_{J_\ell}$ is $2\pi$-periodic (for $1\le\ell\le\varkappa_e$), or the flow of $X_{iJ_\ell}^{\mathbb C}$ is $2\pi$-periodic (for $\varkappa_e<\ell\le \varkappa_e+\varkappa_h$). 
As we observed at the beginning of Step 1, in the former case, $A_\ell$ is {\em elliptic}; all its eigenvalues are pure imaginary and belong to $i\mathbb Z$.
In the latter case, $A_\ell$ is {\em hyperbolic}; all its eigenvalues are integers.
Without loss of generality, we can assume that
$\lambda_{\ell s}\in i{\mathbb Z}_+$ if $A_\ell$ is elliptic (i.e., $1\le\ell\le \varkappa_e$), and
$\lambda_{\ell s}\in {\mathbb Z}_+$ if $A_\ell$ is hyperbolic (i.e., $\varkappa_e<\ell\le \varkappa_e+\varkappa_h$).

Let us fix $a\in\{1,\dots,r\}$. By assumption, the operator $M_a$ is of finite order, hence all its eigenvalues $\mu_{as}^{\pm1},\overline{\mu_{as}^{\pm1}}$ belong to the unit circle in $\mathbb C$ and belong to $\{e^{2\pi i q}\mid q\in{\mathbb Q}\}$.
Thus, without loss of generality, we can assume that
\begin{equation} \label {eq:q:mu}
\mu_{as}=e^{2\pi iq_{as}} \quad \mbox{for some} \quad q_{as}\in{\mathbb Q}\cap[0,\frac12], \qquad
1\le a\le r.
\end{equation}

{\em Substep 1b.} Fix a subspace $V_s\subseteq V$ from the above decomposition.

Consider the set
$$
\M_s:=\{a\in\{1,\dots,r\}\mid\mu_{as}\not\in \{1,-1\} \}.
$$
Thus $q_{as}\in{\mathbb Q}\cap(0,1/2)$ for each $a\in\M_s$ (see (\ref {eq:q:mu})).
Define the linear operators
$$
L_{as}:= \frac{1}{\sin(2\pi q_{a,s})} M_a|_{V_s} - (\cot(2\pi q_{a,s})) \Id_{V_s}, \qquad a \in \M_s,
$$
thus we have
\begin{equation}\label {eq:L:M}
M_a|_{V_s} = (\cos(2\pi q_{a,s})) \Id_{V_s} + (\sin(2\pi q_{a,s})) L_{as} , \qquad 1\le a \le r.
\end{equation}
One can easily show that $L_{as}$ are Hamiltonian (and symplectic), pairwise commute and satisfy the equalities 
\begin{equation} \label {eq:involution:L}
L_{as}^2=-\Id_{V_s}, \qquad a \in \M_s.
\end{equation}
If $\M_s\ne\varnothing$, let us choose $a_s\in\M_s$ and consider the (elliptic Hamiltonian) operator on $V_s$:
\begin{equation} \label {eq:op:ell}
E_s:=L_{a_s,s}.
\end{equation}

Consider the sets
$$
\I_s:=\{\ell\in\{1,\dots,\varkappa_e\}\mid\lambda_{\ell s}\in i{\mathbb Z}\setminus\{0\}\}, \quad
\H_s:=\{\ell\in\{\varkappa_e+1,\dots,\varkappa_e+\varkappa_h\}\mid\lambda_{\ell s}\in {\mathbb Z}\setminus\{0\}\}.
$$
So, all $A_\ell|_{V_s}$ with $\ell\in\I_s$ are elliptic, while all $A_\ell|_{V_s}$ with $\ell\in\H_s$ are hyperbolic.
Thus $\lambda_{\ell s}=p_{s\ell} i$ for some integer $p_{s\ell}>0$ if $\ell\in\I_s$,
$\lambda_{\ell' s}=:p_{s\ell'}>0$ is an integer if $\ell'\in\H_s$.
Define the linear operators
\begin{equation} \label {eq:B:A}
B_{\ell s}:=\frac{1}{p_{s\ell}} A_\ell|_{V_s}, \qquad \ell \in \I_s \cup \H_s.
\end{equation}
Clearly they pairwise commute, are diagonalizable over $\mathbb C$, and hence satisfy the equalities 
\begin{equation} \label {eq:involution:B}
B_{\ell s}^2=-\Id_{V_s} \quad \mbox{for }\ell \in \I_s, \qquad 
B_{\ell' s}^2=\Id_{V_s} \quad \mbox{for }\ell' \in \H_s.
\end{equation}

If $\M_s=\varnothing$ and $\I\ne\varnothing$, let us choose $\ell_s\in\I_s$ and consider the (elliptic Hamiltonian) operator 
\begin{equation} \label {eq:op:ell:}
E_s:=B_{\ell_s,s}
\end{equation}
on $V_s$. If $\H\ne\varnothing$, let us choose $\ell_s'\in\H_s$ and consider the (hyperbolic Hamiltonian) operator
\begin{equation} \label {eq:op:hyp}
H_s:=B_{\ell_s',s}
\end{equation}
on $V_s$. Consider the set of pairwise commuting symplectic involutions (see below)
\begin{equation} \label {eq:involutions}
\{ E_{s} L_{as} \mid a \in \M_s \} \ \cup \
\{ E_{s} B_{\ell s} \mid \ell\in \I_s \} \ \cup \
\{ H_{s} B_{\ell' s} \mid \ell'\in \H_s \}
\end{equation}
on $V_s$. Due to (\ref {eq:involution:L}) and (\ref {eq:involution:B}), these operators are involutions:
$(E_s L_{as})^2=\Id_{V_s}$, $(E_{s} B_{\ell s})^2=\Id_{V_s}$ and $(H_{s} B_{\ell' s})^2=\Id_{V_s}$, hence they are symplectic (indeed: since $A_\ell|_{V_s}$ is Hamiltonian, $B_{\ell s}$ is also Hamiltonian, hence $\Omega(B_{\ell s}B_{\ell' s}u,B_{\ell s}B_{\ell' s}v)=-\Omega(B_{\ell' s}u,B_{\ell s}^2B_{\ell' s}v)=\Omega(u,B_{\ell s}^2B_{\ell' s}^2v)=\Omega(u,v)$ for any $u,v\in V_s$ and either $\ell,\ell'\in\I_s$ or $\ell,\ell'\in\H_s$; the symplecticity of $L_{as} L_{a's}$ and $L_{as} B_{\ell s}$ is proved similarly).

It follows from Linear Algebra that there exists a unique (up to a permutation of terms) decomposition
$$
V_s = \bigoplus\limits_{t} V_{st}
$$
such that each $V_{st}$ is invariant under each symplectic involution from the set (\ref {eq:involutions}),
\begin{equation} \label{eq:EH:LB}
\begin{array}{l}
E_{s} L_{as}|_{V_{st}} = -\varepsilon_{a st} \Id_{V_{st}}, \quad a \in \M_s, \\
E_{s} B_{\ell s}|_{V_{st}} = -\eta_{\ell st} \Id_{V_{st}}, \quad \ell\in \I_s, \qquad
H_{s} B_{\ell' s}|_{V_{st}} = \eta_{\ell'st} \Id_{V_{st}}, \quad \ell'\in \H_s,
\end{array}
\end{equation}
where $\varepsilon_{a st}, \eta_{\ell st},\eta_{\ell'st}\in\{1,-1\}$, $\eta_{a_s st}=1$, $\eta_{\ell_s st}=1$, $\eta_{\ell_s' st}=1$, and for any $t\ne t'$ there exists either $\ell\in\I_s\cup\H_s$ such that $\eta_{\ell st}\ne\eta_{\ell st'}$ or $a\in\M_s$ such that $\varepsilon_{a st}\ne\varepsilon_{a st'}$.

Clearly each subspace $V_{st}$ is symplectic, $E_{s}$-invariant (if $\I_s\ne\varnothing$ or $\M_s\ne\varnothing$) and $H_{s}$-invariant (if $\H_s\ne\varnothing$).

{\em Substep 1c.} Fix a subspace $V_{st}\subseteq V_s$ from the above decomposition.

We have four possibilities for the subspace $V_s$: it is either of elliptic, hyperbolic or focus-focus type, or trivial.

{\em Case 1 (elliptic):} $\M_s\ne\varnothing$ or $\I_s\ne\varnothing$, moreover $\H_s=\varnothing$.
Recall that, in Substep 1b, we fixed the elliptic Hamiltonian operator $E_s$ on $V_s$, see (\ref {eq:op:ell}) and (\ref {eq:op:ell:}).

We have the Hamiltonian (and symplectic) operator $E_{s}|_{V_{st}}$ such that $(E_{s}|_{V_{st}})^2=-\Id_{V_{st}}$.
It follows from Linear Algebra that there exists a basis
$e_1,\dots,e_{\dim V_{st}/2}$, $f_1,\dots,f_{\dim V_{st}/2}$ of $V_{st}$ such that
$$
E_{s}e_{j}=\hat\eta_{stj} f_{j}, \quad E_{s}f_{j}=-\hat\eta_{stj}e_{j}, \quad \Omega(e_{j},f_{j})=1, \qquad 1\le j\le \frac12 \dim V_{st},
$$
where $\hat\eta_{stj}\in\{1,-1\}$ and the planes $\Span\{e_{j},f_{j}\}$ are pairwise skew-orthogonal.

In the above symplectic basis of the ``elliptic'' subspace $V_{st}$, we have from (\ref {eq:EH:LB}) that
$$
\begin{array}{ccc}
L_{as} e_{j} = - \varepsilon_{a st}\hat\eta_{stj} f_{j}, &
  L_{as} f_{j} = \varepsilon_{a st}\hat\eta_{stj} e_{j}, & \\
B_{\ell s} e_{j} = - \eta_{\ell st}\hat\eta_{stj} f_{j}, &
  B_{\ell s} f_{j} = \eta_{\ell st}\hat\eta_{stj} e_{j}, &
\qquad 1\le j\le \frac12 \dim V_{st},
\end{array}
$$
for any $\ell\in \I_s$ and $a\in\M_s$, therefore, we have from (\ref {eq:L:M}) and (\ref {eq:B:A}) that
\begin{equation} \label {eq:A:ell}
\begin{array}{ll}
M_{a} e_{j} = (\cos(2\pi q_{a,s})) e_j - \varepsilon_{a st}\hat\eta_{stj} (\sin(2\pi q_{a,s})) f_j, &
A_{\ell} e_{j} =  \eta_{\ell st}\hat\eta_{stj} p_{s\ell} f_{j}, \\
M_{a} f_{j} = \varepsilon_{a st}\hat\eta_{stj} (\sin(2\pi q_{a,s})) e_j + (\cos(2\pi q_{a,s})) f_j, &
A_{\ell} f_{j} = -\eta_{\ell st}\hat\eta_{stj} p_{s\ell} e_{j},
\end{array}
\end{equation}
$1\le j\le \frac12 \dim V_{st}$, for any $\ell\in \I_s$ and $a\in\{1,\dots,r\}$, while $A_\ell|_{V_{st}}=0$ for any $\ell\in\{1,\dots,\varkappa_e+\varkappa_h\}\setminus\I_s$.

Finally, let us compute the restriction to $V_{st}$ of the second differential $\ddd^2 J_\ell|_{m_0}(\xi_1,\xi_2)=\Omega(\xi_1,A_\ell\xi_2)$ of each function $J_\ell$ at $m_0$.
It has the form
\begin{equation} \label {eq:J:ell}
\ddd^2J_\ell|_{m_0}(e_i,e_j) = \ddd^2J_\ell|_{m_0}(f_i,f_j)
= \eta_{\ell st} \hat\eta_{stj} p_{s\ell}\delta_{ij}, \qquad
\ddd^2J_\ell|_{m_0}(e_i,f_j) = 0,
\end{equation}
$1\le i,j\le \frac12 \dim V_{st}$, for any $\ell\in\I_s$, while $\ddd^2J_\ell|_{V_{st}}=0$ for any $\ell\in\{1,\dots,\varkappa_e+\varkappa_h\}\setminus\I_s$.

{\em Case 2 (hyperbolic):} $\H_s\ne\varnothing$ and $\I_s=\M_s=\varnothing$.
Thus $\mu_{as}\in\{1,-1\}$ and
\begin{equation} \label {eq:M:hyp}
M_{a}|_{V_s} = \mu_{as} \Id_{V_s}, \qquad 1\le a\le r,
\end{equation}
due to (\ref {eq:q:mu}) and (\ref {eq:L:M}).
Recall that, in Substep 1b, we fixed the hyperbolic Hamiltonian operator $H_s=B_{\ell_s' s}$ on $V_s$, see (\ref {eq:op:hyp}).

Consider the Hamiltonian operator $H_{st}:=H_s|_{V_{st}}=B_{\ell_s' s}|_{V_{st}}$. Since $H_{st}^2=\Id_{V_{st}}$, we have a decomposition $V_{st}=V_{st}'\oplus V_{st}''$ where $V_{st}',V_{st}''$ are the eigenspaces of $H_{st}$ with the eigenvalues $1,-1$ respectively. Since $H_{st}$ is Hamiltonian, we have that each of $V_{st}'$ and $V_{st}''$ is isotropic
(in fact: $0=\Omega(H_{st}u,v)+\Omega(u,H_{st}v)=(\lambda+\mu)\Omega(u,v)$, whenever $u,v\in V_s$ are eigenvectors of $H_{st}$ with eigenvalues $\lambda,\mu$, respectively).
Hence $V_{st}',V_{st}''$ are Lagrangian subspaces of $V_{st}$.

Choose a basis $e_1,e_2,\dots,e_{\dim V_{st}/2}$ of the subspace $V_{st}'$. Since $\Omega|_{V_{st}}$ is nondegenerate and vanishes on each $V_{st}'$ and $V_{st}''$, it gives a nondegenerate pairing between the subspaces $V_{st}'$ and $V_{st}''$, so there exists a unique basis $f_1,f_2,\dots,f_{\dim V_{st}/2}$ of the subspace $V_{st}''$ such that $\Omega(e_i,f_j)=\delta_{ij}$.

For any $\ell\in \H_s$, in the above basis of the ``hyperbolic'' subspace $V_{st}$, we have from (\ref {eq:EH:LB}) that
$$
B_{\ell s} e_j = \eta_{\ell st} H_{st} e_j = \eta_{\ell st} e_j, \qquad
B_{\ell s} f_j = \eta_{\ell st} H_{st} f_j = - \eta_{\ell st} f_j, \qquad
1\le j\le \frac12 \dim V_{st},
$$
$\Omega(e_j,f_j)=1$ and the planes $\Span\{e_j,f_j\}$ are pairwise skew-orthogonal. Thus we have from (\ref {eq:B:A}) and (\ref {eq:M:hyp}) that
\begin{equation} \label {eq:A:hyp}
M_{a} e_j = \mu_{as} e_j,\ M_{a} f_j = \mu_{as} f_j , \qquad
A_{\ell} e_j = \eta_{\ell st} p_{s\ell} e_j, \
A_{\ell} f_j = - \eta_{\ell st} p_{s\ell} f_j
\end{equation}
for any $\ell\in \H_s$ and any $a\in\{1,\dots,r\}$, while $A_\ell|_{V_{st}}=0$ for any $\ell\in\{1,\dots,\varkappa_e+\varkappa_h\}\setminus\H_s$.

Finally, let us compute the restriction to $V_{st}$ of the second differential $\ddd^2J_\ell|_{m_0}(\xi_1,\xi_2)=\Omega(\xi_1,A_\ell\xi_2)$ of each function $J_\ell$ at $m_0$.
It has the form
\begin{equation} \label {eq:J:hyp}
\ddd^2J_\ell|_{m_0}(e_i,f_j)
= - \eta_{\ell st} p_{s\ell} \delta_{ij}, \qquad
\ddd^2J_\ell|_{m_0}(e_i,e_j) = \ddd^2J_\ell|_{m_0}(f_i,f_j) = 0,
\end{equation}
$1\le i,j\le \frac12\dim V_{st}$, for any $\ell\in\H_s$, while $\ddd^2J_\ell|_{V_{st}}=0$ for any $\ell\in\{1,\dots,\varkappa_e+\varkappa_h\}\setminus\H_s$.

{\em Case 3 (focus-focus):} $\I_s\ne\varnothing$ or $\M_s\ne\varnothing$, moreover $\H_s\ne\varnothing$.
Recall that, in Substep 1b, we fixed two Hamiltonian operators $E_s$ and $H_s$ on $V_s$, see (\ref {eq:op:ell}), (\ref {eq:op:ell:}) and (\ref {eq:op:hyp}).

Consider two commuting Hamiltonian operators
$$
E_{st}:=E_{s}|_{V_{st}}, \qquad H_{st}:=H_{s}|_{V_{st}} = B_{\ell_s' s}|_{V_{st}},
$$
which are elliptic and hyperbolic, respectively.

Since $H_{st}^2=\Id_{V_{st}}$, we have a decomposition $V_{st}=V_{st}'\oplus V_{st}''$ where $V_{st}',V_{st}''$ are the eigenspaces of $H_{st}$ with the eigenvalues $1,-1$ respectively. Since $H_{st}$ is a Hamiltonian operator, we have that the subspaces $V_{st}',V_{st}''$ are isotropic (see Case 2), so $V_{st}',V_{st}''$ are Lagrangian subspaces of $V_{st}$.

Since the operator $E_{st}$ commutes with $H_{st}$, it leaves invariant each Lagrangian subspace $V_{st}',V_{st}''$. Since $E_{st}^2=-\Id_{V_{st}}$, it follows from Linear Algebra that $\dim V_{st}'=\dim V_{st}/2$ is even and
there exists a basis $e_1,e_2,\dots,e_{\dim V_{st}'}$ of $V_{st}'$ such that
$$
E_{st}e_{2j-1}=e_{2j}, \quad E_{st}e_{2j}=-e_{2j-1},
\qquad 1\le j\le \frac12 \dim V_{st}'.
$$
Since $\Omega|_{V_{st}}$ is nondegenerate and vanishes on each $V_{st}'$ and $V_{st}''$, it gives a nondegenerate pairing between the subspaces $V_{st}'$ and $V_{st}''$, so there exists a unique basis $f_1,f_2,\dots,f_{\dim V_{st}'}$ of the subspace $V_{st}''$
such that $\Omega(e_i,f_j)=\delta_{ij}$.
In this basis of $V_{st}''$, we can compute the operator $E_{st}|_{V_{st}''}$ as follows. Since $E_{st}$ is a Hamiltonian operator, we have
$$
\Omega(e_{2j},E_{st}f_i) = -\Omega(E_{st}e_{2j},f_i) = \Omega(e_{2j-1},f_i) = \delta_{2j-1,i},
$$
$$
\Omega(e_{2j-1},E_{st}f_i) = -\Omega(E_{st}e_{2j-1},f_i) = -\Omega(e_{2j},f_i) = -\delta_{2j,i},
$$
therefore $E_{st}f_{2j-1}=f_{2j}$ and $E_{st}f_{2j}=-f_{2j-1}$, $1\le j\le \dim V_{st}/4$.

In the above symplectic basis of the ``focus-focus'' subspace $V_{st}$, due to the formulae (\ref{eq:EH:LB}),
each elliptic Hamiltonian operator $L_{as}$, $a\in\M_s$, acts by the formulae
$$
L_{as} e_{2j-1} = - \varepsilon_{a st} e_{2j}, \quad
L_{as} e_{2j} =  \varepsilon_{a st} e_{2j-1}, \quad
L_{as} f_{2j-1} = - \varepsilon_{a st} f_{2j}, \quad
L_{as} f_{2j} =  \varepsilon_{a st} f_{2j-1},
$$
each elliptic Hamiltonian operator $B_{\ell s}$, $\ell\in\I_s$, acts by the formulae
$$
B_{\ell s} e_{2j-1} = - \eta_{\ell st} e_{2j}, \quad
B_{\ell s} e_{2j} =  \eta_{\ell st} e_{2j-1}, \quad
B_{\ell s} f_{2j-1} = - \eta_{\ell st} f_{2j}, \quad
B_{\ell s} f_{2j} =  \eta_{\ell st} f_{2j-1},
$$
$1\le j\le \dim V_{st}/4$, while each hyperbolic Hamiltonian operator $B_{\ell' s}$, $\ell'\in\H_s$, acts by the formulae
$$
B_{\ell' s} e_{j} = \eta_{\ell'st} e_{j}, \quad
B_{\ell' s} f_{j} =  -\eta_{\ell'st} f_{j}, \qquad
\Omega(e_j,f_j)=1, \quad 1\le j\le \frac12 \dim V_{st},
$$
and the planes $\Span\{e_j,f_j\}$ are pairwise skew-orthogonal.
Thus, each elliptic symplectic operator $M_a$, $1\le a\le r$, acts (due to (\ref{eq:L:M})) by the formulae
\begin{equation} \label {eq:M:focus}
\begin{array} {ll}
M_{a} e_{2j-1} = (\cos(2\pi q_{a,s}))e_{2j-1} - \varepsilon_{a st}(\sin(2\pi q_{a,s})) e_{2j}, \\
M_{a} e_{2j} =  \varepsilon_{a st} (\sin(2\pi q_{a,s})) e_{2j-1} + (\cos(2\pi q_{a,s}))e_{2j}, \\
M_{a} f_{2j-1} = (\cos(2\pi q_{a,s}))f_{2j-1} - \varepsilon_{a st}(\sin(2\pi q_{a,s})) f_{2j}, \\
M_{a} f_{2j} =  \varepsilon_{a st}(\sin(2\pi q_{a,s})) f_{2j-1} + (\cos(2\pi q_{a,s})) f_{2j}, & \quad 1\le j\le \frac14 \dim V_{st},
\end{array}
\end{equation}
each elliptic operator $A_\ell$, $\ell\in\I_s$, acts (due to (\ref{eq:B:A})) by the formulae
\begin{equation} \label {eq:A:focus:ell}
\begin{array} {lll}
A_{\ell} e_{2j-1} = -\eta_{\ell st} p_{s\ell} e_{2j}, &
A_{\ell} e_{2j} =  \eta_{\ell st} p_{s\ell} e_{2j-1}, & \\
A_{\ell} f_{2j-1} = -\eta_{\ell st} p_{s\ell} f_{2j}, &
A_{\ell} f_{2j} =  \eta_{\ell st} p_{s\ell} f_{2j-1}, & \quad 1\le j\le \frac14 \dim V_{st},
\end{array}
\end{equation}
each hyperbolic operator $A_{\ell'}$, $\ell'\in\H_s$, acts (due to (\ref{eq:B:A})) by the formulae
\begin{equation} \label {eq:A:focus:hyp}
A_{\ell'} e_j = \eta_{\ell' st} p_{s\ell'} e_j, \quad
A_{\ell'} f_j = - \eta_{\ell' st} p_{s\ell'} f_j, \qquad
1\le j\le \frac12 \dim V_{st},
\end{equation}
while $A_\ell|_{V_{st}}=0$ for all remaining $\ell\in\{1,\dots,\varkappa_e+\varkappa_h\}\setminus(\I_s\cup\H_s)$.

Finally, let us compute the restriction to $V_{st}$ of the second differential $\ddd^2J_\ell|_{m_0}(\xi_1,\xi_2)=\Omega(\xi_1,A_\ell\xi_2)$ of each function $J_\ell$ at $m_0$.
For each ``elliptic'' function $J_\ell$, $\ell\in\I_s$, we have
\begin{equation} \label {eq:J:focus:ell}
\begin{array}{lr}
\ddd^2J_\ell|_{m_0}(e_i,f_{2j-1})
= \eta_{\ell st} p_{s\ell} \delta_{i,2j}, &
\ddd^2J_\ell|_{m_0}(e_i,e_j) = \ddd^2J_\ell|_{m_0}(f_i,f_j) = 0, \\
\ddd^2J_\ell|_{m_0}(e_i,f_{2j})
= - \eta_{\ell st} p_{s\ell} \delta_{i,2j-1}, &
1\le i,j\le \frac12 \dim V_{st}.
\end{array}
\end{equation}
For each ``hyperbolic'' function $J_\ell$, $\ell\in\H_s$, we have
\begin{equation} \label {eq:J:focus:hyp}
\ddd^2J_\ell|_{m_0}(e_i,f_j)
= - \eta_{\ell st} p_{s\ell} \delta_{ij}, \qquad
\ddd^2J_\ell|_{m_0}(e_i,e_j) = \ddd^2J_\ell|_{m_0}(f_i,f_j) = 0,
\end{equation}
$1\le i,j\le \frac12\dim V_{st}$. For all remaining $\ell\in\{1,\dots,\varkappa_e+\varkappa_h\}\setminus(\I_s\cup\H_s)$, we have $\ddd^2J_\ell|_{V_{st}}=0$.

{\em Case 4 (trivial):} $\M_s=\varnothing$ and $\I_s=\H_s=\varnothing$. Thus $V_{st}=V_s$.
Similarly to Case 2, we have $\mu_{as}\in\{1,-1\}$ and (\ref {eq:M:hyp}) (e.g. if $\mu_{as}$ equal 1 for all $a\in\{1,\dots,r\}$, then $M_a|_{V_s}=\Id_{V_s}$ for any $a\in\{1,\dots,r\}$; in fact such $s$ is unique if any, let us denote it by $s_0$).
Moreover $A_\ell|_{V_s}=0$ and $\ddd^2J_\ell|_{V_s}=0$ for any $\ell\in\{1,\dots,\varkappa_e+\varkappa_h\}$.
Choose a symplectic basis $e_1,\dots,e_{\dim V_s/2}$, $f_1,\dots,f_{\dim V_s/2}$ of $V_s$, i.e.\ a basis satisfying
$$
\Omega(e_{i},e_{j})=\Omega(f_{i},f_{j})=0, \quad \Omega(e_{i},f_{j})=\delta_{ij}, \qquad 
1\le i,j\le \frac12 \dim V_s.
$$

{\em Substep 1d.}
Now take a symplectic basis $e_1,\dots,e_n$, $f_1,\dots,f_n$ of $V=T_{m_0}M$ formed by the above symplectic bases of the subspaces $V_{st}$ (see Substep 1c, Cases 1--4). In detail: the basis $e_1,\dots,e_{\dim V_{st}/2}$, $f_1,\dots,f_{\dim V_{st}/2}$ of a subspace $V_{st}$ appears as $e_{n_{st}+1},\dots,e_{n_{st}+\dim V_{st}/2}$, $f_{n_{st}+1},\dots,f_{n_{st}+\dim V_{st}/2}$ in the basis of $V$, for some ``shifting'' integers $n_{st}$ s.t.\ ordering of the subspaces $V_{st}$ with integers $n_{st}$ is consistent with the following partial order: 
``focus-focus'' subspaces,
``elliptic'' ones, 
``hyperbolic'' ones, and
``trivial'' ones.

Let $\hat x_1,\hat y_1,\dots,\hat x_{n},\hat y_n$ be the linear coordinates on $V=T_{m_0}M$ such that
\begin{itemize}
\item the coordinates $\hat x_{2j-1},\hat y_{2j-1}$, $\hat x_{2j},\hat y_{2j}$, $1\le j\le k_f$, correspond to the bases $\frac{e_{2j-1}-f_{2j}}{\sqrt2}$, $\frac{e_{2j}+f_{2j-1}}{\sqrt2}$, $\frac{e_{2j-1}+f_{2j}}{\sqrt2}$, $\frac{-e_{2j}+f_{2j-1}}{\sqrt2}$ of the ``focus-focus'' subspaces $V_{st}$ (if any); they were constructed in Substep 1c, Case 3,
\item the coordinates $\hat x_{j},\hat y_{j}$, $2k_f+1\le j\le 2k_f+k_e$,
correspond to the bases $e_j,f_j$ of the ``elliptic'' subspaces $V_{st}$ (if any); they were constructed in Substep 1c, Case 1,
\item the coordinates $\hat x_{j},\hat y_{j}$, $2k_f+k_e+1\le j\le 2k_f+k_e+k_h$,
correspond to the bases $e_j,f_j$ of the ``hyperbolic'' subspaces $V_{st}$ (if any); they were constructed in Substep 1c, Case 2,
\item the remaining coordinates $\hat x_{j},\hat y_{j}$, $2k_f+k_e+k_h+1\le j\le n$,
correspond to symplectic bases of the ``trivial'' subspaces $V_s$ (if any), see Substep 1c, Case 4.
\end{itemize}
In these coordinates, $\Omega|_{m_0}$, $\ddd^2J_\ell|_{m_0}$ and $M_a$ have the desired form (\ref {eq:hat:omega:})--(\ref {eq:twist:::}), due to
(\ref {eq:A:ell}), (\ref {eq:J:ell}),
(\ref {eq:A:hyp})--(\ref {eq:J:focus:hyp}).

{\em Step 2.}
In Step 1, we actually constructed local coordinates on a neighbourhood of $m_0$, that bring $\Omega$ to the desired canonical form (\ref {eq:linear:coord:omega}) at $m_0$, and bring the functions $J_\ell$ to the desired form (\ref {eq:linear:coord:ell}) and (\ref {eq:linear:coord:hyp}) up to 3rd order terms.
We want to deform these coordinates a little bit, in order to achieve the equalities 
(\ref {eq:linear:coord:omega})--(\ref {eq:linear:coord:hyp}) exactly.

Consider the group
$$
G := \Gamma\times(S^1)^{\varkappa_e+\varkappa_h},
$$
where $(S^1)^{\varkappa_e+\varkappa_h}$ is the $(\varkappa_e+\varkappa_h)$-torus with a usual group structure.
Since $G$ is a compact Lie group acting analytically on a neighbourhood $U^{\mathbb C}\subset M^{\mathbb C}$ of $m_0\in M$ with a fixed point $m_0$, it acts by linear transformations w.r.t.\ to some holomorphic local coordinates on a $G$-invariant neighbourhood $U_0^{\mathbb C}$ of $m_0$ (\cite{boc}, \cite[Sec 3.1.4]{cha}).

Let us explicitely construct (in Substeps 2a--2c below) {\em real-analytic} {\em symplectic} coordinates $x_1,y_1,\dots,x_n,y_n$ on a neighbourhood $U\subset U_0$ of $m_0$ such that $G$ acts by linear transformations on a $G$-invariant neighbourhood $U^{\mathbb C}$ of $m_0$ w.r.t.\ the coordinates $x_1^{\mathbb C},y_1^{\mathbb C},\dots,x_n^{\mathbb C},y_n^{\mathbb C}$.

{\em Substep 2a.} Let $U_0$ be a small neighbourhood of $m_0$ in $M$, and $U_0^{\mathbb C}$ its  small open $G$-invariant complexification.

Take any real-analytic (pseudo-)Riemannian metric $ds^2$ on $U_0$ such that
$$
\ddd s^2|_{m_0}
= \sum\limits_{j=1}^{k_f}
a_j(\ddd\hat x_{2j-1}^2+\ddd\hat y_{2j-1}^2-\ddd\hat x_{2j}^2-\ddd\hat y_{2j}^2)
+\sum\limits_{j=1}^{k_e} a_{k_f+j} (\ddd\hat x_{2k_f+j}^2+\ddd\hat y_{2k_f+j}^2)
+
$$
$$
+\sum\limits_{j=1}^{k_h} a_{k_f+k_e+j} \ddd\hat x_{2k_f+k_e+j} \ddd\hat y_{2k_f+k_e+j}
+\sum\limits_{j=k_e+k_h+2k_f+1}^{n} a_{j-k_f} (\ddd\hat x_{j}^2+\ddd\hat y_{j}^2)
$$
(compare (\ref {eq:d2J:ell}), (\ref {eq:d2J:hyp})). Here $\hat x_j,\hat y_j$ are linear coordinates on $T_{m_0}M$ constructed in Step 1, and $a_1,\dots,a_{n-k_f}\in\R\setminus\{0\}$ are fixed real numbers (e.g.\ $a_j=1$). If $a_j>0$, this pseudo-Riemannian metric has signature $(2n-k_h-2k_f,k_h+2k_f)$, so it is a Riemannian metric if $k_h=k_f=0$ (elliptic case).

Construct a $G$-invariant (pseudo-)Riemannian metric $\langle (\ddd s^2)^{\mathbb C}\rangle$ on $U_0^{\mathbb C}$ by averaging $(\ddd s^2)^{\mathbb C}$ over $G$:
$$
\langle (\ddd s^2)^{\mathbb C}\rangle := \frac1{(2\pi)^\varkappa|\Gamma|}
\sum\limits_{\psi\in\Gamma}\int\limits_0^{2\pi}\dots\int\limits_0^{2\pi}
(\phi_{iJ_\varkappa^{\mathbb C}}^{t_\varkappa} \circ\dots\circ \phi_{iJ_{\varkappa_e+1}^{\mathbb C}}^{t_{\varkappa_e+1}}\circ
\phi_{J_{\varkappa_e}^{\mathbb C}}^{t_{\varkappa_e}}\circ\dots\circ \phi_{J_1^{\mathbb C}}^{t_1}\circ\rho(\psi)
)^*(\ddd s^2)^{\mathbb C} \ddd t_1\dots \ddd t_\varkappa,
$$
where $\varkappa:=\varkappa_e+\varkappa_h$, $\phi_{f}^t$ denotes the flow of the vector field $X_{f}$.
We claim that 
\begin{itemize}
\item[(i)] $\langle (\ddd s^2)^{\mathbb C}\rangle$ is $G$-invariant;
\item[(ii)] $\langle (\ddd s^2)^{\mathbb C}\rangle$ is real-analytic;
\item[(iii)] $\langle (\ddd s^2)^{\mathbb C}\rangle|_{m_0}=(\ddd s^2)^{\mathbb C}|_{m_0}$ and, in particular, it is a nondegenerate quadratic form.
\end{itemize}
The property (i) is obvious.

For proving the property (ii), observe that $\langle (\ddd s^2)^{\mathbb C}\rangle$ can be obtained from $(\ddd s^2)^{\mathbb C}$ in $\varkappa+1$ steps, where the first step performs averaging
over $\Gamma$, and the $(\ell+1)$-st step performs averaging
over the subtorus $\{1\}^{\ell-1}\times S^1\times\{1\}^{\varkappa-\ell}$ of the torus $(S^1)^\varkappa$, $1\le \ell\le \varkappa+1$. After the first $\varkappa_e+1$ steps, the resulting average will be real-analytic, since $\ddd s^2$ is real-analytic and we averaged it over a real-analytic action. The $(\ell+1)$-st step with $\ell>\varkappa_e$ performs an average over the 1-parameter family of diffeomorphisms 
$\phi_{iJ_{\ell}^{\mathbb C}}^t:U^{\mathbb C}\to U^{\mathbb C}$, $0\le t\le 2\pi$, such that
\begin{equation} \label {eq:flow:ell}
\phi_{iJ_{\ell}^{\mathbb C}}^0 = \Id, \qquad
\frac{\ddd}{\ddd t}\phi_{iJ_{\ell}^{\mathbb C}}^t (m)=iX_{J_\ell}^{\mathbb C}|_{\phi_{iJ_{\ell}^{\mathbb C}}^t(m)}.
\end{equation}
This implies that, at each real point $m\in U$, we have
$$
\frac{d}{dt}\overline{\phi_{iJ_{\ell}^{\mathbb C}}^{2\pi-t}(m)}
=i\ \overline{X_{J_\ell}^{\mathbb C}}|_{\phi_{iJ_{\ell}^{\mathbb C}}^{2\pi-t}(m)}
=iX_{J_\ell}^{\mathbb C}|_{\overline{\phi_{iJ_{\ell}^{\mathbb C}}^{2\pi-t}(m)}},
$$
therefore 
$\phi_{iJ_{\ell}^{\mathbb C}}^t(m)=\overline{\phi_{iJ_{\ell}^{\mathbb C}}^{2\pi-t}(m)}$. Hence, averaging over the family of diffeomorphisms $\phi_{iJ_{\ell}^{\mathbb C}}^t$ is the same as averaging over the family of diffeomorphisms $\overline{\phi_{iJ_{\ell}^{\mathbb C}}^t}$.
But, if $\ddd s^2_1$ is a real-analytic (pseudo-)Riemannian metric, then
$(\overline{\phi_{iJ_{\ell}^{\mathbb C}}^t})^*(\ddd s^2_1)^{\mathbb C}
= \overline{(\phi_{iJ_{\ell}^{\mathbb C}}^t)^*(\ddd s^2_1)^{\mathbb C}}$ on $U$.
Hence, the resulting average coincides with its $\mathbb C$-conjugate, which shows that it is real-valued on $U$.

For proving the property (iii), it is enough to check that $(\ddd s^2)^{\mathbb C}|_{m_0}$ is invariant under the linearized $G$-action at $m_0$.
Observe that $\ddd s^2|_{m_0} = \sum\limits_{j=1}^{k_f} a_j \ddd^2h_{2j-1}|_{m_0} + \sum\limits_{j=1}^{n-2k_f} a_{k_f+j} \ddd^2h_{2k_f+j}|_{m_0}$. Here $h_1,\dots,h_{2k_f+k_e+k_h}$ denote the quadratic functions (\ref {eq:mz:hj:}) w.r.t.\ the coordinates $\hat x_j,\hat y_j$, $h_j:=(\hat x_j^2+\hat y_j^2)/2$ for $2k_f+k_e+k_h+1\le j\le n$. But the quadratic forms $\ddd^2h_j|_{m_0}$ on $T_{m_0}M$ pairwise Poisson commute with respect to the symplectic form $\Omega|_{m_0}$. Therefore, each $\ddd^2h_j|_{m_0}$ is invariant under the Hamiltonian flows generated by the quadratic functions (\ref {eq:d2J:ell}) and (\ref {eq:d2J:hyp}), which coincide with the second differentials $\ddd^2J_\ell$ by Step 1. Therefore, each $\ddd^2h_j^{\mathbb C}|_{m_0}$ is invariant under the linearized $(S^1)^\varkappa$-action at $m_0$.
Furthermore, each $\ddd^2h_j|_{m_0}$ is $M_a$-invariant due to (\ref {eq:twist::}) and (\ref {eq:twist:::}) proved in Step 1. Hence each $\ddd^2h_j^{\mathbb C}|_{m_0}$ (and therefore $(\ddd s^2)^{\mathbb C}|_{m_0}$) is invariant under the linearized $G$-action at $m_0$, as required.

By the properties (ii) and (iii), $\langle (ds^2)^{\mathbb C}\rangle$ is the holomorphic extension to $U_0^{\mathbb C}$ of some real-analytic (pseudo-)Riemannian metric on $U_0$, which we denote by $\langle \ddd s^2\rangle$. We have $\langle (\ddd s^2)^{\mathbb C}\rangle=\langle \ddd s^2\rangle^{\mathbb C}$, and it is $G$-invariant due to the property (i).

{\em Substep 2b.}
Let us identify $U_0$ with a small neighbourhood $\widehat U$ of the origin in $T_{m_0}M$ via the exponential map
\begin{equation} \label {eq:exp}
\exp_{m_0}|_{\widehat U}:\widehat U  \stackrel \approx \longrightarrow U_0
\end{equation}
corresponding to the (pseudo-)Riemannian metric $\langle \ddd s^2\rangle$ on $U_0$ constructed in Substep 2a.

Let us transfer the linear coordinates $\hat x_1,\hat y_1,\dots,\hat x_n,\hat y_n$ from $\widehat U$ to $U_0$ by means of the identification $U_0\approx \widehat U$ in (\ref {eq:exp}).
We will obtain some coordinates on $U_0$, which we denote by $u_1,v_1,\dots,u_n,v_n$.

Due to Substep 2a, the holomorphic extension of the map (\ref {eq:exp}) is $G$-equivariant.

{\em Substep 2c.}
Now we want to ``deform'' the coordinates $u_1,v_1,\dots,u_n,v_n$ on $U_0$ (by some $G$-equivariant transformation), in order to make $\Omega$ being constant.
We will achieve this by means of Moser's path method \cite[Theorem 2]{moser}, more precisely its equivariant version, as follows. Let $\hat\Omega$ be the constant symplectic 2-form on $T_{m_0}M$ coinciding with $\Omega|_{m_0}$ at the origin. Denote by $\Omega_0$ the symplectic structure on $U_0$ obtained from $\hat\Omega$ under the identification (\ref {eq:exp}), thus $\Omega_0|_{m_0}=\Omega|_{m_0}$.
Since $\hat\Omega^{\mathbb C}$ is $G$-invariant and the exponential map (\ref {eq:exp})
is $G$-equivariant (by Substep 2b), we conclude that $\Omega_0^{\mathbb C}$ is also $G$-invariant.

Denote
$$
\Omega_t:=(1-t)\Omega_0+t\Omega, \qquad 0\le t\le1.
$$
Choose a real-analytic 1-form $\alpha$ on $U_0$ such that $\alpha^{\mathbb C}$ is $G$-invariant, $\alpha|_{m_0}=0$ (moreover, $\alpha|_m=O(|m-m_0|^2)$ w.r.t.\ some local coordinates on $U_0$; this property will be achieved and used in Step 4 below), and
$$
\ddd\alpha= \Omega-\Omega_0 \equiv \frac{\ddd}{\ddd t}\Omega_t
$$
(such a $G$-invariant 1-form can be obtained by averaging some 1-form with the above properties over $G$, due to the $G$-invariance of $\Omega^{\mathbb C}$ and $\Omega_0^{\mathbb C}$, then the resulting average will be real-analytic, due to the same arguments as in Substep 2a).
Define a 1-parameter family of vector fields $\xi_t$ on $U_0$, $0\le t\le 1$, by Moser's equation
$$
i_{\xi_t}\Omega_t=-\alpha
$$
(such a vector field exists and is unique, perhaps in a smaller neighbourhood of $m_0$, since $\Omega_t|_{m_0}=\Omega|_{m_0}$ is a nondegenerate 2-form).
Define a 1-parameter family of diffeomorphisms $\phi_t:U_t\to U_0$, $0\le t\le1$, such that $\phi_0=\Id_{U_0}$ and
$$
\frac{\ddd}{\ddd t}\phi_t = \xi_t \circ \phi_t
$$
in $U_t$.
Observe that $m_0$ is fixed under each $\phi_t$, since $\xi_t|_{m_0}=0$ because $\alpha|_{m_0}=0$.
Clearly, $\Omega_t^{\mathbb C}$, $\xi_t^{\mathbb C}$ and $U_t^{\mathbb C}$ are $G$-invariant, thus each $\phi_t^{\mathbb C}$ is $G$-equivariant. We have
$$
\frac{\ddd}{\ddd t}(\phi_t^*\Omega_t) = \phi_t^*(\Lie_{\xi_t}\Omega_t + \frac{\ddd\Omega_t}{\ddd t})
= \phi_t^*((i_{\xi_t}\ddd+\ddd i_{\xi_t})\Omega_t + \Omega-\Omega_0)
= \phi_t^*(\ddd i_{\xi_t}\Omega_t + \Omega-\Omega_0)
= \phi_t^*(-\ddd\alpha + \Omega-\Omega_0)
= 0.
$$
Since the above equalities hold for any $t\in[0,1]$, we conclude that $\phi_t^*\Omega_t\equiv\phi_0^*\Omega_0=\Omega_0$ for any $t$, in particular for $t=1$ we obtain
$$
\phi_1^*\Omega=\Omega_0.
$$
Define on $U:=\phi_1(U_1\cap U_0)$ the coordinates
$$
x_j:=u_j\circ\phi_1^{-1}|_U, \quad y_j:=v_j\circ\phi_1^{-1}|_U, \qquad 1\le j\le n,
$$
i.e.\ $x_1,y_1,\dots,x_n,y_n$ are induced from the coordinates $u_1,v_1,\dots,u_n,v_n$ on $U_0':=U_1\cap U_0$ by means of the identification $U_0'\stackrel\approx\longrightarrow U$ via the diffeomorphism $\phi_1:U_1\to U_0$. Clearly, $x_j(m_0)=y_j(m_0)=0$ for any $j=1,\dots,n$ (since $m_0$ is fixed under each $\phi_t$).

{\em Step 3.}
We claim that
\begin{equation} \label {eq:quadratic::}
\Omega|_{U}=\sum\limits_{j=1}^n\ddd x_{j}\wedge \ddd y_{j}, \qquad
J_\ell|_{U}=c_\ell+\sum\limits_{a,b=1}^{2n}c_{ab}^\ell w_a w_b,
\end{equation}
for some real constants $c_\ell,c_{ab}^\ell$. Here we denoted $(w_1,\dots,w_{2n}):=(x_1,y_1,\dots,x_n,y_n)$.

Recall (see Substep 2b and the beginning of Substep 2c) that
$$
\Omega_0=((\exp_{m_0}|_{\widehat U})^{-1})^*\hat\Omega, \qquad
u_j:=\hat x_j\circ(\exp_{m_0}|_{\widehat U})^{-1}, \quad
v_j:=\hat y_j\circ(\exp_{m_0}|_{\widehat U})^{-1}.
$$
Thus
$$
\Omega|_{U}=((\phi_1\circ\exp_{m_0}|_{\widehat U'})^{-1})^*\hat\Omega, \quad
x_j=\hat x_j\circ(\phi_1\circ\exp_{m_0}|_{\widehat U'})^{-1}, \
y_j=\hat y_j\circ(\phi_1\circ\exp_{m_0}|_{\widehat U'})^{-1},
$$
where $\widehat U':=(\exp_{m_0}|_{\widehat U})^{-1}(U_0')$.
We conclude from (\ref {eq:hat:omega:}) and from the beginning of Substep 2c that $\Omega|_{U}$ has the desired form as in (\ref {eq:quadratic::}).

Observe that the above diffeomorphism
\begin{equation} \label {eq:change}
\phi_1\circ\exp_{m_0}|_{\widehat U'}:\widehat U' \stackrel \approx \longrightarrow U
\end{equation}
has a $G$-equivariant holomorphic extension, since $\phi_1^{\mathbb C}$ and $\exp_{m_0}^{\mathbb C}$ are $G$-equivariant.
Since the coordinates $w_1,\dots,w_{2n}$ on $U$ are induced from the linear coordinates $\hat x_1,\hat y_1,\dots,\hat x_n,\hat y_n$ on $T_{m_0}M$ under the diffeomorphism $U\approx\widehat U'\subset T_{m_0}M$ in (\ref {eq:change}), moreover the group $G$ acts on $T_{m_0}M^{\mathbb C}$ by linear transformations w.r.t.\ $\hat x_1^{\mathbb C},\hat y_1^{\mathbb C},\dots,\hat x_n^{\mathbb C},\hat y_n^{\mathbb C}$, we conclude that it acts on $U^{\mathbb C}$ by linear transformations w.r.t.\ $w_1^{\mathbb C},\dots,w_{2n}^{\mathbb C}$.
Therefore the Hamiltonian vector fields $X_{J_\ell}$ on $U$ (which generate an infinitesimal action of the Lie algebra of $G$) are also linear w.r.t.\ $w_1,\dots,w_{2n}$, $1\le\ell\le k$. Since $\Omega|_{U}$ has constant components w.r.t.\ $w_1,\dots,w_{2n}$, as in (\ref {eq:quadratic::}),
we conclude that $\ddd J_\ell|_{U}$ has linear components w.r.t.\ $w_1,\dots,w_{2n}$, i.e.\
$\ddd J_\ell|_{U}=2\sum\limits_{a,b=1}^n c_{ab}^\ell w_a \ddd w_b$ for some constants $c_{ab}^\ell=c_{ba}^\ell\in\R$.
This immediately gives us the desired quadratic form for $J_\ell|_{U}$ as in (\ref {eq:quadratic::}).

{\em Step 4.} Let us show that the quadratic forms in (\ref {eq:quadratic::}) have the special form as in (\ref {eq:linear:coord:ell}), (\ref {eq:linear:coord:hyp}).
For this, it is enough to make sure that $\ddd x_j|_{m_0}=\ddd\hat x_j|_{m_0}$ and $\ddd y_j|_{m_0}=\ddd\hat y_j|_{m_0}$, $1\le j\le n$, due to (\ref {eq:d2J:ell}) and (\ref {eq:d2J:hyp}).
So, it is enough to achieve that $\ddd\phi_t|_{m_0}=\Id$ for each $t\in[0,1]$.
The latter is equivalent to the equality $\ddd\xi_t|_{m_0}=0$, which in turn is equivalent to the fact that each component of $\alpha|_m$ w.r.t.\ some coordinates on $U$ is of order $O(|m-m_0|^2)$.

So, it is enough to show that, in Substep 2c, we can choose a 1-form $\alpha$ on $U_0$ satisfying the above condition.
We have $\Omega|_{U_0}-\Omega_0=\ddd(\alpha_1-\phi^*\alpha_1)$, where $\alpha_1$ is any 1-form on $U_0$ such that $\ddd\alpha_1=\Omega$,
and $\phi$ is a diffeomorphism such that $\phi(m_0)=m_0$ and $\ddd\phi|_{m_0}=\Id$.
In canonical coordinates $(p,q)=(p_1,q_1,\dots,p_n,q_n)$ such that $p_j|_{m_0}=q_j|_{m_0}=0$, we have $\Omega|_{U_0}=\ddd p\wedge \ddd q:=\sum\limits_{j=1}^n \ddd p_j\wedge \ddd q_j$.
Put $\alpha_1=p\ddd q:=\sum\limits_{j=1}^np_j\ddd q_j$, $\alpha:=\alpha_1-\phi^*\alpha_1$.
Since $\phi(m_0)=m_0$ and $\ddd\phi|_{m_0}=\Id$, it follows from Hadamard's lemma that $\phi^*p_i=p_i+Q_{i}'(p,q)$, $\phi^*q_j=q_j+Q_{j}''(p,q)$, where $Q_i'(p,q)$ and $Q_j''(p,q)$ are some quadratic forms whose coefficients are real-analytic functions in $(p,q)$. Thus
$\alpha=\alpha_1-\phi^*\alpha_1=p\ddd q-\phi^*(p\ddd q)=p\ddd q-(\phi^*p)\ddd(\phi^*q)=\sum\limits_{j=1}^n p_j\ddd q_j-(p_j+Q_j'(p,q))\ddd(q_j+Q_j''(p,q)) = -\sum\limits_{j=1}^n (p_j+Q_j'(p,q))\ddd Q_j''(p,q) - Q_j'(p,q)\ddd q_j$. Each component of the latter 1-form is of order $O(|p|^2+|q|^2)$, as required.

Lemma \ref {lem:period:} (a), and hence Lemma \ref {lem:period} (a), is proved.
\qed

\subsection{Proof of part (b) of Lemmata \ref {lem:period} and \ref {lem:period:}} \label {subsec:proof:lem:period:}

Similarly to the previous subsection, we will give a proof of Lemma \ref {lem:period:} (b).
Due to the Cauchy integral formula for holomorphic functions, we can and will assume that 
$\|\tilde F^{\mathbb C}-F^{\mathbb C}\|_{C^k}
+\|\tilde \Omega^{\mathbb C}-\Omega^{\mathbb C}\|_{C^{k+n-2}}
+\|\tilde \rho^{\mathbb C}-\rho^{\mathbb C}\|_{C^{k-1}}<\varepsilon$ for some $k\ge4$
(one has similar inequalities for the real objects on $M$ in assumptions of Lemma \ref {lem:period} (b)).
If $\varkappa_h=0$ (i.e., all functions $J_\ell$ generate $2\pi$-periodic flows), all our arguments and constructions will literally work both in the real-analytic and $C^\infty$ cases. This will give a proof of Lemma \ref {lem:period} (b) too.

Consider the subgroup $G_0:=\Gamma\times(S^1)^{\varkappa_e}\times\{1\}^{\varkappa_h}$ of the group $G=\Gamma\times(S^1)^{\varkappa_e+\varkappa_h}$.
Recall that $U$ is $G_0$-invariant, and $U^{\mathbb C}$ is $G$-invariant.

We will identify $U$ with its image under the coordinate map $(x_1,y_1,\dots,x_n,y_n):U\to\R^{2n}$ from (a).
We will equip $\R^{2n}$ with cartesian coordinates $(x_1,y_1,\dots,x_n,y_n)$ and with the standard symplectic structure $\Omega$ as in (\ref {eq:quadratic::}).

Consider the ``perturbed'' $G$-action generated by $\tilde\rho(\Gamma)$ and by the ``perturbed'' functions $\tilde J_1,\dots,\tilde J_{\varkappa_e},i\tilde J_{\varkappa_e+1},\dots,i\tilde J_{\varkappa_e+\varkappa_h}$ w.r.t.\ the ``perturbed'' symplectic structure $\tilde\Omega$. By abusing language, we will call this action the {\em $\tilde G$-action}, in order to distinguish it from the ``unperturbed'' $G$-action.

Divide the proof into several steps. 

{\em Step 1.}
In this step, we prove that the $\tilde G$-action has a fixed point $\tilde m_0\in U'$ $O(\varepsilon)-$close to $m_0$.
We will use the same notations (e.g.\ $A_\ell, M_a, V_s$) as in the proof of (a).

{\em Substep 1a.}
Denote
$$
R:=\bigcap_{a=1}^r \ker (M_a - \Id), \qquad
K:= R \cap \bigcap_{\ell=1}^{\varkappa_e+\varkappa_h} \ker A_\ell 
= \bigcap_{\ell=1}^{\varkappa_e+\varkappa_h} \ker A_\ell|_R
$$
where $\Id$ is the identity operator in $\R^{2n}$ (here we used that $M_a$ and $A_\ell$ commute, hence the set $R$ of fixed points of $M_a$ is $A_\ell$-invariant).
By (a), the group $G$ acts on $U$ linearly w.r.t.\ our coordinates $(x_1,y_1,\dots,x_n,y_n)$. We conclude that
$R\cap U$ is the set of fixed points of the (``unperturbed'') $\Gamma$-action,
 and $(K\cap U)^{\mathbb C}$ is the set of fixed points of the (``unperturbed'') $G$-action.

Clearly, there exists a linear combination
$$
A=\sum\limits_{\ell=1}^{\varkappa_e+\varkappa_h} c_\ell A_\ell
$$
such that $R\cap \ker A=K$.

We claim that $K$ is a symplectic subspace of $\R^{2n}$.
Indeed, each $\ker A_\ell$ and $\ker(M_a-\Id)$ is the direct product of several subspaces $V_s$. Hence $K$ is also the direct product of some subspaces $V_s$ (actually, it is one of the ``trivial'' subspaces $V_s$, see the proof of (a), Case 4 of Substep 1c), therefore $K$ is symplectic.

{\em Substep 1b.} Since $K$ is symplectic (by Substep 1a), we have
$$
\R^{2n}=W\times K
$$
where $W:=K^\perp$ denotes the skew-orthogonal complement of $K$ w.r.t.\ $\Omega$.
In fact, $W$ and $K$ are the coordinate subspaces of $\R^{2n}$:
$$
\R^{2n} = \R^{2k}\times\R^{2(n-k)} = W \times K
$$
where $2k_f+k_e+k_h\le k\le n$, with coordinates 
$$
\sigma := (x_1,y_1,\dots,x_k,y_k)
 \quad \mbox{on $W$}, \qquad \quad 
\tau := (x_{k+1},y_{k+1},\dots,x_n,y_n) \quad \mbox{on $K$}.
$$
Denote by
$$
Pr_K:\R^{2n} = W\times K \to K, \qquad (\sigma,\tau) \mapsto \tau,
$$
the projection along $W$. Clearly, the subspace $W^{\mathbb C}$ and the map $Pr_K^{\mathbb C}$ are $G$-invariant, furthermore every point of $K$ is fixed under the $G$-action.

Choose
$G_0$-invariant neighbourhoods
$$
U':=U_W'\times U_{K}'\subset U_W\times U_{K}\subset W\times K
$$
of $m_0=0$ in $\R^{2n}$ such that $\overline{U'}\subset U_W\times U_K$ and $\overline{U_W\times U_K}\subset U$.
Then we can choose $G$-invariant neighbourhoods
$$
(U')^{\mathbb C}:=(U_W')^{\mathbb C}\times (U_{K}')^{\mathbb C}\subset U_W^{\mathbb C}\times U_{K}^{\mathbb C}\subset W^{\mathbb C}\times K^{\mathbb C}
$$
of $m_0=0$ in ${\mathbb C}^{2n}$ such that $\overline{(U')^{\mathbb C}}\subset U_W^{\mathbb C}\times U_K^{\mathbb C}$ and $\overline{U_W^{\mathbb C}\times U_K^{\mathbb C}}\subset U^{\mathbb C}$.

Consider the (``perturbed'') $\tilde G_0$-invariant neighbourhood $\tilde U:=\tilde G_0(U_W\times U_{K})$ of $m_0$ in $\R^{2n}$
and the (``perturbed'') $\tilde G$-invariant neighbourhood ${\tilde U}^{\mathbb C}:=\tilde G(U_W^{\mathbb C}\times U_{K}^{\mathbb C})$ of $m_0$ in ${\mathbb C}^{2n}$.
Clearly, $U'\subset\tilde U\subset U$ and ${U'}^{\mathbb C}\subset{\tilde U}^{\mathbb C}\subset U^{\mathbb C}$ if the perturbation is small enough.

Denote by
$$
\widetilde{Pr_K^{\mathbb C}}: (U_W')^{\mathbb C}\times U_K^{\mathbb C}\to K^{\mathbb C}
$$
the map obtained by averaging the map $Pr_K^{\mathbb C}|_{\tilde U^{\mathbb C}}$ over the $\tilde G$-action:
$$
\widetilde{Pr_K^{\mathbb C}}(m) := \sum\limits_{\psi\in\Gamma}\frac1{(2\pi)^\varkappa|\Gamma|}
\int\limits_0^{2\pi}\dots\int\limits_0^{2\pi}
Pr_K^{\mathbb C}\circ
\phi_{i\tilde J_\varkappa^{\mathbb C}}^{t_\varkappa}
\circ\dots\circ
\phi_{i\tilde J_{\varkappa_e+1}^{\mathbb C}}^{t_{\varkappa_e+1}} \circ
\phi_{\tilde J_{\varkappa_e}^{\mathbb C}}^{t_{\varkappa+e}}
\circ\dots\circ
\phi_{\tilde J_1^{\mathbb C}}^{t_1}
\circ \tilde\rho(\psi)
(m)
\ddd t_1\dots \ddd t_\varkappa,
$$
$m\in (U_W')^{\mathbb C}\times U_K^{\mathbb C}\subset (\tilde U)^{\mathbb C}\subset{\mathbb C}^{2n}$, where $\varkappa:=\varkappa_e+\varkappa_h$, $\phi_{f}^t$ denotes the flow of the vector field $X_{f}$;
we use the identification of $U$ with its image under the coordinate map $(x_1,y_1,\dots,x_{n},y_n):U\to\R^{2n}$ from (a).
Clearly, $\widetilde{Pr_K^{\mathbb C}}$ is $\tilde G$-invariant and $O(\varepsilon)$-close to $Pr_K^{\mathbb C}|_{(U_W')^{\mathbb C}\times U_K^{\mathbb C}}$ in $C^{k-1}-$norm (since the $G$-action on $W\times K$ is componentwise, and the $G$-action on $K$ is trivial), thus it is a $\tilde G$-invariant submersion.
The map $\widetilde{Pr_K^{\mathbb C}}$ is real-valued on $U_W'\times U_K$ by the same arguments as in Step 1 of the proof of (a). Thus, $\widetilde{Pr_K^{\mathbb C}}=\tilde Pr_K^{\mathbb C}$ for some real-analytic map (submersion)
$$
\tilde Pr_K: U_W' \times U_K\to K.
$$

Due to the Inverse Functions Theorem, for each point $\tau\in U_K'$, its pre-image
$$
\tilde W_{\tau}:= \tilde Pr_K^{-1}(\tau) \subset U_W' \times U_K \subset W\times K
$$
is a real-analytic submanifold $O(\varepsilon)$-close to the submanifold $W_{\tau}:=U_{W}'\times \{\tau\}$ in $C^{k-1}-$norm. Clearly $\tilde W_{\tau}^{\mathbb C}$ is $\tilde G$-invariant (since $\tilde Pr_K^{\mathbb C}$ is $\tilde G$-invariant). Hence the submanifold $\tilde W_{\tau}$ is symplectic and has the form
$$
\tilde W_{\tau} = \{(\sigma,\tilde T(\sigma,\tau)) \mid \sigma\in U_{W}'\}  \subset
W \times K, \qquad \tau\in U_K',
$$
for some real-analytic map $\tilde T: U'\to U_K$ close to $Pr_K|_{U'}$ in $C^{k-1}-$norm.

Denote $\tilde U':=\tilde Pr_K^{-1}(U_K')$, so $\tilde U'=\bigcup\limits_{\tau\in U_K'}\tilde W_{\tau}$ is close to $U'$.

{\em Substep 1c.}
Choose $a\in\{1,\dots,r\}$ and consider the cyclic subgroup $\Gamma_a$ of $\Gamma$ generated by $\psi^a$ (recall that $\psi^a\in\Gamma$ are the generators of the group $\Gamma$, as in the proof of (a), Step 1). Denote $R_a:=\ker(M_a-\Id)$, so $R_a\cap U$ is the set of all fixed points of the (``unperturbed'') $\Gamma_a$-action.
Clearly, $R_a$ is the direct product of several coordinate subspaces $Ox_jy_j$, so it is symplectic and has the form $R_a=W_a\times K$, where $W_a$ is the skew-orthogonal complement of $K$ in $R_a$. Denote the skew-orthogonal complement of $R_a$ by $W^a:=R_a^\perp$, so it is symplectic too. We have 
$$
\R^{2n} = W \times K = W^a \times (W_a \times K) = W^a \times R_a , 
$$
with coordinates
$$
\sigma^a \quad \mbox{on $W^a$}, \qquad \qquad 
\sigma_a \quad \mbox{on $W_a$}.
$$

By applying the arguments from Substep 1b to the $\Gamma_a$-action and the (``perturbed'') $\tilde\Gamma_a$-action, we construct a real-analytic $\tilde\Gamma_a$-invariant map (submersion)
$$
\tilde Pr_{R_a}:U'_{W^a}\times U_{R_a}\to R_a
$$
close in $C^{k-1}-$norm to the projection $Pr_{R_a}$ along $R_a$ restricted to $U'_{W^a}\times U_{R_a}$, and 
a family of real-analytic $\tilde\Gamma_a$-invariant submanifolds 
$$
\tilde W^a_{\sigma_a,\tau}:=\tilde Pr_{R_a}^{-1}(\sigma_a,\tau) \subset U_{W^a}' \times U_{R_a} \subset W^a\times R_a
$$
close to the submanifold $W^a_{\sigma_a,\tau}:=U_{W^a}'\times \{(\sigma_a,\tau)\}$ in $C^{k-1}-$norm, $(\sigma_a,\tau)\in U'_{R_a}=U'_{W_a}\times U'_K$. This submanifold has the form
$$
\tilde W^a_{\sigma_a,\tau} = \{ ( \sigma^a,\tilde S_a(\sigma^a,\sigma_a,\tau),\tilde T_a(\sigma^a,\sigma_a,\tau) ) \mid \sigma^a\in U_{W^a}'\} \subset
W^a \times R_a, \qquad (\sigma_a,\tau)\in U_{R_a}'=U'_{W_a}\times U'_K,
$$
for some real-analytic maps $\tilde S_a: U'\to U_{W_a}$ 
and $\tilde T_a: U'\to U_K$ 
close in $C^{k-1}-$norm to $Pr_{W_a}|_{U'}$, $(\sigma^a,\sigma_a,\tau)\mapsto\sigma_a$, and $Pr_{K}|_{U'}$,  respectively.

Observe that $\ker((M_a-\Id)|_{W^a})=R_a\cap W^a=\{0\}$, 
so $0$ is a unique fixed point of the (``unperturbed'') map $M_a|_{W^a}$.
It follows from the Inverse Functions Theorem that the (``perturbed'') map
 $\tilde\rho(\psi^a)|_{\tilde W^a_{\sigma_a,\tau}}$ has a unique fixed point 
$$
\tilde m^a_{\sigma_a,\tau} = ( \tilde S^a(\sigma_a,\tau),\tilde S_a(\tilde S^a(\sigma_a,\tau),\sigma_a,\tau),\tilde T_a(\tilde S^a(\sigma_a,\tau),\sigma_a,\tau) )
 \in \tilde W^a_{\sigma_a,\tau},
$$
for some real-analytic function $\tilde S^a: U_{R_a}'\to U_{W^a}'$ close to $S^a\equiv 0$ in $C^{k-1}-$norm. Hence, 
$$
\tilde R_a := \{\tilde m^a_{\sigma_a,\tau}\mid (\sigma_a,\tau)\in U_{R_a}'\}
$$
is the set of all fixed points in $\tilde Pr_{R_a}^{-1}(U'_{R_a})$ of the $\tilde\Gamma^a$-action. Clearly, $\tilde R_a$ is a real-analytic submanifold close to $U_{R_a}' = R_a\cap U'$ in $C^{k-1}-$norm (and, hence, it is symplectic).

Since the group $\Gamma$ is Abelian, it follows that, for any nonempty subset $\M\subseteq\{1,\dots,r\}$, the set 
$\bigcap\limits_{a\in\M}\tilde R_a$ is a $(\dim(\bigcap\limits_{a\in\M}R_a))$-dimensional submanifold close to $(\bigcap\limits_{a\in\M}R_a)\cap U'$.
Thus
$$
\tilde R := \tilde R_1\cap\dots\cap\tilde R_r
$$
is a $(\dim R)$-dimensional submanifold close to $R\cap U'=U_R'$ in $C^{k-1}-$norm, where $R$ is the same as in Substep 1a. In particular, it is symplectic.
Clearly, $\tilde R$ is nothing else than the set of all fixed points in $\tilde Pr_R^{-1}(U'_R)$ of the $\tilde\Gamma$-action.

{\em Substep 1d.}
Due to Substep 1c, the symplectic submanifold $\tilde R$ is the set of all fixed points  in $\tilde Pr_R^{-1}(U'_R)$ of the $\tilde\Gamma$-action. Hence, it is invariant under the (local) Hamiltonian flows generated by $\tilde J_1,\dots,\tilde J_\varkappa$ (since the $\tilde\Gamma$-action commutes with these flows).

Since $R$ is symplectic and contains a symplectic subspace $K$, it has the form 
$$
R=R'\times K,
$$
where $R'$ is the skew-orthogonal complement of $K$ in $R$ (and, hence, symplectic too).
Moreover, its skew-orthogonal complement $W':=(R')^\perp$ in $W$ is also the direct product of several coordinate subspaces $Ox_jy_j$, so it is symplectic too. We have 
$$
\R^{2n} = W \times K = W' \times (R' \times K) = W' \times R , 
$$
with coordinates
$$
\sigma'  \quad \mbox{on $W'$}, \qquad \qquad 
\sigma'' \quad \mbox{on $R'$}.
$$
Thus, the submanifold $\tilde R$ has the form
$$
\tilde R = \{(\tilde S'(\sigma'',\tau),\sigma'',\tau) \mid (\sigma'',\tau)\in R\cap U'\},
$$
for some real-analytic map $\tilde S':R\cap U' \to W'\cap U$ close to $S'\equiv 0$ in $C^{k-1}-$norm.

By applying the arguments from Substep 1b to the $G$-action on $(R\cap U')^{\mathbb C}$ and to the (``perturbed'') $\tilde G$-action on $\tilde R^{\mathbb C}$, we can construct a real-analytic map (submersion)
$$
\tilde Pr_K': \tilde R\to K
$$
close in $C^{k-1}-$norm to the projection $Pr_K':=Pr_K|_{R\cap U'}$ along $R'$, such that $(\tilde Pr_K')^{\mathbb C}: \tilde R^{\mathbb C}\to K^{\mathbb C}$ is $\tilde G$-invariant. Further, we construct a family of real-analytic submanifolds
$$
\tilde R'_{\tau}:=(\tilde Pr_K')^{-1}(\tau) \subset \tilde R
$$
close to the submanifold $R'_{\tau}:=(R'\cap U')\times \{\tau\} \subset R\cap U'$, $\tau\in U'_K$. Clearly, the submanifolds $(\tilde R'_{\tau})^{\mathbb C}$ are $\tilde G$-invariant, $\tau\in U'_K$.
In particular, each of these submanifolds is symplectic and has the form
$$
\tilde R'_{\tau} = \{(\tilde S'(\sigma'',\tilde T'(\sigma'',\tau)),\sigma'',\tilde T'(\sigma'',\tau)) \mid \sigma'' \in R'\cap U' \} \subset
W' \times R, \qquad \tau\in U'_K,
$$
for some real-analytic map $\tilde T': R\cap U'\to U_K$ close to $Pr_{K}|_{R\cap U'}$ in $C^{k-1}-$norm.

Consider the ``unperturbed'' and the ``perturbed'' functions
$$
f:=\sum\limits_{\ell=1}^\varkappa c_\ell J_\ell, \qquad \tilde f:=\sum\limits_{\ell=1}^\varkappa c_\ell \tilde J_\ell,
$$
see Substep 1a. It follows from Substep 1a that the ``unperturbed'' function $f|_{R'_\tau}$ is quadratic and has a unique critical point, namely at the ``origin'' $(\sigma',\sigma'',\tau)=(0,0,\tau)\in R'_\tau$, $\tau\in U_K'$. 
Since $\tilde R'_{\tau}$ is close to $R'_{\tau}=(R'\cap U')\times \{\tau\}$, and (due to the Inverse Functions Theorem) after a small perturbation, nondegenerate critical points are preserved and deformed a little bit, we conclude that the ``perturbed'' function $\tilde f|_{\tilde R'_\tau}$ has a unique critical point 
$$
\tilde m_{\tau}
= (\tilde S'(\tilde S''(\tau),\tilde T'(\tilde S''(\tau),\tau)),\tilde S''(\tau),\tilde T'(\tilde S''(\tau),\tau)), \qquad \tau\in U_K',
$$
for some real-analytic function $\tilde S'':U_K' \to R'$ close to $S''\equiv 0$ in $C^{k-1}-$norm.

We claim that the point $\tilde m_{\tau}$ is fixed under the $\tilde G$-action, $\tau\in U_K'$. Indeed, the function $(\tilde f|_{\tilde R'_\tau})^{\mathbb C}$ is $\tilde G$-invariant, hence its (unique) critical point $\tilde m_{\tau}$ is also $\tilde G$-invariant.

Thus, the intersection of $U'$ with the fixed points set of the $\tilde G$-action on $(U')^{\mathbb C}$ coincides with
$$
\tilde K := \{\tilde m_{\tau} \mid \tau\in U_K' \}
= \{ (\tilde S'(\tilde S''(\tau),\tilde T'(\tilde S''(\tau),\tau)),\tilde S''(\tau),\tilde T'(\tilde S''(\tau),\tau)) \mid \tau\in U_K'\}
\subset \tilde U' ,
$$
which is a $(\dim K)$-dimensional submanifold close to $\{0\}\times U_K'$.
Consider the point
$$
\tilde m_0 := (\tilde S'(\tilde S''(0),\tilde T'(\tilde S''(0),0)),\tilde S''(0),\tilde T'(\tilde S''(0),0) )\in \tilde K
$$
of $\tilde K$ corresponding to the origin $\tau=0$ of $K$.
It is a desired fixed point of the $\tilde G$-action.

{\em Step 2.} 
Now we will apply (a) to the ``perturbed'' $\tilde G$-action and its fixed point $\tilde m_0$.

It follows from the proof of (a), Step 1, that one can provide an algorithm for constructing a symplectic basis $e_1,f_1,\dots,e_n,f_n$ of $V:=T_{m_0}M$ in which $\Omega|_{m_0}$, $\ddd^2J_\ell|_{m_0}$ and $M_a$ have the form (\ref {eq:hat:omega:})--(\ref {eq:twist:::}).
Let us perform the same algorithm for constructing a symplectic basis of $\tilde V:=T_{\tilde m_0}M$.
Similarly to the proof of (a), Step 1, denote by $\tilde A_\ell$ the linearization of the vector field $X_{\tilde J_\ell}$ at the equilibrium point $\tilde m_0$, $1\le \ell\le \varkappa_e+\varkappa_h$. Furthermore, denote by $\tilde M_a$ the linearization of the symplectomorphism $\tilde\rho(\psi^a)$ at the fixed point $\tilde m_0$, $1\le a\le r$ (where $\psi^a\in\Gamma$ are the generators of the group $\Gamma$).

Since the ``perturbed'' operators $\tilde A_\ell$ and $\tilde M_a$ pairwise commute and are close to the ``unperturbed'' operators $A_\ell$ and $M_a$ (respectively), it follows from Linear Algebra that there exists a unique decomposition
$$
\R^{2n}=\bigoplus\limits_{s}\tilde V_s
$$
such that each $\tilde V_s$ is close to $V_s$ and is invariant under each $\tilde A_\ell$ and each $\tilde M_a$. Clearly, the subspaces $\tilde V_s$ are symplectic and pairwise skew-orthogonal.

Since the ``perturbed'' flow of each ``perturbed'' vector field $X_{\tilde J_\ell}$, $1\le \ell\le\varkappa_e$, and $X_{i\tilde J_\ell}$, $\varkappa_e+1\le \ell\le \varkappa_e+\varkappa_h$, is $2\pi$-periodic, the operator $\tilde A_\ell$ is diagonalizable over $\mathbb C$ and each its eigenvalue belongs either to $i\mathbb Z$ for $1\le \ell\le\varkappa_e$, or to $\mathbb Z$ for $\varkappa_e+1\le \ell\le \varkappa_e+\varkappa_h$. Since the eigenvalues of $\tilde A_\ell|_{\tilde V_s}$ are close to the eigenvalues of $A_\ell|_{V_s}$ and belong to ${\mathbb Z}\cup i\mathbb Z$, we conclude that $\spec(\tilde A_\ell|_{\tilde V_s})=\spec(A_\ell|_{V_s})=\{\pm\lambda_{\ell s},\pm\overline{\lambda_{\ell s}}\}$. 
Thus $\tilde A_\ell$ and $A_\ell$ have the same spectrum. 
Similarly, since the ``perturbed'' operators $\tilde M_a$, $1\le a\le r$, are of finite order, the operator $\tilde M_a$ is diagonalizable over $\mathbb C$ and each its eigenvalue belongs to the unit circle of $\mathbb C$. Since the eigenvalues of $\tilde M_a|_{\tilde V_s}$ are close to the eigenvalues of $M_a|_{V_s}$ and belong to the unit circle in $\mathbb C$, moreover $\tilde M_a^{|\Gamma|}=\Id$, we conclude that $\spec(\tilde M_a|_{\tilde V_s})=\spec(M_a|_{V_s})=\{\mu_{as}^{\pm1},\overline{\mu_{as}^{\pm1}}\}$. 
Thus $\tilde M_a$ and $M_a$ have the same spectrum.

For each subspace $\tilde V_s$ from the above decomposition, we define (similarly to the proof of (a), Substep 1b) commuting symplectic linear operators $\tilde L_{as}$, $a\in\M_s$, and $\tilde B_{\ell s}$, $\ell\in\I_s\cup\H_s$, on $\tilde V_s$, and obtain a unique decomposition
$$
\tilde V_s=\bigoplus\limits_{t}\tilde V_{st}
$$
such that each $\tilde V_{st}$ is close to $V_{st}$ and is invariant under each $\tilde A_\ell$ and each $\tilde M_a$. Clearly, the subspaces $\tilde V_{st}$ are symplectic and pairwise skew-orthogonal.

Following the proof of (a), Substeps 1c, 1d, one can provide an algorithm for constructing a symplectic basis of each subspace $\tilde V_{st}$ close to the corresponding symplectic basis of $V_{st}$ and satisfying analogues of
(\ref {eq:A:ell})--(\ref {eq:J:focus:hyp}) w.r.t.\ $\tilde M_a,\tilde B_{\ell s},\tilde\Omega|_{\tilde m_0},\ddd^2\tilde J_\ell|_{\tilde m_0}$. The resulting symplectic basis on $\tilde V=T_{\tilde m_0}M$ brings the symplectic structure $\tilde\Omega|_{\tilde m_0}$, each operator $\tilde M_a$ and each quadratic form $\ddd^2\tilde J_\ell|_{\tilde m_0}$ to the desired diagonal forms as in 
(\ref {eq:hat:omega:})--(\ref {eq:twist:::}).

{\em Step 3.}
Similarly to the proof of (a), Step 2, one constructs ``perturbed'' coordinates $(\tilde x,\tilde y)$ on $\tilde U$ that are close to $(x,y)$ in $C^{k-4}-$norm, and such that $\tilde x_j(\tilde m_0)=\tilde y_j(\tilde m_0)=0$.
Indeed, $(\tilde x,\tilde y)=(\tilde u,\tilde v)\circ\tilde\phi_1^{-1}$ and $\|\ddd\tilde s^2-\ddd s^2\|_{C^{k-2}}
+\|(\tilde u,\tilde v)-(u,v)\|_{C^{k-3}} 
+\|\tilde\Omega_0-\Omega_0\|_{C^{k-4}}
+\|(\tilde p,\tilde q)-(p,q)\|_{C^{k-3}}
+\|\tilde\alpha-\alpha\|_{C^{k-4}}
+\|\tilde\phi_t-\phi_t\|_{C^{k-4}}=O(\varepsilon)$.
Similarly to the proof of (a), Step 3, one proves that these coordinates satisfy an analogue of (\ref {eq:quadratic::}), where the constants $c_j,c_{ab}^j$ are replaced by some real constants $\tilde c_j,\tilde c_{ab}^j$ close to $c_j,c_{ab}^j$, respectively.
Similarly to the proof of (a), Step 4, one has $\tilde c_{ab}^j=c_{ab}^j$.

Lemma \ref {lem:period:} (b) (and hence Lemma \ref {lem:period} (b)) is proved.
\qed

\section{Proof of Theorems \ref {thm:period} and \ref {thm:period:}} \label {sec:app:}

(a) Similarly to the proof of Lemmata \ref {lem:period} and \ref {lem:period:} (cf.\ \S \ref{subsec:proof:lem:period} and \S \ref{subsec:proof:lem:period:}), we will give a proof of Theorem \ref {thm:period:} (a).
If $\varkappa_h=0$ (i.e., all functions $J_\ell$ generate $2\pi$-periodic flows), all our arguments and constructions literally work both in the real-analytic and $C^\infty$ cases. This will give a proof of Theorem \ref {thm:period} (a) too.

{\em Step 1.}
Under the hypotheses of Theorem \ref {thm:period:} (a), consider the Hamiltonian action of the subtorus $(S^1)^r\subset(S^1)^{r+\varkappa_e+\varkappa_h}$ on a small neighbourhood $U$ of the orbit ${\cal O}_{m_0}$ generated by the functions $I_1,\dots,I_r$. 
An element $\psi=(\psi_1,\dots,\psi_r)\in(S^1)^r=(\R/2\pi{\mathbb Z})^r$ acts by the transformation $\phi_{I_1}^{\psi_1}\circ\dots\circ\phi_{I_r}^{\psi_r}$,
where $\phi_f^t$ denotes the Hamiltonian flow generated by the function $f$. Denote
$$
\Gamma:=\left\{ \psi \in(S^1)^r \mid \phi_{I_1}^{\psi_1}\circ\dots\circ\phi_{I_r}^{\psi_r} (m_0) = m_0 \right\},
$$
so $\Gamma$ is the isotropy subgroup of the point $m_0$ (and, hence, of each point of the orbit ${\cal O}_{m_0}$) w.r.t.\ the $(S^1)^r$-action.
It is a closed subgroup of the torus $(S^1)^r$. It is commutative, since the torus $(S^1)^r$ is commutative. It is discrete, since the $(S^1)^r$-action is locally free by assumption. Therefore $\Gamma$ is a finite commutative group.

Since the $(S^1)^r$-action is transitive on the orbit ${\cal O}_{m_0}$, this orbit is a torus 
$$
{\cal O}_{m_0} \approx (S^1)^r / \Gamma \approx \R^r / p^{-1}(\Gamma)
$$
where $p:\R^r\to(\R/2\pi{\mathbb Z})^r=(S^1)^r$ denotes the projection.

Choose a set of generators $\gamma^1,\dots,\gamma^r\in\R^r$ of the lattice $p^{-1}(\Gamma)\subset\R^r$, so they are basic cycles of the homology group $H_1({\cal O}_{m_0})$. Then $\{\psi^a=p(\gamma^a)\in\Gamma \mid 1\le a\le r\}$ is a generating set (perhaps, not minimal) of the finite Abelian group $\Gamma$.

Following the arguments of \cite[\S4]{zung:mir04}, we conclude that there exist 
a $(S^1)^{r+\varkappa_e}$-invariant tubular neighbourhood $U({\cal O}_{m_0})$ of ${\cal O}_{m_0}$ and a normal finite covering $\widehat{U({\cal O}_{m_0})}$ of $U({\cal O}_{m_0})$  such that the (locally-free) $(S^1)^r$-action on $U({\cal O}_{m_0})$ can be pulled back to a free $(S^1)^r$-action on $\widehat{U({\cal O}_{m_0})}$.
The symplectic form $\Omega$, the momentum map $F$ and its corresponding singular Lagrangian fibration, and the action functions $I_1,\dots,I_r$ can be pulled back to $\widehat{U({\cal O}_{m_0})}$. We will use \ $\widehat{}$ \ to denote the pull-back: for example, the pull-back of ${\cal O}_{m_0}$ is denoted by ${\widehat {\cal O}}_{m_0}$, and the pull-back of $I_1$ is denoted by $\widehat I_1$.
The free action of $\Gamma$ on $\widehat{U({\cal O}_{m_0})}$ (this action will be denoted by $\rho$) commutes with the free $(S^1)^r$-action. By cancelling out the translations (symplectomorphisms given by the $(S^1)^r$-action are called translations), we get another action of $\Gamma$ on $\widehat{U({\cal O}_{m_0})}$ that fixes ${\widehat {\cal O}}_{m_0}$. We will denote this latter action by $\rho'$.

Take a point $\widehat m_0\in\widehat{\cal O}_{m_0}$ (a pullback of $m_0$), a local disk $\widehat P$ of dimension $2n-r$ that intersects $\widehat{\cal O}_{m_0}$ transversally at $\widehat m_0$ and that is preserved by $\rho'$. Denote by $\check\varphi_1,\dots,\check\varphi_r$ the uniquely defined functions modulo $2\pi$ on $\widehat{U({\cal O}_{m_0})}$ that vanish on $\widehat P$ and such that $\widehat X_{I_i}(\check\varphi_j)=\delta_{ij}$, $1\le i,j\le r$. Then each local disk $\{\widehat I_1=\const,\dots,\widehat I_r=\const\}\cap\widehat P$ near $\widehat m_0$ has an induced symplectic structure, induced functions $\widehat J_1,\dots,\widehat J_{\varkappa_e+\varkappa_h}$ that pairwise Poisson commute, and an induced Hamiltonian $(S^1)^{\varkappa_e+\varkappa_h}$-action generated by $\widehat J_1,\dots,\widehat J_{\varkappa_e},i\widehat J_{\varkappa_e+1},\dots,i\widehat J_{\varkappa_e+\varkappa_h}$.

{\em Step 2.}
Applying Lemma \ref {lem:period:} (b) in the case with a fixed point, finite symmetry group $\Gamma$, and parameters $I_1,\dots,I_r$, we can define local functions $x_1,y_1,\dots,x_{n-r},y_{n-r}$ on $\widehat P$, such that they form a local symplectic coordinate system on each local disk $\{\widehat I_1=\const,\dots,\widehat I_r=\const\}\cap\widehat P$, with respect to which the induced Hamiltonian $(S^1)^{\varkappa_e+\varkappa_h}$-action is linear and the action $\rho'$ of $\Gamma$ is linear. We extend $x_1,y_1,\dots,x_{n-r},y_{n-r}$ to functions on $\widehat{U({\cal O}_{m_0})}$ by making them invariant under the action of $(S^1)^r$.

It follows from \cite[Lemma 4.2]{zung:mir04} that the symplectic structure $\widehat \Omega$ on $\widehat{U(\mathcal{O}_{m_0})}$ has the form 
$$
\widehat\Omega = \sum_{s=1}^r\ddd \widehat I_s\wedge\ddd(\check\varphi_s + g_s) + \sum_{j=1}^{n-r} \ddd x_j\wedge\ddd y_j ,
$$
for some real-analytic functions $g_s$ in a neighbourhood of $\widehat{\cal O}_{m_0}$ in $\widehat{U({\cal O}_{m_0})}$, that are invariant under the $(S^1)^r$-subaction.

Define $\varphi_s:=\check\varphi_s+g_s$. Then with respect to the coordinate system $(\widehat I_s,\varphi_s,x_j,y_j)$, the symplectic form $\widehat \Omega$ has the standard form, the Hamiltonian $(S^1)^{r+\varkappa_e+\varkappa_h}$-action is linear, and the free action $\rho$ of $\Gamma$ is also linear.

This yields Theorem \ref {thm:period:} (a), and hence Theorem \ref {thm:period} (a).
\qed

\medskip
(b) For a ``perturbed'' system, we follow the same proof as for the ``unperturbed'' system (cf.\ Steps 1 and 2 from above), with the only difference that we apply Lemma \ref {lem:period:} (b) in Step 2 to the ``perturbed'' system. We can and will assume that 
$\|\tilde F^{\mathbb C}-F^{\mathbb C}\|_{C^k} + \|\tilde\Omega^{\mathbb C}-\Omega^{\mathbb C}\|_{C^{k+n-2}} < \varepsilon$ (one has similar inequalities for the real objects on $M$ in assumptions of Theorem \ref {thm:period} (b)).

{\em Step 3.}
Let us use \ $\widetilde{}$ \ for all ``perturbed'' objects. 
Denote by $U'({\cal O}_{m_0})\subset U({\cal O}_{m_0})$ a tubular neighbourhood invariant under the ``unperturbed'' $(S^1)^{r+\varkappa_e}$-action, and by
$\widetilde U({\cal O}_{m_0})\supset U'({\cal O}_{m_0})$ a tubular neighbourhood invariant under the ``perturbed'' $(S^1)^{r+\varkappa_e}$-action.

Clearly, we can lift the ``perturbed'' $(S^1)^r$-action to a covering $\widehat{\widetilde U({\cal O}_{m_0})}\subset\widehat{U({\cal O}_{m_0})}$.
Consider the ``perturbed'' $(S^1)^{\varkappa_e+\varkappa_h}$-action and the ``perturbed'' $\Gamma$-actions $\tilde\rho$ and $\tilde\rho'$ on $\widetilde U({\cal O}_{m_0})$; they are close to the ``unperturbed'' ones in $C^{k-1}-$norm. By the same arguments as in Step 2, it follows from Lemma \ref{lem:period:} (b) that there exists 
a coordinate system $(\widehat {\widetilde I}_s,\tilde \varphi_s,\tilde x_j,\tilde y_j)$ on $\widehat{\widetilde U({\cal O}_{m_0})}$ close to $(\widehat I_s,\varphi_s,x_j,y_j)$ in $C^{k-4}-$norm (if $k\ge5$), in which the ``perturbed'' symplectic form $\widehat {\widetilde\Omega}$ has the standard form, the ``perturbed'' Hamiltonian $(S^1)^{r+\varkappa_e+\varkappa_h}$-action is linear, and the ``perturbed'' action $\tilde\rho'$ of $\Gamma$ is also linear, moreover the ``perturbed'' $\Gamma\times(S^1)^{r+\varkappa_e+\varkappa_h}$-action on the $(\tilde x,\tilde y)$-component is the same as the ``unperturbed'' $\Gamma\times(S^1)^{r+\varkappa_e+\varkappa_h}$-action on the $(x,y)$-component.
In particular, the ``perturbed'' free action $\tilde\rho$ of $\Gamma$ is also linear.

By above, the coordinate system $(\widehat {\widetilde I}_s,\tilde \varphi_s,\tilde x_j,\tilde y_j)$ is close to $(\widehat I_s,\varphi_s,x_j,y_j)$, hence the ``perturbed'' normal form is close to the ``unperturbed'' one.
Since the group $\Gamma$ is finite and the ``perturbed'' and the ``unperturbed'' $\Gamma\times(S^1)^{r+\varkappa_e+\varkappa_h}$-actions on the $\tilde\varphi$- and $\varphi$-components are close to each other, we conclude that these actions are the same.

This yields Theorem \ref {thm:period:} (b), and hence Theorem \ref {thm:period} (b).
\qed

\medskip
%\section{Declarations}
{\bf Akcnowledgements }
The author is grateful to Alexey Bolsinov for helpful comments on Cartan subalgebras of the Lie algebra $sp(2n,\R)$ and valuable suggestions on a preliminary version of the paper, to Andrey Oshemkov for useful discussion on proving extendability of homomorphisms to a circle from a finite subgroup of a torus (cf.\ (\ref {eq:extension})), to Stefan Nemirovski for helpful comments on topologies on the spaces of analytic functions, and to the referees for useful comments which helped to improve the paper.

\medskip
{\bf Funding }
This work was supported by the Russian Science Foundation (project 17-11-01303).
%\subsection{Conflicts of interest/Competing interests}
%Not applicable.
%\subsection{Availability of data and material (data transparency)}
%Not applicable.
%\subsection{Code availability}
%Not applicable.
%\subsection{Authors' contributions}% (optional: please review the submission guidelines from the journal whether statements are mandatory)
%Not applicable.

\end{document}